\documentclass[final,onefignum,onetabnum,letterpaper]{siamonline190516}
\usepackage[scale=0.8]{geometry} % Reduce document margins

\usepackage{soul}
\usepackage{lipsum}
\usepackage{amsfonts}
\usepackage{epstopdf}
\usepackage{epsfig}
\ifpdf
  \DeclareGraphicsExtensions{.eps,.pdf,.png,.jpg}
\else
  \DeclareGraphicsExtensions{.eps}
\fi

\usepackage{datetime}
\newdateformat{monthyeardate}{%
	\monthname[\THEMONTH] \THEDAY, \THEYEAR}

\usepackage{academicons}
\usepackage{xcolor}
\renewcommand{\orcid}[1]{\href{https://orcid.org/#1}{\textcolor[HTML]{A6CE39}{orcid.org/#1}}}

\usepackage{amsmath}
\allowdisplaybreaks
\usepackage{amssymb}
\usepackage{commath}
\usepackage{mathtools}
\usepackage{bbm}
\usepackage{bm}
\usepackage{upgreek}

\usepackage{color}
\usepackage{graphicx}
\usepackage[small]{caption}
\usepackage{subcaption}

\usepackage{relsize}
\usepackage{adjustbox}
\usepackage{algorithmic}
\usepackage{booktabs}
\usepackage{tikz}
\usepackage{pifont}% http://ctan.org/pkg/pifont for cmark
\usetikzlibrary{bayesnet} % for Bayesian network (graphical models)

\usepackage{enumitem}
\setlist[enumerate]{leftmargin=.5in}
\setlist[itemize]{leftmargin=.5in}

\newsiamremark{remark}{Remark}
\newsiamremark{example}{Example}
\newsiamremark{hypothesis}{Hypothesis}
\crefname{hypothesis}{Hypothesis}{Hypotheses}
\newsiamthm{claim}{Claim}
\headers{A Bayesian framework for spectral reprojection}{T.\ Li and A.\ Gelb}

\title{A Bayesian framework for spectral reprojection \footnote{This work is partially supported by the NSF grant DMS \#1912685, AFOSR grant \#F9550-22-1-0411, DOE ASCR \#DE-ACO5-000R22725, and DOD ONR MURI grant \#N00014-20-1-2595.}
\thanks{
\monthyeardate\today 
}}

\author{ 
Tongtong Li\thanks{Department of Mathematics, Dartmouth College, Hanover, NH 03755, USA (\email{tongtong.li@dartmouth.edu, tongtong.li@umbc.edu}, \orcid{0000-0002-7664-4764} and \email{Anne.E.Gelb@Dartmouth.edu}, \orcid{0000-0002-9219-4572})}
\and 
Anne Gelb\footnotemark[2]
}

\usepackage{amsopn}

\usepackage{comment}
\newcommand{\ds}{\displaystyle}

\usepackage[normalem]{ulem}

\newcommand{\beps}{{\boldsymbol\epsilon}}
\newcommand{\bmu}{{\boldsymbol\mu}}

\newcommand{\bdelta}{{\boldsymbol\delta}}
\newcommand{\f}{\mathbf{f}}
\newcommand{\g}{\mathbf{g}}
\newcommand{\bb}{\mathbf{b}}

\newcommand{\btheta}{\boldsymbol{\theta}}
\newcommand{\0}{{\mathbf{0}}}
\def\bF{\mathbf{F}}

\def\bH{\mathbf{H}}
\def\bI{\mathbf{I}}
\def\bW{\mathbf{W}}
\newcommand{\bbC}{\mathbb{C}}
\newcommand{\bbR}{\mathbb{R}}
\newcommand{\bbZ}{\mathbb{Z}}
\newcommand{\cB}{\mathcal{B}}
\newcommand{\cC}{\mathcal{C}}
\newcommand{\cE}{\mathcal{E}}
\newcommand{\cF}{\mathcal{F}}
\newcommand{\cG}{\mathcal{G}}
\newcommand{\cJ}{\mathcal{J}}
\newcommand{\cN}{\mathcal{N}}
\newcommand{\cX}{\mathcal{X}}
\newcommand{\cA}{\mathcal{A}}
\def\wt{\widetilde}
\def\wh{\widehat}

\newcommand{\hf}{\widehat{f}}
\newcommand{\hfl}{\widehat{f}^\lambda}
\newcommand{\hgl}{\widehat{g}^\lambda}

\newcommand{\cll}{C^{\lambda}_l}

\newcommand{\lm}{\lambda}
\newcommand{\AG}[1]{{\color{red}#1}}

\usepackage[normalem]{ulem}

\newtheorem{rem}{Remark}[section]

\newtheorem{defi}{Definition}[section]

\numberwithin{equation}{section}
\numberwithin{figure}{section}
\numberwithin{table}{section}

\allowdisplaybreaks

\usepackage{lineno}
\linenumbers

\ifpdf
\hypersetup{
  pdftitle={BSR},
  pdfauthor={Tongtong Li}
}
\fi

\begin{document}

\nolinenumbers

\maketitle

% REQUIRED
\begin{abstract}
	Fourier partial sum approximations yield exponential accuracy for smooth and periodic functions, but produce the infamous Gibbs phenomenon for non-periodic ones. Spectral reprojection  resolves the Gibbs phenomenon by projecting the Fourier partial sum onto a Gibbs complementary basis, often prescribed as the Gegenbauer polynomials. Noise in the Fourier data and the Runge phenomenon both degrade the quality of the Gegenbauer reconstruction solution, however.  Motivated by its theoretical convergence properties, this paper proposes a new Bayesian framework for spectral reprojection, which allows a greater understanding of the impact of noise on the reprojection method from a statistical point of view. We are also able to improve the robustness with respect to the Gegenbauer polynomials parameters.  Finally, the framework provides a mechanism to quantify the uncertainty of the solution estimate. % two Bayesian approaches to construct a posterior using a prior based on spectral reprojection. Each method has its own advantages in certain environments. Both provide the insights of connections between numerical analysis and Bayesian framework, and have the advantage of quantifying the uncertainty of the estimate.
\end{abstract}

% REQUIRED

% REQUIRED

% Once the paper is published
\begin{DOI}
	Not yet assigned
\end{DOI}

%\section*{Notes}
%\begin{itemize}
%	\item 
%\end{itemize}

\section{Introduction}
\label{sec:introduction}

Fourier samples model data acquisitions  in applications such as magnetic resonance imaging (MRI) and synthetic aperture radar (SAR). Indeed, recovering images from either MRI $k$-space data or SAR phase history data most often involves recasting the problem as a linear inverse problem with the forward operator given by the discrete (non-uniform) Fourier transform (DFT) matrix (see e.g.~\cite{ACGR,AG, gorham2010sar, jakowatz2012spotlight,PipeMenon}).  Compressive sensing (CS) algorithms that promote sparse solutions in a known sparse domain \cite{ArchibaldGelb,Cetin2005,DLP,Moses2004} have become increasingly widespread in providing point estimate image recoveries.  More recently, Bayesian inference methods have been developed to also quantify the uncertainty of the estimate.

This investigation develops a new Bayesian framework for recovering smooth but non-periodic functions from given noisy Fourier data. Uncertainty quantification can also be achieved when the hyperparameters in the posterior density function are fixed.  Importantly, however, rather than use the often employed {\em sparse prior} (or sparse penalty term in CS), here we construct a prior based on {\em spectral reprojection} \cite{GSSV92,Gottlieb97}.    

The spectral reprojection method is a {\em forward} approach designed to {\em reproject} the observable Fourier data onto a Gibbs complementary basis. It is sometimes referred to as Gegenbauer reconstruction when the Gibbs complementary basis is comprised of Gegenbauer polynomials. The reprojection eliminates the Gibbs phenomenon and restores the exponential convergence (hence the use of {\em spectral} in its name)  in the maximum norm. For self-containment, we summarize spectral reprojection in \cref{sec:spectralreprojection}.  
 
Although Gegenbauer reconstruction has been successfully used in applications where the observable Fourier data have complex additive Gaussian noise of  mean zero \cite{ACGR,vgcr}, it was also demonstrated in \cite{ArchGelb2002} that while the estimator is unbiased, its variance is spatially dependent. Nevertheless, theoretical results in the seminal work \cite{GSSV92} provide key insights that inspire us to develop a new Bayesian inference method for the corresponding {\em inverse} problem. Namely, the derivation of the error terms in the exponential convergence proof naturally motivates the choices of the likelihood and prior terms in the Bayesian method. In particular, the prior used for the construction of the posterior should be designed to favor solutions whose orthogonal polynomial partial sum expansion yields good approximations. Such an assumption is consistent for recovering (discretized) functions that are smooth but not periodic, and is arguably more appropriate than using a sparsifying operator such as first order differencing, which by design assumes that the underlying function is piecewise constant, or Tikhonov regularization, which is even more restrictive. Moreover, when coupled with the likelihood term, the common kernel condition \cite{kaipio2006statistical} is automatically satisfied so that a unique minimum for the corresponding \emph{maximum a posteriori} (MAP) estimate may be obtained.  We call this approach the {\em Bayesian spectral reprojection} (BSR) method, and its point estimate solution, which is consistent with Gegenbauer reconstruction, is the MAP estimate of the corresponding posterior density and is determined through optimization.  As already noted, we are also able to provide uncertainty quantification for fixed hyperparameters.

%Will say that we are inspired by the rigorous error analysis of the spectral reprojection error in \cite{Gottlieb97}, and we observe that its formulation can be used to consider an analogous optimization problem, which also looks like a maximum a posterior (map) estimate. \TL{The one inspired by spectral reprojection yield similar result, yet has the advantage of providing uncertainty quantification. 

We further propose a {\em generalized} Bayesian spectral reprojection (GBSR)  method, which modifies the BSR by formulating the likelihood so that the observables are not first transformed to the Gibbs complementary basis.  Rather, we directly use the observable Fourier data for this purpose.  Removing this restriction allows us to explore a larger space that still assumes that the function is well-approximated by the Gegenbauer partial sum expansion, but also seeks a good fit to the observable data.  As our numerical results demonstrate, the point estimate obtained for GSBR is more robust in low SNR environments to different parameter choices.  It is also possible to quantify the uncertainty when the posterior density function uses fixed hyperparameters. %  from the method above (by changing the likelihood to remove the restriction of limiting space of solutions), and could provide more accurate and robust result while still having the advantage of uncertainty quantification. (Although we motivate this method from the forward spectral reprojection method, we note that this method could be understand from inverse problem point of view. Furthermore, unlike the inverse methods in \cite{ShizgalJung} which is essentially least squares and requires large data, we do not have this restrictions.)

\subsection*{Paper organization}
The rest of this paper is organized as follows.  \Cref{sec:preliminaries} describes the underlying problem and summarizes the spectral reprojection method \cite{Gottlieb97}.  The problem is then discretized in \Cref{sec:MVformulation} to enable the Bayesian approach used in \Cref{sec:Bayesspectral}.  There we describe how the theoretical results proving exponential convergence for spectral reprojection  inspire the construction of a new prior, leading to the Bayesian spectral reprojection (BSR) method.  These ideas are then modified in \Cref{sec:Bayes} for the {\em generalized} Bayesian spectral reprojection GBSR) method.  Numerical examples in \Cref{sec:numerical} show the efficacy of our new methods for recovering 1D smooth functions from noisy Fourier data.  \Cref{sec:summary} provides some concluding remarks.

\section{Preliminaries}
\label{sec:preliminaries}

%This investigation adapts the {\em spectral reprojection} approach to a Bayesian framework, and in so doing we are able to improve the point estimate approximation of the non-periodic function given noisy Fourier data as well as quantify the uncertainty.  For purposes of self-containment, we review the basic properties of spectral reprojection below. \AG{I need to rewrite this.}

%\subsection{Problem set up}
%\label{sec:problemsetup}
 Let $f$ be a real-valued analytic $L^2$ integrable function on $[-1,1]$. %that is, 
%$$f: [-1,1] \rightarrow \bbR \quad \text{s.t. } \int_{-1}^{1} |f(x)|^2 dx < \infty.$$
We are given its first $N$ (for convenience we choose $N$ to be even) noisy Fourier measurements $\wh{\bb}  \in \bbC^N$ with components 
\begin{equation}
\label{eq:noisyfourier1d}
\wh{b}_k = \wh{f}_k + \wh{\epsilon}_k, \quad\quad k = -\frac{N}{2}, \cdots, \frac{N}{2} -1, 
\end{equation}
where %the $k$th Fourier coefficient is defined as
\begin{equation}
\label{eq:fourcoeff1D}
    \wh{f}_k = \frac{1}{2}\int_{-1}^1 f(x) e^{-ik\pi x} \,dx,
\end{equation}
and the components $\wh{\epsilon}_k \sim \cC\cN(0,\alpha^{-1})$ of $\wh{\bm \epsilon} \in \bbC^N$ each follows a circularly-symmetric complex normal distribution with unknown variance $\alpha^{-1}>0$. %which we define as 
%\[\sigma^2 = \frac{1}{N} \sum_{k = -\frac{N}{2}}^{\frac{N}{2}}(\epsilon_k - \mu)^2 = \frac{1}{N} \sum_{k = -\frac{N}{2}}^{\frac{N}{2}}\epsilon_k^2 \]
%We explicitly determine $\sigma^2$ for different choices of signal to noise ratio (SNR) values given by
%\begin{equation} \label{eq:SNR}
%   \text{SNR}_{dB} =  10\log_{10}{\left(\frac{\hat{||\bm f||}^2}{\hat{||\bm \epsilon||}^2}\right)} =10\log_{10}{\left(\frac{\hat{||\bm f||}^2}{N\sigma^2}\right)}. 
%   \end{equation}
%Using various SNR values allows us to evaluate the effectiveness of our method in noisy environments.    
The noiseless Fourier partial sum approximation % given \eqref{eq:fourcoeff1D},
\begin{equation} \label{eq:partial_four_sum1}
S_N f(x) = \sum_{k = -\frac{N}{2}}^{\frac{N}{2}-1} \wh{f}_ke^{ik\pi x}
\end{equation}
yields the infamous Gibbs phenomenon for non-periodic $f$ \cite{orszag},  notably manifested by the spurious oscillations at the boundaries $x = \pm 1$ as well as an overall reduced order of convergence away from the boundaries.   We note that the Gibbs phenomenon also occurs whenever $f$ is piecewise smooth with spurious oscillations forming near internal discontinuity locations.  When these locations are known in advance  the techniques introduced in this investigation can be straightforwardly adapted via linear transformation. When they are unknown,  methods such as those in  \cite{GT1,stefan2012sparsity} can be first applied to find  them. Thus to simplify our presentation we will only consider  non-periodic analytic $L^2$ integrable  functions on $[-1,1]$.  

\subsection{Resolution of the Gibbs phenomenon}
\label{sec:Gibbs}

The Gibbs phenomenon has been extensively studied for a variety of applications.  While the purpose of this investigation is neither to survey nor compare the myriad techniques developed to resolve it, we describe a few well-known methods for context.    

Spectral filters (see e.g.~\cite{Vandeven}) provide the most straightforward and cost efficient approach for mitigating the Gibbs phenomenon, but cause fine details of the solution to be ``smoothed out'' at the boundaries of non-periodic smooth functions. Another forward method is {\em spectral reprojection}, first introduced in \cite{GSSV92}, which can more sharply resolve the approximation by first reprojecting \eqref{eq:partial_four_sum1} into a more suitable basis. The inverse polynomial reconstruction method in \cite{ShizgalJung} and least squares method developed in \cite{Adcock2012} (which is designed for general samples in a separable Hilbert space) can sometimes better resolve the Gibbs phenomenon, but there is a trade-off between stability and increased points-per-wave requirement.  Convex optimization approaches using $\ell_1$ regularization, commonly referred to as compressive sensing algorithms (e.g. \cite{ArchibaldGelb}) have also been used successfully. Each of these methods requires  parameter tuning, which  becomes increasingly difficult as more noise is prevalent in the data. Moreover, they are not designed to quantify the uncertainty of the solution.  This is particularly problematic since different parameter choices can sometimes lead to wildly different approximations. Finally we mention the statistical filter method developed in \cite{Solomonoff} which constructs the approximation from the given Fourier data by minimizing the error over a set of possible reconstructions.  Although the approach is Bayesian, it does not provide any uncertainty quantification. Our goal here is therefore to develop a Bayesian approach that (1) accurately approximates a smooth but non-periodic function in $[-1,1]$ on a chosen set of pixels given a finite number of its noisy Fourier coefficients, and (2) quantifies the uncertainty in the approximation.  %That is,  we seek to recover ${\bm f} \in \mathbb{R}^N$ where each component approximates $f(x_j)$, $j = 1,\dots,N$, in $[-1,1]$, a

\subsection{Spectral reprojection}
\label{sec:spectralreprojection}

In a nutshell, spectral reprojection  {\em reprojects} $S_N f$ in \eqref{eq:partial_four_sum1} onto a more suitable {\em Gibbs complementary} basis.  Under certain parameter constraints, the Gegenbauer polynomials serve as a reprojection basis that ensures exponential accuracy in the maximum norm {\em without} significantly increasing the points-per-wave requirement. These constraints eliminate both Chebyshev and Legendre polynomials as suitable.\footnote{Other Gibbs complementary bases, such as the Freud polynomials, may be more robust \cite{GelbTanner}.   As here we are focused on using a Bayesian framework, the closed form solutions provided by the Gegenbauer polynomials makes the analysis and implementation more convenient.}     Below we provide a brief synopsis as it pertains to our investigation.   An in depth  review of spectral reprojection and its convergence analysis may be found in \cite{Gottlieb97}. 

In order to achieve exponential accuracy, the reprojection basis must have the property that the high modes of the Fourier basis have exponentially small projections onto its  low modes. Let the polynomial basis $\left \lbrace \psi^\lambda_l \right \rbrace$, with the parameter $\lambda \in \bbR^+$ corresponding to a prescribed weight function $w_\lambda$ and index $l \in \bbZ^+ \cup\{0\}$ corresponding to its degree, be such a basis. The $m$th-degree polynomial approximation of a function $f$ using $\left \lbrace \psi^\lambda_l \right \rbrace_{l = 0}^m$ is given by
\begin{equation}
P_m^\lambda f(x) := \sum_{l=0}^m \wh{f}_l^\lambda \psi_l^\lambda(x). 
\label{eq:orthogsum}
\end{equation}
To determine the coefficients $\wh{f}_l^\lambda$, consider
the standard definition for a weighted inner product with respect to the prescribed weight function $w_\lambda$
\begin{equation}
\langle f,g \rangle_{w_\lambda} = \int_{-1}^1 w_\lambda(x) f(x) g(x) dx.
\label{eq:weightedinnerproduct}
\end{equation}
We choose $\left \lbrace \psi^\lambda_l \right \rbrace_{l = 0}^m$ to be orthogonal with respect to $w_\lambda$ so that
\begin{equation}\label{eq:orthog}
\langle \psi^\lambda_l, \psi^\lambda_{l^\prime} \rangle_{w_\lambda} = 
\int_{-1}^1 w_\lambda(x) \psi^\lambda_l(x) \psi^\lambda_{l^\prime}(x)  dx = 
\begin{cases}
\left \Vert \psi^\lambda_l \right \Vert^2_{w_\lambda} :=h_l^\lambda & \mbox{if $ l = l^\prime$ },  \\
0 & \mbox{otherwise},
\end{cases}
\end{equation} 
where for ease of notation we have defined 
$h_l^\lambda = \Vert \psi^\lambda_l  \Vert^2_{w_\lambda}.$
Hence by construction the coefficients $\wh{f}_\ell^\lambda$ in \eqref{eq:orthogsum} are determined as
\begin{equation}
\label{eq:gegcoeff}
\wh{f}^\lambda_l = \frac{\langle f, \psi^\lambda_l\rangle_{w_\lambda}}{\left \Vert \psi^\lambda_l \right \Vert^2_{w_\lambda}} = \frac{1}{h_l^\lambda} \int_{-1}^1 w_\lambda(x) f(x) \psi^\lambda_{l}(x)  dx.
\end{equation}

% The approximation \eqref{eq:orthogsum} is on the space of orthogonal polynomials defined by the weighted inner product

An admissible Gibbs complementary basis must first and foremost yield exponential convergence (see e.g.~\cite{Boyd,HGG}) in the maximum norm as a partial sum expansion for analytic $L^2$ integrable $f \in L^2[-1,1]$. That is,
\begin{equation}
\label{eq:gegexperror}
\max_{-1 \leq x \leq 1} \left| f(x) - P_m^\lambda f(x) \right| \leq Const_1 \cdot e^{-m^r}
\end{equation}
for $m >> 1$ and $r > 0$.  The significance of the {\em reprojection} method lies in the observation that $f$ in \eqref{eq:orthogsum} can be replaced by its $N$-term Fourier expansion $S_Nf$ in \eqref{eq:partial_four_sum1} yielding 
\begin{equation} P_m^\lambda S_Nf(x):= \sum_{l=0}^m \wh{g}_l^\lambda \psi_l^\lambda(x) \label{eq:pmsnf} \end{equation}
{where
$$\wh{g}_l^\lambda  =  \frac{\langle S_Nf, \psi_l^\lambda\rangle_{w_\lambda}}{\left \Vert \psi_l^\lambda \right \Vert^2_{w_\lambda}} = \frac{1}{h_l^\lambda}\int_{-1}^1 w_\lambda(x) S_Nf(x) \psi^\lambda_{l}(x)  dx,$$
while still retaining exponential reconstruction accuracy in the maximum norm on $[-1,1]$. That is, 
\begin{equation}
\label{eq:gegexperror2}
\max_{-1 \leq x \leq 1} \left| f(x) - P_m^\lambda S_Nf(x) \right| \leq Const_2 \cdot e^{-m^p},
\end{equation}
for $m >> 1$ and $p > 0$. Satisfying \eqref{eq:gegexperror2} requires that $m = \kappa N$ for some admissible $\kappa < 1$ (see \Cref{thm:GSSV92} below). 
This is a profound statement which says that a poorly converging $N$-term Fourier expansion contains sufficient information to recover and approximate the function with exponential accuracy. The form of \eqref{eq:pmsnf} gives the method its name -- the Fourier partial sum approximation  is {\em reprojected} to a Gibbs complementary basis yielding {\em spectral} accuracy.

\subsection{Gegenbauer polynomials}
\label{sec:gegenbauer}
Mainly due to straightforward analysis and convenient closed form equations that can be used directly in numerical implementation, the most commonly employed Gibbs complementary basis is the family of Gegenbauer polynomials, $\{C^\lambda_l(x)\}_{l = 0}^m$.  We  use them here as well, and note that for robustness purposes, in some cases it might be beneficial to instead use the Freud polynomials \cite{GelbTanner} (see  \cref{rem:freudpolys}). We also point out that there is no inherent constraint that the reprojection basis need be polynomials.
%The Gegenbauer polynomials are defined by 
We now proceed to provide the definition of the Gegenbauer polynomials. As in \cite{GSSV92}, we rely on \cite{bateman} for several of the corresponding explicit formulae.
%orthogonal over the range $x\in [-1,1]$  as defined by %with the weight function $(1-x^2)^{\lambda-\frac{1}{2}}$.
\begin{defi}
\label{def:gegpolys}
The Gegenbauer polynomial  $C_l^{\lambda}(x)$,  $\lambda \geq 0$, is the polynomial of degree $l$ that satisfies
$$ \int_{-1}^{1} w_\lambda(x) C_k^{\lambda}(x) C_l^{\lambda}(x) dx=0, \qquad k\neq l,$$
where $w_\lambda(x) = (1-x^2)^{\lambda-\frac{1}{2}}$ is the corresponding weight function.  The normalizing term $h_l^\lambda$ in \cref{eq:gegcoeff} %is then
%$$\AG{h_l^\lambda}  =\int_{-1}^{1} w_\lambda(x) \left( C_l^{\lambda}(x) \right)^2 dx,$$
%which 
is explicitly given by 
\begin{equation}\label{eq: h}
h_l^\lambda=\pi^{\frac{1}{2}}C_l^{\lambda}(1)\frac{\Gamma(\lambda+\frac{1}{2})}{\Gamma(\lambda)\Gamma(l+\lambda)}.
\end{equation}
Here
 $C_l^{\lambda}(1)=\frac{\Gamma(l+2\lambda)}{l!\Gamma(2\lambda)}$
and  $\Gamma (\lm)$ denotes the usual gamma function.
\end{defi}

By substituting $C^\lambda_l(x)$ into \eqref{eq:orthogsum}, we obtain the  Gegenbauer partial sum expansion 
\begin{equation}
\label{eq:gegsum}
P_m^{\lambda}f(x) = \sum_{l=0}^{m} \wh{f}^{\lambda}_{l}\, C^\lambda_l(x),
\end{equation}
with Gegenbauer coefficients in \eqref{eq:gegcoeff} provided by
\begin{equation}\label{eq: Gegen0}
\wh{f}^{\lambda}_{l}=\frac{1}{h^\lambda_l} \int_{-1}^{1}(1-x^2)^{\lambda-\frac{1}{2}}\,f(x)\,C^\lambda_l(x)dx.  
\end{equation}
Note that \eqref{eq:gegsum} converges exponentially in the sense of \eqref{eq:gegexperror} for $f$ analytic in $[-1,1]$ and 
 $\lambda \ge 0$ \cite{Boyd,HGG}.
In particular, Chebyshev ($\lambda = 0$) and Legendre ($\lambda = \frac{1}{2}$) are classical examples of Gegenbauer polynomials. %, a

Because it inspires the development of our new Bayesian approach to spectral reprojection, we quote the main theorem in \cite[Theorem 5.1]{GSSV92}: %for when the Gibbs complementary basis is given by the Gegenbauer polynomials:
\begin{theorem} [Removal of the Gibbs Phenomenon]
\label{thm:GSSV92}
Consider an analytic and non-periodic function $f(x)$ on $[-1,1]$, satisfying
\begin{equation}
\label{eq:analrequirement}
\max_{-1 \leq x \leq 1} \left| \frac{d^k f}{dx^k} (x) \right|
\leq C(\rho ) \frac{k!}{\rho^k}, \qquad \rho \geq 1.
\end{equation}
Assume that the first $N$ Fourier coefficients in \eqref{eq:fourcoeff1D} are known.
Let $\hgl_l$, $0 \leq l \leq m$, be the Gegenbauer reprojection coefficients in \eqref{eq:pmsnf} %the Gegenbauer expansion coefficients of $f_N (x) = \skn \hf_k \eikpx$ explicitly given by
\begin{equation}
\hgl_l =  \frac{1}{h_l^\lambda} \int_{-1}^{1}(1-x^2)^{\lambda-\frac{1}{2}}\,S_Nf(x)\,C^\lambda_l(x)dx
= \delta_{0l} \hf_0 + \Gamma (\lm ) i^l (l+\lm )
\sum_{k = -\frac{N}{2}}^{\frac{N}{2}-1} J_{l+\lm}(\pi k) \left( \frac{2}{\pi k} \right)^\lm
\hf_k,
\label{eq:gegcoeffs}
\end{equation}
where $S_Nf$ is the Fourier partial sum given by \eqref{eq:partial_four_sum1}, $\delta_{0l}$ is the Kronecker delta function, $\Gamma (\lm)$ denotes the usual gamma function, and $J_{l+\lm}(\pi k)$ is the Bessel function of the first kind.\footnote{The explicit formula in \eqref{eq:gegcoeffs} is discussed in \cite{GSSV92} and results from properties of orthogonal polynomials, see e.g. \cite{abramowitzstegun}.} 
Then, for $\lm = m =\kappa N$ where $\kappa < \frac{\pi e}{27}$, we have
\begin{equation}
\label{eq:gegbound}
\max_{-1 \leq x \leq 1} \left| f(x) - \sum_{l=0}^{m} \hgl_l \cll (x) \right| \leq A_D N^2 q_D^N +  A_R q_R^N,
\end{equation}
where
$$
q_D = \left( \frac{27 \kappa}{\pi e} \right)^{\kappa},
 \qquad
q_R = \left( \frac{27 }{32 \rho} \right)^{\kappa}.
$$
\end{theorem}
The proof found in \cite{GSSV92,Gottlieb97} makes use of the explicit form in \eqref{eq:gegcoeffs} and follows the triangular inequality 
\begin{eqnarray}
\max_{-1 \le x \le 1} \left| f(x) - \sum_{l=0}^{m} \hgl_l \cll (x) \right| &  \leq & \max_{-1 \le x \le 1} \left|  \sum_{l=0}^{m} \hfl_l \cll (x)  -\sum_{l=0}^{m} \hgl_l \cll (x) \right|\nonumber\\
&+&\max_{-1 \le x \le 1} \left| f(x) - \sum_{l=0}^{m} \hfl_l \cll (x) \right| \label{eq:triangle}\\
%  \max_{-1 \leq x \leq 1} \left|  \sum_{l=0}^{m} \hfl_l \cll (x)  -\sum_{l=0}^{m} \hgl_l \cll (x) \right|\\
&=& ||E_{D}||_\infty + ||E_{R}||_\infty\nonumber\\
&\le& A_D N^2q_D^N + A_Rq_R^N.\nonumber
\end{eqnarray}
The first term, called the {\em truncation} error in \cite{GSSV92,Gottlieb97}, is a measure of how well the projected data fits the true data (hence we refer to this term as the {\em data} error and use $E_D$) in the Gibbs complementary basis (here the Gegenbauer polynomials with $\lambda = \kappa N$). The regularization error ($E_R$) measures how accurate the $m$th-degree approximation is in the Gibbs complementary basis given the true coefficients $\{\hfl_l\}_{l = 0}^m$.  Importantly, the error bound is provided in the {\em maximum} norm, proving that the spectral reprojection method in \eqref{eq:pmsnf} using the Gegenbauer polynomials as the reprojection basis with appropriate $\lambda$ and $m$  truly resolves the Gibbs phenomenon.   We also note that referring to the two error components as respectively the data and regularization errors provides some intuition for the Bayesian approach we adopt in  \Cref{sec:Bayesspectral}.

To summarize, given the first $N$ Fourier coefficients \eqref{eq:fourcoeff1D} of a non-periodic analytic function on $[-1,1]$, \Cref{thm:GSSV92} says that the Gegenbauer polynomials can be used as the Gibbs complementary basis in \eqref{eq:pmsnf}, with coefficients $\{\wh{g}_l^\lambda\}_{l = 0}^m$ defined in \eqref{eq:gegcoeffs}. This yields the celebrated {\em Gegenbauer reconstruction method} 
\begin{equation}\label{eq: fmNlambda1}
f_{m,N}^{\lambda}(x) := P_m^\lambda S_N f(x) = \sum_{l=0}^{m} \wh{g}^{\lambda}_{l}\, C^\lambda_l(x), 
\end{equation}
which converges exponentially to $f$ for any $x \in [-1,1]$ in the maximum norm.

\begin{rem}
\Cref{thm:GSSV92} does not suggest that $m = \lambda$, only that they both grow linearly with $N$.  Indeed, while the regularization error $E_R$ decreases with larger $m$, the data error $E_D$ increases.  Larger $\lambda$ decreases $E_D$, but creates a narrower weight function, causing the approximation to be more extrapolatory.  A more thorough investigation for optimizing these  parameters, based on the underlying regularity of $f$, can be found in \cite{GelbJack2006}.  Moreover, the theorem rules out using Legendre ($\lambda = \frac{1}{2}$) or Chebyshev polynomials ($\lambda = 0$), since $\lambda$ must grow linearly with $N$ to ensure small  $E_D$.  Thus we see that although they both yield exponentially convergent approximations for $E_R$ according to  \eqref{eq:gegexperror}, they do {\em not} qualify as admissible Gibbs complementary bases.  Finally, it is important to note that  \Cref{thm:GSSV92} does not account for the additive Gaussian noise in \eqref{eq:noisyfourier1d}, and while the general estimate still holds \cite{ArchGelb2002}, parameter optimization is less clear.  %Since this is not the focus of our investigation, in our numerical experiments we select $\lambda$ and $m$ to satisfy  \Cref{thm:GSSV92} and do not attempt to optimize them further.  
This will be further discussed  in \cref{rem:prior1}, \cref{rem:prior2}, and \cref{sec:SNRcase}.
\end{rem}

Although \eqref{eq: fmNlambda1} can approximate $f$ at any point location $x$ in $[-1,1]$, for simplicity in this investigation we discretize the spatial domain uniformly using the same number of points as given Fourier data, so that 
\begin{equation}
x_j = -1 + \frac{2(j-1)}{N}, \qquad j= 1, ..., N.
\label{eq:uniformgrid}
\end{equation}
Uniform discretization is also desirable as it is standard for pixelated images. We then seek to recover an approximation to the vectorized solution $\f \in \bbR^N$, where each component of $\f$ is given by $f_j = f(x_j)$, $j = 1,\dots, N$.
\Cref{alg:gegapprox} describes how to approximate $f$ at the prescribed gridpoints \eqref{eq:uniformgrid} satisfying the conditions in \eqref{eq:analrequirement} given its noisy Fourier coefficients \eqref{eq:noisyfourier1d} using the blueprint provided by  \Cref{thm:GSSV92}. Results for the {\em noiseless} case first described in \cite{GSSV92} are displayed in \cref{fig: 2}.

\begin{algorithm}[h!]
\caption{Gegenbauer reprojection from noisy Fourier data}% for a Real-Valued Signal $\bm x$ and sparsifying operator $L$}
\label{alg:gegapprox}
\begin{algorithmic}
\item[] {\bf Input} Data vector $\wh{\bb} = \{ \wh{b}_k \}_{k=-\frac{N}{2}}^{\frac{N}{2}-1} \in \bbC^N$ in \eqref{eq:noisyfourier1d}, %for known noise variance $\sigma^2$,
choices for $\lambda = \kappa_1 N$ and $m = \kappa_2 N$, $0 < \kappa_1, \kappa_2 < 1$, and grid points $x_j$, $j = 1,\dots,N$.
\item[] {\bf Output} Approximation $f_{m,N}^\lambda(x_j)$, $j = 1,\dots, {N}$.
\item[Step 1.] Calculate the approximate Gegenbauer coefficients as %\TL{(Do we need to distinguish the notation as this is obtained by noisy data $\wh{b}_k$ while \eqref{eq:gegcoeffs} is obtained by $\wh{f}_k$?)}
\begin{equation}\label{eq: approx Gegen Bessel} 
\wh{g}^{\lambda}_{l}=\delta_{0l}\wh{b}_0 + \Gamma(\lambda) i^{l}(l+\lambda) \sum_{k = -\frac{N}{2}}^{\frac{N}{2}-1} J_{l+\lambda}(\pi k) \left(\frac{2}{\pi k}\right)^{\lambda} \wh{b}_k,\quad\quad l = 0,\dots, m,  
\end{equation}
where $J_{\nu}(\cdot)$ is the Bessel function of the first kind. 
\item[Step 2.] Construct the Gegenbauer partial sum approximation
\begin{equation}\label{eq: fmNlambda}
f_{m,N}^{\lambda}(x_j) = \sum_{l=0}^{m} \wh{g}^{\lambda}_{l}\, C^\lambda_l(x_j), \quad\quad j = 1,\dots, N.
\end{equation}
\end{algorithmic}
\end{algorithm}

\section{Model discretization}
\label{sec:MVformulation}
Since we will be formulating our problem as a discrete linear inverse problem, we now describe how the results in \Cref{sec:preliminaries} can be discretized and expressed using matrix-vector notation. Following standard terminology for recovery problems, we will refer to ${\f} = \{f(x_j)\}_{j= 1}^N \in \mathbb{R}^N$ as the underlying signal of interest.

Recall that for uniform grid points \eqref{eq:uniformgrid}, the composite trapezoidal rule can be used to approximate $\wh{\f} = \{\wh{f}_k\}_{k = -\frac{N}{2}}^{\frac{N}{2}-1}$  in \eqref{eq:fourcoeff1D} as
\begin{equation}\label{eq: trapezoid}
\wt{f}_k = \frac{1}{N}\sum_{j = 1}^N f(x_j) e^{-ik\pi x_j}, \quad -\frac{N}{2} \le k \le  \frac{N}{2} -1,
\end{equation}
which we write in matrix-vector form as
\begin{equation}
    \label{eq:DFT}
\wt{\f} = \bF \f.
\end{equation}
Here $\wt{\f} = \{\wt{f}_k\}_{k = -\frac{N}{2}}^{\frac{N}{2}-1}\in \bbC^N$, and $\bF \in \bbC^{N \times N}$ is the discrete Fourier transform (DFT) matrix with entries 
$$\bF(n,j) = \frac{1}{N} e^{-i(n - \frac{N}{2} -1) \pi x_j},   \quad 1\leq n,j \leq N.$$
%Although \eqref{eq:DFT} is not equivalent to \eqref{eq:fourcoeff1D}, the resulting aliasing error is typically smaller than the error coming from the noise in \eqref{eq:noisyfourier1d}. Hence the model discrepancy is reasonable (see e.g.~\cite{ArchibaldGelb}).

Analogously, given Fourier coefficients $\wh{\f}= \{\wh{f}_k\}_{k = -\frac{N}{2}}^{\frac{N}{2} -1} \in \bbC^N$ in \eqref{eq:fourcoeff1D}, we write the corresponding discrete Fourier partial sum for \eqref{eq:partial_four_sum1} as 
\begin{equation}\label{eq:partial_four_sum1 mv}
\f_N = \wh{\bF}^{\ast} \,\wh{\f},
\end{equation}
where $\wh{\bF}^{\ast} \in \bbC^{N \times N}$ is the (un-normalized) complex conjugate transpose of $\bF$ in \eqref{eq:DFT}, also referred to as the discrete inverse Fourier transform operator with entries\footnote{The basis $\{e^{i k \pi x}\}_{k \in \bbZ}$ is orthogonal but not normalized, and there are different conventions for defining the discrete Fourier and inverse Fourier transform to account for the normalizing term (see e.g.~\cite{HGG}).  With a slight abuse of notation, we  adopt the definitions of $\bF$ and $\wh{\bF}^\ast$ in \eqref{eq:DFT} and \eqref{eq:partial_four_sum1 mv} which are consistent with our numerical implementation.}
$$\wh{\bF}^{\ast}(j,n) = e^{i (n-\frac{N}{2}-1) \pi x_j}, \quad\quad 1\leq j,n \leq N.$$ 
Observe that  $\f_N = \{S_N f(x_j)\}_{j = 1}^N  \in \bbR^N$ due to cancellation of the imaginary terms.   
\begin{rem}
\label{rem:aliasing}
Due to aliasing,  \eqref{eq:DFT} is not equivalent to \eqref{eq:fourcoeff1D}, with the explicit relationship given by (see e.g.~\cite{HGG})
$$\wt{f}_k = \wh{f}_k + \sum_{\substack{|q| \le \infty\\ q \ne 0}}\hat{f}_{k + \frac{Nq}{2}}, \quad -\frac{N}{2} \le k \le \frac{N}{2} -1,$$
which we will write simply as
\begin{equation}
\label{eq: Fourier discrete}
\wh{\f} = \wt{\f} + \bdelta_{f} =  \bF \f +  \bdelta_{f},
\end{equation}
with $\bdelta_{f} \in \bbC^N$ representing the aliasing error. It is typically assumed that the magnitude components of $\bdelta_{f}$ are smaller than those in $\wh{\bm \epsilon} = \{\epsilon_k\}_{k = -\frac{N}{2}}^{\frac{N}{2}-1}$ defined in \eqref{eq:noisyfourier1d}, and therefore the model discrepancy is reasonable (see e.g.~\cite{ArchibaldGelb}). Indeed, the traditional observation data model corresponding to \eqref{eq:noisyfourier1d} follows \eqref{eq: Fourier discrete} and is defined as
\begin{equation}
    \label{eq:observationFourier}
    \wh{\bb} = \bF \f + \wt{\beps},
\end{equation}
where by construction $\wt{\beps} = \wh{\beps} + \bdelta_{f} \in \mathbb{C}^N$ follows a circularly-symmetric complex normal distribution.  We will use \cref{eq:observationFourier} in \Cref{sec:Bayes}  and further characterize the noise distribution there.
\end{rem}

\begin{rem}
For simplicity we also use the same notation in \cref{eq:partial_four_sum1 mv} for the discretized Fourier partial sum  given noisy Fourier coefficients $\wh{\bb}= \{\wh{b}_k\}_{k = -\frac{N}{2}}^{\frac{N}{2} -1} \in \bbC^N$ in \eqref{eq:noisyfourier1d}. That is, we write
\begin{equation}\label{eq:partial_four_sum1 mv2}
\f_N = \wh{\bF}^{\ast} \,\wh{\bb}.
\end{equation}
Which of \cref{eq:partial_four_sum1 mv} or \cref{eq:partial_four_sum1 mv2} is employed will be contextually evident.

\end{rem}

The corresponding discrete Gegenbauer partial sum  \eqref{eq:gegsum}  can be written as
\begin{equation}\label{eq:gegsum mv}
\f_m^\lambda = \bF^{Geg} \,\wh{\f}^{\lambda},
\end{equation}
where $\f_{m}^{\lambda}=\{P_{m}^{\lambda} f(x_j)\}_{j = 1}^{N} \in \bbR^N$ and $\wh{\f}^{\lambda}= \{\wh{f}^{\lambda}_{l}\}_{l = 0}^m \in \bbR^{m+1}$. Here $\bF^{Geg} \in \bbR^{N \times (m+1)}$ with entries
\begin{equation}\label{F^Geg}
\bF^{Geg}(j,n) = C_{n-1}^{\lambda}(x_j), \quad\quad 1\leq j\leq N, \quad\quad 1 \leq n \leq m+1.
\end{equation}
Since the Gegenbauer coefficients $\wh{\f}^{\lambda}$ in \eqref{eq: Gegen0} are not explicitly known for our problem,  we use the composite trapezoidal rule to obtain the approximation 
\begin{equation}\label{eq:quadrature}
\wt{f}^{\lambda}_{l} = \frac{1}{h^\lambda_l} \cdot \frac{2}{N} \sum_{j=1}^{N} (1-x_j^2)^{\lambda-\frac{1}{2}} f(x_j) C^\lambda_l(x_j),\quad\quad l = 0,\dots,m,
\end{equation}
as the components of $\wt{\f}^{\lambda} \in \bbR^{m+1}$, noting that the integrand evaluation vanishes at both boundaries for $\lambda > \frac{1}{2}$ and that the smoothness of the integrand increases with $\lambda$. Combining \eqref{F^Geg} and \eqref{eq:quadrature} yields
\begin{equation}\label{eq: Gegen1 mv}
\wt{\f}^{\lambda} = \bF^{CT} \,\f,\quad \ds \bF^{CT}=\frac{2}{N}\bH \left(\bF^{Geg} \right)^T \bW.
\end{equation}
Here $\bH \in  \bbR^{(m+1)\times(m+1)}$ is the diagonal normalizing matrix  with entries $\bH(l,l) =\frac{1}{h^\lambda_l}$ and $\bW \in \bbR^{N \times N}$ is the  diagonal weighting matrix with entries $\bW(j,j) = (1-x_j^2)^{\lambda-\frac{1}{2}}$. 

%Moreover, the error incurred using \eqref{eq:quadrature} (assuming large enough $N$) is typically smaller than the overall error in the approximation due to the additive noise in \eqref{eq:noisyfourier1d}. 
%\TL{(TL: we want to be careful about ``large $N$" which is not our case here. One advantage of our method (induced from spectral reprojection) is we don't need very large $N$, unlike the inverse methods.)} \AG{[AG: I understand your concern but I think it's fine, the $N$ just has to do with the resolution of the function, which has to be large enough.]} 
Related to the discussion in \Cref{rem:aliasing}, there is a quadrature error associated with the Gegenbauer coefficients which we simply write as %Similar as \eqref{eq: Fourier discrete}, we have 
\begin{equation}
\label{eq: Gegen discrete}
\wh{\f}^{\lambda} = \wt{\f}^{\lambda} + \bdelta_{g} =  \bF^{CT} \, \f +  \bdelta_{g},
\end{equation}
so that the corresponding model for \eqref{eq:gegsum mv} is
\begin{equation}\label{eq:gegsum mv2}
%\f_m^\lambda \approx \bF^{Geg} \,\wt{\f}^{\lambda} = \bF^{Geg}\bF^{CT}\f.
\f_m^\lambda = \bF^{Geg}\wt{\f}^\lambda + \bF^{Geg}\bdelta_g = \bF^{Geg}\bF^{CT}\f + \wt{\bdelta}_g,
\end{equation}
where $\wt{\bdelta}_g = \bF^{Geg}\bdelta_g$ can be considered as noise resulting from not having the exact Gegenbauer coefficients $\wh{\f}^\lambda$. Moreover, assuming large enough $N$, the error $\bdelta_g$ incurred using \eqref{eq:quadrature} is typically smaller than the overall error in the approximation due to the additive noise in \eqref{eq:noisyfourier1d}.

Finally, we provide the discretized model for spectral reprojection in \eqref{eq: fmNlambda}.  This is accomplished by  first defining $\wh{\g}^{\lambda}= \{\wh{g}^{\lambda}_{l}\}_{l = 0}^m$, with each component $\wh{g}^\lambda_l$  given by \eqref{eq: approx Gegen Bessel}, which we write as
\begin{equation}\label{eq: approx Gegen Bessel mv}
\wh{\g}^{\lambda} = \bF^{Bessel} \,\wh{\bb}.
\end{equation}
Here $\bF^{Bessel} \in \bbC^{(m+1) \times N}$ has entries  
\begin{equation}
\bF^{Bessel}(j,n) = 
\begin{cases}
1 & \mbox{if $ j = 1$ and $ n=\frac{N}{2}+1$ },  \\[1ex]
0 & \mbox{if $ j > 1$ and $ n=\frac{N}{2}+1$ }, \\
\ds \Gamma(\lambda) i^{j-1}(j-1+\lambda) J_{j-1+\lambda}(\pi k) \left(\frac{2}{\pi k}\right)^{\lambda} & \mbox{otherwise},
\end{cases}
\end{equation}
for  $1\leq j \leq m+1$, $1 \leq n \leq N$ and $k = n - \frac{N}{2} - 1$.

By defining $\f_{m,N}^{\lambda}= \{f_{m,N}^{\lambda}(x_j)\}_{j = 1}^N$, we have the discretized formulation of  \eqref{eq: fmNlambda}:
\begin{equation}\label{eq: fmNlambda mv}
\f_{m,N}^\lambda = \bF^{Geg} \,\wh{\g}^{\lambda} = \bF^{Geg} \bF^{Bessel}\wh{\bb}.
\end{equation}
Note that \eqref{eq: fmNlambda mv} is {\em directly} obtained from the given measurements $\wh{\bb}$ in \eqref{eq:noisyfourier1d}, that is, there is no additional model error (noise) incurred by implementation.  Therefore, when comparing \eqref{eq: fmNlambda mv} to \cref{eq:gegsum mv2}, we obtain the  observation model
\begin{equation}
    \f_{m,N}^\lambda = \f_m^\lambda + {\bm \epsilon}^\lambda
    = \bF^{Geg}\bF^{CT}\f + {\bm \epsilon}^\lambda, \label{eq:gegobservationmodel}
\end{equation}
where the additive noise ${\bm \epsilon}^\lambda$ has contributions both from $\wh{\bm \epsilon}$ in \eqref{eq:noisyfourier1d} and $\wt{\bm \delta}_g$ in \eqref{eq:gegsum mv2}.  
\section{A Bayesian framework for spectral reprojection}
\label{sec:Bayesspectral}

Spectral reprojection has been successfully used to recover non-periodic (or piecewise smooth) functions from Fourier data in a variety of contexts, including post-processing  solutions for numerical conservation laws \cite{HGG} and recovering images for segmentation purposes \cite{ACGR}. Although some noise is always present, either in the form of numerical simulation error or data acquisition error, most studies have not formally analyzed its contribution to the overall error.  However, in \cite{ArchGelb2002} it was shown that while the Gegenbauer approximation is an unbiased estimator for additive Gaussian noise, the variance depends on the particular approximation location, with  increased variance near the boundaries (or edge locations of piecewise smooth functions). Filtering can reduce the overall variance, but requires further tuning, and moreover leads to a loss of resolution. These issues motivate us to gain a better understanding of the spectral reprojection method using a statistical framework, as we now describe.

We first note that the proof of Theorem \ref{thm:GSSV92}, and specifically the triangular inequality \eqref{eq:triangle}, provides a blueprint for formulating a Bayesian approach. In particular the inequality also holds for $\Vert \cdot \Vert_2$ up to a constant.\footnote{Using $\Vert \cdot \Vert_2$ is convenient because it provides easier analysis and implementation for the corresponding inverse problem. The Bayesian approach can also be used for other norms.} Second, the data error term $E_D$ in \eqref{eq:triangle} is analogous to the likelihood term in a \emph{maximum a posterior} (MAP) estimate for a posterior distribution of the unknown signal ${\f}$ in the discretized framework. 
%\AG{\sout{Indeed, the  spectral reprojection method, which reprojects the non-converging Fourier partial sum onto a Gibbs complementary basis, was specifically developed because the more natural recovery basis (here the partial Fourier sum) yields poor results.}} 
Finally, the regularization error term $E_{R}$ in \eqref{eq:triangle} measures the difference between the function and its Gegenbauer partial sum approximation \eqref{eq:gegsum}, which can be viewed as  our prior belief about  $\f$. That is, we expect the signal $\f$ to be close to $\f_m^\lambda$ in \eqref{eq:gegsum mv2}. In a nutshell, these insights allow us to view spectral reprojection from a Bayesian perspective for fixed values of $m$ and $\lambda$. We will refer to this approach as the {\em Bayesian Spectral Reprojection (BSR)} method.  

\begin{remark}
\label{rem:bayesspectral} As our numerical experiments will show in \Cref{sec:numerical}, the Bayesian spectral reprojection method yields estimates that are similar to those determined from \Cref{alg:gegapprox}, which is to be expected based on its design. An advantage to our new Bayesian formulation is that unlike other traditionally chosen priors (e.g.~sparse priors), its derivation is motivated by rigorous error approximation analysis. Moreover, as we will describe below, it is also possible to partially quantify the uncertainty of the underlying approximation.
\end{remark}

\subsection{Statistical formulation for spectral reprojection}\label{sec:inverseGeg}
%\AG{Rather than approaching the recovery of $f$ directly from the given data \cref{eq:noisyfourier1d},} 

Inspired by the spectral reprojection method as discretized by \cref{eq:gegsum mv2} and \eqref{eq: fmNlambda mv}, we now seek a statistical formulation for the recovery of the unknown signal $\f$.  We begin by considering the linear inverse problem 
\begin{equation}
{\cG}={\cA}\cX+{\cE}, 
\label{eq:genbayesmodel}
\end{equation}
where $\cX$, ${\cG}$, and ${\cE}$ are random variables defined over a common probability space, and $\cX$ and ${\cE}$ are assumed to be independent. In this framework $\cX$ represents the unknown which we seek to recover, ${\cG}$ denotes the observables, ${\cA}$ is a known linear operator, and ${\cE}\sim \cN(\0,\gamma^{-1} \bI_N)$ is a vector of independent and identically distributed (i.i.d.) Gaussian noise with mean $0$ and unknown variance $\gamma^{-1}$. Finally throughout this manuscript we use $\bI_N$ to denote the $N\times N$ identity matrix. 

\begin{rem}
\label{rem:gaussiannoise}  %\AG{We will have to say something here that $\cE$ is not really Gaussian, and its distribution was studied (although not complete) in \cite{ArchGelb2002}.  However we are just using it as a model.}
%\TL{(TL: I would like to point out that the noise studied in \cite{ArchGelb2002} is in the Gegenbauer coefficient space, not physical space. However, since we only have an extra linear operator $\bF^{Geg}$, the conclusion is informational.)}
As already noted, the variance of the Gegenbauer approximation has non-trivial expression and depends on the particular approximation location \cite{ArchGelb2002}. The assumption about the covariance matrix of $\cE$ being a constant multiple of the identity matrix  therefore generally does not hold. We nevertheless adopt it for simplicity in \cref{eq:genbayesmodel}.
\end{rem}
With \eqref{eq:gegobservationmodel} representing the observational model system, the realizations of each random variable in \eqref{eq:genbayesmodel} respectively are the unknown vector of interest $\f \in \bbR^N$, the observable data $\f_{m,N}^\lambda \in \bbR^N$ in \eqref{eq: fmNlambda mv}, and noise vector $\beps^\lambda \in \bbR^N$. A hierarchical Bayesian model based on Bayes' theorem then provides an estimate of the full \emph{posterior} distribution given by
\begin{equation}\label{eqn:bayes}
p_{\cX, {\Theta} \,| \,{\cG}}(\f,{\btheta} \,|\,\f_{m,N}^\lambda) \propto p_{{\cG}\,|\,\cX,{\Theta}}(\f_{m,N}^\lambda \, | \,\f,{\btheta}) \, p_{\cX \,| \,{\Theta}}(\f \,| \,{\btheta}) \, p_{{\Theta}}({\btheta}),
\end{equation}
where $p_{{\cG}\,|\,\cX,{\Theta}}(\f_{m,N}^\lambda \, | \,\f,{\btheta})$ is the likelihood density function determined by relationship between $\f$ and $\f_{m,N}^\lambda$ and assumptions on ${\cE}$,  $p_{\cX \,| \,{\Theta}}(\f \,| \,{\btheta})$ is the density of the prior distribution encoding \emph{a priori} assumptions about the solution, and $p_{{\Theta}}({\btheta})$ is the hyperprior density for all other involved parameters, which we denote collectively as ${\btheta}$. In our investigation $p_{{\Theta}}({\btheta}) = p_\upgamma(\gamma) \,p_\upbeta(\beta)$, where $\gamma$, $\beta$ corresponds to the inverse noise and inverse prior variances respectively (see Section \ref{sec:hyperprior}). Note that the right hand side of \eqref{eqn:bayes} is not normalized (hence $\propto$ is used), but as the normalized version is often not computationally tractable, it is the posterior typically sought in Bayesian inference methods.  In what follows we describe how to construct the corresponding likelihood, prior and hyperprior density functions. 

%%%%%%%%%%%%%%%%%%%%%%%%%%%%%%

\subsection{Formulating the likelihood} \label{sec:likelihood1}
 
The likelihood density function $p(\f_{m,N}^\lambda \,| \,\f, \gamma)$ in \eqref{eqn:bayes} models the connection between the signal and observation vectors.\footnote{To avoid cumbersome notation, henceforth we adopt simplified} expressions for each density function in \eqref{eqn:bayes}. As already noted, the data error term $E_D$ in \eqref{eq:triangle} states intuitively that $\f_{m,N}^{\lambda}$ is close to $\f_m^\lambda$, while $\f_m^\lambda$ is related to the underlying signal $\f$ via \eqref{eq:gegsum mv2}. This insight allows us to define the observation model as \eqref{eq:gegobservationmodel}, and more specifically to prescribe the distribution for the likelihood as
\begin{equation}
\label{eq:likelihooddistribution1}
\f_{m,N}^{\lambda} \sim \cN (\bF^{Geg}\bF^{CT} \f , \gamma^{-1}\bI),
\end{equation}
with the interpretation that the data error ${{\beps}^\lambda} =  \f_{m,N}^{\lambda} - \bF^{Geg}\bF^{CT} \f $ corresponding to $E_D$ follows the distribution $\cN(\0,\gamma^{-1} \bI_N)$ for unknown variance $\gamma^{-1}$. Thus
\begin{equation} \label{eq:EDnoise}
p({\beps}^\lambda \,| \, \gamma)= \left(\frac{\gamma}{2\pi}\right)^{\frac{N}{2}}\exp \left(-\frac{\gamma}{2} \Vert {\beps}^\lambda \Vert _2^2 \right).
\end{equation}
The likelihood density function is determined entirely from \eqref{eq:EDnoise}, yielding
\begin{equation}\label{eq:likelihood}
 p(\f_{m,N}^\lambda\,|\,\f,\gamma)   =
\left(\frac{\gamma}{2\pi}\right)^{\frac{N}{2}}\exp(-\frac{\gamma}{2} \Vert \f_{m,N}^\lambda - \bF^{Geg}\bF^{CT}\, \f \Vert_2^2).
\end{equation}
Importantly, we observe that $\f_{m,N}^{\lambda} = \bF^{Geg} \bF^{Bessel}\wh{\bb}$ is obtained {\em directly} by linear transformation of the measurement data $\wh{\bb}$ via \eqref{eq: fmNlambda mv}. %\AG{Hence the likelihood \eqref{eq:likelihood} corresponds to a simple linear transformation of \eqref{eq:genbayesmodel} for which the noise can be reasonably modeled as Gaussian with unknown variance.} \TL{(TL: Here I think \eqref{eq:genbayesmodel} directly gives \eqref{eq:likelihood}, there is no linear transformation.)}

%%%%%%%%%%%%%%%%%%%%%%%%%%%%%%

\subsection{Formulating the prior}
\label{sec:prior1}
The prior function models our prior belief about the unknown signal $\f$.  In this investigation $\f = \{f(x_j)\}_{j = 1}^N$ where $f$ is analytic in $[-1,1]$, and hence, due to \eqref{eq:gegexperror}, a distribution for the true solution $\f$ can be motivated using \eqref{eq:gegsum mv2}, resulting in  
\begin{equation}\label{eq:prior1}
\cX \sim \cN ({\f}^{\lambda}_{m} , \beta^{-1} \bI_{N} ),    
\end{equation}
%equivalently,
%\begin{equation}\label{eq: prior2}
%\f \sim \cN (\bF_{Geg}\bF_{CT} \, \f , \beta^{-1} \bI_{n} ),    
%\end{equation}
where ${\f}^{\lambda}_{m}$ is the  Gegenbauer partial sum approximation in \eqref{eq:gegsum mv2} with a hyperprior $\beta>0$ being inverse variance.
Essentially, since the regularization error term $E_{R}$ measures the difference between the function and its approximation with Gegenbauer basis, we expect the signal $\f$ to be close to  $\f^{\lambda}_{m}$, leading to the conditionally Gaussian prior function
\begin{equation}\label{eq: prior}
p(\f \,| \,\beta) = \left(\frac{\beta}{2\pi}\right)^{\frac{N}{2}}\exp\left(-\frac{\beta}{2}\Vert(\bI_{N} - \bF^{Geg}\bF^{CT})\f\Vert_2^2 \right).
\end{equation}
A few remarks are in order:
\begin{rem}[Choice of $\lambda$]\label{rem:prior1} 
Our choice for $\lambda$ is motivated by \Cref{thm:GSSV92}, even though our underlying assumption about $f$ being analytic in $[-1,1]$ means that it can be approximated equally well  by any orthogonal polynomial expansion \eqref{eq:orthogsum}, that is for all $\lambda \ge 0$. For the trapezoidal rule to be accurate using equally spaced nodes (\eqref{eq:quadrature} and later \eqref{eq: prior})  a sufficient number of vanishing derivatives of the integrand function at $\pm 1$ is needed, further motivating the use of  ``large enough'' $\lambda$.  By contrast, it was demonstrated in \cite{Boyd2005} that if $\lambda$ is too large, the Gegenbauer reprojection method yields poor results, mainly due to the narrow region of support of $w_\lambda(x)$ (\cref{def:gegpolys}). To this end we note that other Gibbs complementary bases, such as the Freud polynomials used in \cite{GelbTanner}, are more robust. The benefits in this case are only realized when $N$ is very large, a consequence of the resolution requirements for the corresponding polynomial approximation.   As it is important in the context of  inverse problems that the solution be obtainable without a large number of measurements, we restrict our study to Gegenbauer polynomials. Our specific choices for $\lambda$, especially as a function of signal-to-noise ratio SNR, is discussed in more detail in \cref{sec:numparams}. %and choose $\lambda \AG{=}  4$ which is large enough to maintain accuracy in the quadrature error (compared to noise contribution) and small enough to avoid these other issues.  \AG{We also employ $\lambda = 1$ in our numerical experiments for comparison purposes.} 
\end{rem}
\begin{rem}[Choice of  $m$]
\label{rem:prior2} We also choose $m$ to satisfy Theorem \ref{thm:GSSV92}, which follows from \eqref{eq:gegexperror}. Although it is possible to make both $\lambda$ and $m$ unknown random variables,  this will greatly increase the computational cost owing to the computation of  $\bF^{Geg}$ and $\bF^{CT}$ which would no longer be deterministic. In particular the dimensions of $\bF^{Geg}$ would change for each iteration.  Since our prior pertains to the analyticity of the underlying function $f$, adopting $m$  according to \eqref{eq:gegexperror} is reasonable. 
\end{rem}
\begin{rem} [Sparse priors]
\label{rem:prior3}
Priors are often constructed to satisfy the assumption that the underlying solution is sparse in some transform domain, such as the gradient or wavelet domain, see e.g.~\cite{calvetti2019hierachical,glaubitz2023generalized,he2009exploiting,ji2008bayesian} and references therein.  From this perspective, since the decay of the orthogonal polynomial coefficients for $f$ is exponential, we could construct a sparse prior where the transform domain is computed in \eqref{eq: Gegen1 mv} for $\lambda \ge 0$. The resulting prior would be similar to those used in the above references with the notable exception that for $m \sim \kappa N$ satisfying \cref{thm:GSSV92}, we are allowing the underlying function to have considerable variability, which is much less restrictive than using the gradient or low-order wavelet transforms.  Since our goal is to express spectral reprojection in a Bayesian framework, we do not compare our new method to formulations using sparse priors in this investigation.
\end{rem}

%%%%%%%%%%%%%%%%%%%%%%%%%%%%%%

\subsection{Formulating the hyperpriors}
\label{sec:hyperprior}
In \eqref{eqn:bayes} we denoted $p({\btheta})$ as the hyperprior density function for any additional parameters needed to describe the hierarchical model.  As noted there, this amounts to  $p({\btheta}) = p(\gamma)\, p(\beta)$, with $\gamma, \beta$ being the inverse noise variance and inverse prior variance respectively. Choosing $\gamma$ and $\beta$ to be Gamma distributions is convenient since their corresponding  densities  are   conjugate for the corresponding conditionally Gaussian distributions \eqref{eq:likelihood} and \eqref{eq: prior}. 

The Gamma distribution has probability density function 
\begin{equation}\label{eq:pdf_gamma}
\Gamma( \eta \,|\, c,d ) 
= \frac{d^c}{\Gamma(c)} \eta^{c-1} e^{-d\eta},
\end{equation}  
where $c$ and $d$ are positive shape and rate parameters.  
For $c \rightarrow 1$ and $d \rightarrow 0$, \eqref{eq:pdf_gamma} is an uninformative prior, which is desirable since we do not want these parameters to have undue influence on the solution. Thus here for simplicity we choose
\begin{subequations} 
\begin{equation}
\label{eq:gammaprior} 
p(\gamma)= \Gamma(\gamma \, \vert \, c, d),   
\end{equation}
\begin{equation} 
\label{eq:betaprior}
p(\beta)= \Gamma(\beta \, \vert \, c, d).
\end{equation}
\end{subequations}
with $c=1$ and $d=10^{-4}$ for all of our numerical examples.\footnote{For more detailed explanations regarding this choice of hyperprior, see \cite{calvetti2019hierachical,glaubitz2023generalized} and references therein.}  
 Following the above discussion, we then replace $p({\btheta})$ in \eqref{eqn:bayes} with
\begin{equation} \label{eq: hyper-prior}
p({\btheta}) = p(\gamma)p(\beta) \propto \gamma^{c-1} \exp(-d \, \gamma) \,
\beta^{c-1} \exp (-d \, \beta).
\end{equation}

%%%%%%%%%%%%%%%%%%%%%%%%%%%%%%

\subsection{Bayesian inference}\label{sec:bayesinference}
We now have all of the ingredients needed to construct \eqref{eqn:bayes}.  Specifically, substituting in  the likelihood \eqref{eq:likelihood}, the prior \eqref{eq: prior} and the hyperprior \eqref{eq: hyper-prior} density functions into \eqref{eqn:bayes}, we obtain the posterior density function approximation
\begin{gather}\label{eq:posterior}
p(\f, \gamma, \beta \, \vert \, \f_{m,N}^\lambda) \propto 
\gamma^ {c+\frac{N}{2}-1} \beta^{c+\frac{N}{2} -1} \exp \left\{ -\frac{\gamma}{2} \Vert \f_{m,N}^\lambda - \bF^{Geg}\bF^{CT}\, \f \Vert_2^2 
-\frac{\beta}{2} \Vert (\bI_{N} - \bF^{Geg}\bF^{CT})\f \Vert_2^2 -d \, \gamma
-d \, \beta
\right\}.
\end{gather}
From here we can  obtain the MAP estimate, which is equivalent to minimizing the negative logarithm of the right hand side of \eqref{eq:posterior}, yielding the objective function 
\begin{eqnarray}
\label{eq:objective1}
{\cJ}(\f, \gamma, \beta)&=&-(c+\frac{N}{2}-1) \log(\gamma) -(c+\frac{N}{2}-1) \log(\beta) + \frac{\gamma}{2} \Vert \f_{m,N}^\lambda - \bF^{Geg}\bF^{CT}\, \f \Vert_2^2   \nonumber \\ 
& + & \frac{\beta}{2} \Vert (\bI_{N} - \bF^{Geg} \bF^{CT})  \f  \Vert_2^2
 +d \, \gamma  +d \,\beta   
\end{eqnarray}
up to a constant that is independent of $\f$, $\gamma$, and $\beta$.  In what follows we discuss how to minimize \eqref{eq:objective1}.
%We call our new approach the {\em Bayesian Spectral Reprojection (BSR)} method.
\subsubsection{The block-coordinate descent approach} \label{sub:BCD}
We use a block-coordinate descent (BCD) approach \cite{wright2015coordinate,beck2017first} to approximate the MAP estimate of \eqref{eq:posterior}. In particular, BCD minimizes  $\cJ$ in \eqref{eq:objective1} by alternating minimization  w.r.t.\ $\f$ for fixed $\gamma, \beta$,  minimization w.r.t.\ $\gamma$ for fixed $\f,\beta$, and  minimization  w.r.t.\ $\beta$ for fixed $\f,\gamma$.   That is, given an initial guess for the hyperparameters $\gamma$ and $\beta$, the BCD algorithm proceeds through a sequence of updates of the form 
\begin{subequations}
\label{eq:BCD} 
\begin{equation}
\f = \underset{\f}{\text{arg min}} \left\{ {\cJ}(\f,\gamma, \beta ) \right\}, \label{eq:BCD_f}
\end{equation}
\begin{equation}
\gamma  = \underset{\gamma}{\text{arg min}} \left\{ {\cJ}(\f,\gamma, \beta ) \right\}, \label{eq:BCD_gamma}
\end{equation}
\begin{equation}
\beta = \underset{\beta}{\text{arg min}} \left\{ {\cJ}(\f,\gamma, \beta ) \right\},\label{eq:BCD_beta}
\end{equation}
\end{subequations}
until a convergence criterion is met.
In our implementation we initialize the hyperparameter vectors as $\gamma = 1$ and $\beta = 1$, and we stop when either the relative change of the objective function value  or the relative change in $\f$ falls below a given threshold.   Efficient implementation of the three update steps in \eqref{eq:BCD} is described below.

\subsubsection{Updating the unknown signal $\f$} 
\label{sub:f_update}
Due to our intentionally chosen conditionally Gaussian prior \eqref{eq: prior} and Gamma hyperprior \eqref{eq:gammaprior} -- \eqref{eq:betaprior}, for fixed $\gamma$ and $\beta$ we have 
\begin{equation} \label{eq: f density}
p(\f , \gamma, \beta \,|\, \f_{m,N}^\lambda) \propto p(\f_{m,N}^\lambda \, | \, \f, \gamma) \, p(\f \, | \, \beta) \propto \cN (\bmu, C), 
\end{equation}
where
\begin{equation} \label{eq: f mu}
\bmu = \gamma \, C ({\bF^{Geg} \bF^{CT}})^{T} \, \f_{m,N}^\lambda
\end{equation}
and
\begin{equation} \label{eq: f C}
C = \left(\gamma \, ({\bF^{Geg} \bF^{CT}})^{T}  ({\bF^{Geg} \bF^{CT}}) + \beta \,(\bI_{N} - \bF^{Geg} \bF^{CT}) ^ T (\bI_{N} - \bF^{Geg} \bF^{CT})\right)^{-1}.
\end{equation}
The update step for $\f$ in \eqref{eq:BCD_f} therefore reduces to solving the linear system
\begin{equation} \label{eq: f_eqn}
\left(\gamma \, ({\bF^{Geg} \bF^{CT}})^{T}  ({\bF^{Geg} \bF^{CT}}) + \beta \,(\bI_{N} - \bF^{Geg} \bF^{CT}) ^ T (\bI_{N} - \bF^{Geg} \bF^{CT})\right) \f = \gamma \, ({\bF^{Geg} \bF^{CT}})^{T} \, \f_{m,N}^\lambda.
\end{equation}
We also note that by construction the \emph{common kernel condition} (see e.g.~\cite{kaipio2006statistical})
\begin{equation}\label{eq:common_kernel}
\text{kernel}( \bF^{Geg} \bF^{CT}) \cap \text{kernel}(\bI_{N} - \bF^{Geg} \bF^{CT})  = \{ \0 \}
\end{equation} 
is automatically satisfied, and hence \eqref{eq: f_eqn} has a unique solution which can be directly obtained.

\begin{rem}
\label{rem:common kernel}
Recall that the kernel of an operator $G: \bbC^N \to \bbC^M$ is given by $\text{kernel}(G) = \{ \, \f \in \bbC^N \mid G \f = \0 \, \}$, that is, the set of vectors that are mapped to zero.
The common kernel condition \eqref{eq:common_kernel} ensures that the prior (regularization) introduces a sufficient amount of complementary information to the likelihood (the given measurements) to make the problem well-posed, which is a standard assumption in regularized inverse problems \cite{kaipio2006statistical,tikhonov2013numerical}. 
\end{rem}

Analogous to the objective function \eqref{eq:objective1}, we also note that \eqref{eq: f_eqn} is equivalent to updating $\f$ via the quadratic optimization problem 
\begin{equation}\label{eq:f_update} 
\f = \underset{\f}{\text{arg min}} \left\{  \frac{\gamma}{2} \Vert \f_{m,N}^\lambda - \bF^{Geg}\bF^{CT}\, \f \Vert_2^2   +\frac{\beta}{2} \Vert (\bI_{N} - \bF^{Geg} \bF^{CT})  \f  \Vert_2^2 \right\}.
\end{equation}
Various methods can be used to  efficiently solve  \eqref{eq:f_update},  including the fast iterative shrinkage-thresholding (FISTA) algorithm \cite{beck2009fast}, the preconditioned conjugate gradient (PCG) method \cite{saad2003iterative}, potentially combined with an early stopping based on Morozov's discrepancy principle \cite{calvetti2015hierarchical,calvetti2018bayes,calvetti2020sparse}, and the gradient descent approach \cite{glaubitz2023generalized}. 
There is no general advantage of one method over another, and the choice should be made based on the specific problem at hand.

\subsubsection{Updating the hyperparameters ${\gamma}$ and ${\beta}$}
\label{sub:alpha_update}

We now address the update for the hyperparameter $\gamma$  for fixed $\f$ and $\beta$. Substituting \eqref{eq:objective1} into \eqref{eq:BCD_gamma} and ignoring all terms that do not depend on $\gamma$ yields
\begin{equation}\label{eq:update_alpha2} 
\gamma =\underset{\gamma}{\text{arg min}}\left\{\cJ_\gamma\right\} = \underset{\gamma}{\text{arg min}} \left \{ -(c+\frac{N}{2}-1) \log(\gamma)  + \frac{\gamma}{2} \Vert \f_{m,N}^\lambda - \bF^{Geg}\bF^{CT}\, \f \Vert_2^2   +d \, \gamma   \right \}. 
\end{equation} 
Solving $\ds \frac{d\cJ_\gamma}{d\gamma} = 0$ then provides
\begin{equation}\label{eq:update_alpha} 
\gamma = \frac{2c+N-2}{\Vert \f_{m,N}^\lambda - \bF^{Geg}\bF^{CT}\, \f \Vert_2^2+2d}
\end{equation} 
for the $\gamma$-update.  The $\beta$-update is similarly obtained as
\begin{equation}\label{eq:update_beta} 
\beta = \frac{2c+N-2}{ \Vert (\bI_{N} - \bF^{Geg} \bF^{CT})  \f  \Vert_2^2+2d}.
\end{equation}
\Cref{alg:bayesapprox} summarizes the iterative process for finding the point estimate solution $\f$ in \eqref{eq:BCD} along with precision values $\gamma$ and $\beta$.

\begin{algorithm}[h!]
\caption{The Bayesian spectral reprojection (BSR) MAP estimate}% from noisy Fourier data using objective function \eqref{eq:objective1}}
\label{alg:bayesapprox}
\begin{algorithmic}
\item[] {\bf Input} Data vector $\wh{\bb} \in {\mathbb C}^N$ in \eqref{eq:noisyfourier1d}, %for known noise variance $\sigma^2$,
choices for $\lambda = \kappa_1 N$ and $m = \kappa_2 N$, $0 < \kappa_1, \kappa_2 < 1$, and grid points $x_j$, $j = 1,\dots,N$.
\item[] {\bf Output} MAP estimates for $\f$, $\gamma$, and $\beta$.
\item[Initialize $\gamma$ and $\beta$. We choose $\gamma = \beta = 1$.]
\item[Calculate the likelihood $\f_{m,N}^\lambda$ according to \eqref{eq: fmNlambda mv}.]
\item[\bf{repeat}] \quad \\
 Update $\f$ according to \eqref{eq: f_eqn} \\
 Update $\gamma$ according to \eqref{eq:update_alpha}\\
 Update $\beta$ according to \eqref{eq:update_beta} 
\item[\bf{until}] convergence or the maximum number of iterations is reached
\end{algorithmic}
\end{algorithm}

\begin{remark}[Fixed hyperparameters]
\label{rem:fixedgammabeta} We note that the posterior formulation in \cref{eq:posterior} allows for uncertainty quantification for the recovered solution $\f$ for {\em fixed} $\gamma$ and $\beta$.  Specifically, from  \eqref{eq: f density} we have
\begin{equation} \label{eq: f density 2}
p(\f \,|\, \f_{m,N}^\lambda) \propto p(\f_{m,N}^\lambda \, | \, \f) \, p(\f ) \propto \cN (\bmu, C), 
\end{equation}
where the mean $\bmu$ and the covariance matrix $C$ are again respectively given by \eqref{eq: f mu} and \eqref{eq: f C}. We can then  sample directly from  $\cN (\bmu, C)$ to obtain credible intervals\footnote{Related to the (frequentist) confidence interval, a (Bayesian) credible interval is used to characterize a probability distribution, and is also sometimes referred to as a Bayesian posterior interval \cite{Laud1999BayesianSA, gelman2003bayesian}.} for every component of the solution $\f$. When $\gamma$ and $\beta$ are themselves random variables, computational sampling is needed to recover \eqref{eq:posterior}, for instance, full posterior sampling using a Markov chain Monte Carlo (MCMC) method \cite{mcbook}, typically requiring considerable computational expense.  
\end{remark}

\section{Generalized Bayesian spectral reprojection}
\label{sec:Bayes}

By design, and as will be demonstrated via numerical experiments in \Cref{sec:numerical}, the BSR method in  \Cref{alg:bayesapprox} yields comparable point estimates as the original spectral reprojection method in \Cref{alg:gegapprox}.  We are furthermore able to quantify the uncertainty of the  approximation for fixed hyperparameters $\gamma$ and $\beta$ (see \cref{rem:fixedgammabeta}). Several potential shortcomings of the BSR approach  are readily apparent, however.  Specifically, methods for computational inverse problems are designed to promote solutions that yield close approximations to the observable information while also being constrained by assumptions regarding the underlying signal. In addition to mitigating the effects of ill-posedness, the regularization (prior) term also serves to avoid over-fitting. By construction, the likelihood and prior terms in posterior density function  \cref{eq:posterior} are  {\em both}  based on the assumption that $f \in \mathbb{P}_m^\lambda$, that is, the space of Gegenbauer polynomials given by \cref{eq:orthogsum}. While such an assumption is suitable for solving the {\em forward} problem in \Cref{alg:gegapprox}, it inherently limits the space of possible solutions when considering the {\em inverse} problem.  A more natural approach for solving linear inverse problems with i.i.d.~additive (complex) Gaussian noise is to {\em directly} fit the given data model in the construction of the likelihood, and then determine a suitable prior (regularization) for fitting the data model. In this regard, our insights in using the spectral reprojection prior \eqref{eq: prior} still remain, in particular pertaining to the discussion in \cref{rem:prior1} and \cref{rem:prior2}.

Given these insights, we now introduce the {\em generalized} Bayesian spectral reprojection (GBSR) method.  We demonstrate its advantages in generating the MAP estimate with regard  to accuracy and robustness in our numerical experiments in \cref{sec:numerical}, where we also display credible intervals for fixed hyperparameters.

%%%%%%%%%%%%%%%%%%%%%%%%%%%%%%%%%%%

\subsection{Formulating the inverse problem}
\label{sec:invproblemforFourier}
We begin by considering the linear inverse problem
\begin{equation}
\wh\cB=\cF\cX+ \wt\cE, 
\label{eq:genbayesmodel 2}
\end{equation}
where $\cX$, $\wh\cB$, and $\wt\cE$ are random variables defined over a common probability space, and $\cX$ and $\wt\cE$ are assumed to be independent. As before, $\cX$ represents the unknown which we seek to recover, $\wh\cB$ denotes the observables, $\cF$ is a known linear operator, and $\wt\cE\sim \cC\cN(\0,\wt{\gamma}^{-1} \bI_N)$ is a vector of independent and identically distributed (i.i.d.) circularly symmetric complex Gaussian noise with mean $0$ and unknown variance $\wt{\gamma}^{-1}$.\footnote{Unlike the BSR formulation in \cref{eq:genbayesmodel}, since no transformation is used in \cref{eq:genbayesmodel 2}, the noise does truly follow $\cC\cN(\0,\alpha^{-1} \bI_N)$ for $\alpha^{-1}$ defined in \eqref{eq:noisyfourier1d}, that is, $\wt{\gamma}=\alpha$.} The realizations of the random variables in \eqref{eq:genbayesmodel 2} are modeled in \eqref{eq:observationFourier}, respectively by $\f \in \bbR^N$, $\wh{\bb} \in \bbC^N$, and $\wt\beps \in \bbC^N$.

%the discretized version of \eqref{eq:fourcoeff1D}, that is
%\begin{equation*}\label{eq: model}
%\wh{\bb}=\bF \, \f + \beps,
%\end{equation*}
%once again noting that $\bF \in \bbC^{N\times N}$ the DFT matrix \eqref{eq:DFT}, $\f$ is the unknown solution at the uniform grid points given by  \eqref{eq:uniformgrid}. We remark that the choice of DFT matrix \eqref{eq:DFT} for the linear operator $\bF$ in \eqref{eq: model} is reasonable as we want our reconstruction to be smooth and live on equally spaced points. 
We now proceed to formulate the corresponding posterior density function given by
\begin{equation}\label{eqn:bayes2}
p(\f,\wt\btheta \,|\,\wh{\bb}) \propto p(\wh{\bb} \, | \,\f,\wt\btheta) \, p(\f \,| \,\wt\btheta) \, p(\wt\btheta),
\end{equation}
where, analogously to \eqref{eqn:bayes}, $p(\wh{\bb} \,| \,\f,\wt\btheta)$ is the likelihood density function determined by $\bF$ in \eqref{eq:observationFourier} and assumptions on $\wt\cE$,  $p(\f\,|\,\wt\btheta)$ is the density of the prior distribution encoding \emph{a priori} assumptions about the solution, and $p(\wt\btheta)$ is the hyperprior density for all other involved parameters, which we denote collectively as $\wt\btheta$.  As was also done in \Cref{sec:Bayesspectral}, we write  $p(\wt\btheta) = p(\wt{\gamma}) \, p(\beta)$, where $\wt{\gamma}$ and $\beta$ correspond to the inverse noise and inverse prior variances respectively.

%To this end, we recall that the prior remain the same as described in Section \ref{sec:prior1}, but we need to formulate a likelihood density function and corresponding hyperprior density function.

%%%%%%%%%%%%%%%%%%%%%

\subsection{Formulating the likelihood} \label{sec:likelihood2}
%The likelihood density function $p(\wh{\bb} \,| \,\f, \alpha)$ in \eqref{eqn:bayes2} models the connection between the signal and observation vectors.  
Since $\wt\cE \sim\cC \cN(\0,\wt{\gamma}^{-1} \bI_N)$ for unknown variance $\wt{\gamma}^{-1}$, we have
\[
p(\wt\beps \,| \, \wt{\gamma})= \left(\frac{\wt{\gamma}}{\pi}\right)^{N}\exp \left(-\wt{\gamma} \Vert \wt\beps \Vert _2^2 \right).
\]
It follows from \eqref{eq:observationFourier} that the likelihood density function is then
\begin{equation}\label{eq:likelihood2}
 p(\wh{\bb}\,|\,\f,\wt{\gamma})   =
\left(\frac{\wt{\gamma}}{\pi}\right)^{N}\exp(-\wt{\gamma} \Vert \wh{\bb} - \bF\, \f \Vert_2^2),
\end{equation}
with corresponding distribution 
\begin{equation}
\label{eq:likelihooddistribution2}
\wh\cB \sim \cC \cN (\bF \, \f, \wt{\gamma}^{-1} \bI_N ).    
\end{equation}

%%%%%%%%%%%%%%%%%%%%%

\subsection{Formulating the prior and the hyperpriors}
\label{sec:prior2}

For reasons already mentioned, we again choose the prior density function to be used in  \eqref{eqn:bayes2}  as \eqref{eq: prior}.  We also note that in the context of the likelihood density function given by \eqref{eq:likelihood2}, our choice of prior means that the resulting objective function yields enough overlap between the data (likelihood) and regularization (prior) terms to determine a meaningful solution, beyond satisfying the common kernel condition.\footnote{Since ${\bF}$ is invertible, $ker(\bF) = \{ \0\}$ and the kernel condition for  GBSR is automatically satisfied.} 
Specifically,  for $r > 0$ and constants $Const_1, C > 0$, \eqref{eq:gegexperror} yields
$$\langle P_{m}^{\lambda} f - f, e^{i k \pi x} \rangle \le Const_1\int_{-1}^1  e^{-m^r} e^{-ik\pi x}dx  < Ce^{-rm},$$ and in fact is only nonzero for $k = 0$. Since $m \propto N$, which is assumed to adequately resolve the underlying signal $\f$ from $\wh\bb$, we conclude that the solution $\f$ can be recovered for our choice of prior. Furthermore, although only  $\lambda \ge 0$ is required to satisfy \eqref{eq:gegexperror}, as discussed in \cref{rem:prior1}, $\lambda$ still must be  large enough for \eqref{eq:quadrature} to accurately approximate $\{\wh{f}^\lambda_\ell\}_{\ell = 1}^m$  in \cref{eq:gegcoeff}, while small enough to maintain robustness. A complete discussion for choosing in the context of noise is provided in \cref{sec:numparams}.

Finally, the hyperpriors $\wt{\gamma}$ and $\beta$ are again modeled as Gamma distributions,  yielding the hyperprior density function \eqref{eq: hyper-prior} density function with $\wt{\gamma}$ replacing $\gamma$ as the likelihood hyperprior, i.e.
\begin{equation} \label{eq: hyper-prior_copy}
p(\wt \btheta) = p(\wt{\gamma})p(\beta) \propto \wt{\gamma}^{c-1} \exp(-d \, \wt{\gamma}) \,
\beta^{c-1} \exp (-d \, \beta).
\end{equation}

%%%%%%%%%%%%%%%%%%%%%

\subsection{Bayesian inference}\label{sec:bayesinference2}

Substituting the likelihood \eqref{eq:likelihood2}, the prior \eqref{eq: prior} and the hyperprior \eqref{eq: hyper-prior_copy} density functions into \eqref{eqn:bayes2}, we obtain the posterior density function approximation
\begin{gather}\label{eq:posterior2}
p(\f, \wt{\gamma}, \beta \, \vert \, \wh{\bb} ) \propto 
\wt{\gamma}^ {c+N-1} \beta^{c+\frac{N}{2} -1} \exp \left\{ -\wt{\gamma} \Vert \wh{\bb} - \bF\, \f \Vert_2^2 
-\frac{\beta}{2} \Vert (\bI_{N} - \bF^{Geg}\bF^{CT})\f \Vert_2^2 -d \, \wt{\gamma} 
-d \, \beta
\right\}.
\end{gather}
The MAP estimate of \eqref{eq:posterior2} is the minimizer of the objective function
\begin{align}
\label{eq:objective2}
\wt\cJ(\f, \wt{\gamma}, \beta)=-(c+N-1) \log(\wt{\gamma}) -(c+\frac{N}{2}-1) \log(\beta) + \wt{\gamma} \Vert \wh{\bb} - \bF\, \f \Vert_2^2   +\frac{\beta}{2} \Vert (\bI_{N} - \bF^{Geg} \bF^{CT})  \f  \Vert_2^2
 +d \, \wt{\gamma}  +d \,\beta    
\end{align}
up to a constant that is independent of $\f$, $\wt{\gamma}$, and $\beta$. The BCD algorithm iteratively updates the solution according to
\begin{subequations}
\label{eq:updates2}
\begin{equation}
\ds \left(2\wt{\gamma} \, {\bF}^{\ast}  {\bF} + \beta \,(\bI_{N} - \bF^{Geg} \bF^{CT}) ^ T (\bI_{N} - \bF^{Geg} \bF^{CT})\right) \f = 2\wt{\gamma} \, {\bF}^{\ast} \, \wh{\bb}, \label{eq:f_update2}
\end{equation}
\begin{equation}
\ds \wt{\gamma} = \frac{c+N-1}{\Vert \wh{\bb} - \bF\, \f \Vert_2^2+d}, \label{eq:alpha_update2}
\end{equation}
\begin{equation}
\ds \beta = \frac{2c+N-2}{ \Vert (\bI_{N} - \bF^{Geg} \bF^{CT})  \f  \Vert_2^2+2d}.
\label{eq:beta_update2}
\end{equation}
\end{subequations}

%\AG{Put algorithm here.}

\Cref{alg:bayesapprox2} summarizes the iterative process for minimizing the objective function \eqref{eq:objective2}, which provides the point estimate solution $\f$ along with precision values $\wt{\gamma}$ and $\beta$. Analogous to the formulation obtained in \cref{rem:fixedgammabeta}, for fixed $\wt{\gamma}$ and $\beta$ it is possible to directly sample from the resulting posterior to obtain credible intervals.

\begin{algorithm}[h!]
\caption{The Generalized Bayesian spectral reprojection (GBSR) MAP estimate }
\label{alg:bayesapprox2}
\begin{algorithmic}
\item[] {\bf Input} Data vector $\wh{\bb} \in {\mathbb C}^N$ in \eqref{eq:noisyfourier1d}, 
choices for $\lambda = \kappa_1 N$ and $m = \kappa_2 N$, $0 < \kappa_1, \kappa_2 < 1$, and grid points $x_j$, $j = 1,\dots,N$.
\item[] {\bf Output} MAP estimates for $\f$, $\wt{\gamma}$, and $\beta$.
\item[Initialize $\wt{\gamma}$ and $\beta$ as $\wt{\gamma} = \beta = 1$.]
\item[\bf{repeat}] \quad \\
 Update $\f$ according to \eqref{eq:f_update2} \\
 Update $\wt{\gamma}$ according to \eqref{eq:alpha_update2}\\
 Update $\beta$ according to \eqref{eq:beta_update2} 
\item[\bf{until}] convergence or the maximum number of iterations is reached
\end{algorithmic}
\end{algorithm}

\section{Numerical experiments}
\label{sec:numerical}

We now demonstrate the efficacy of the  BSR (\cref{alg:bayesapprox}) and GSBR (\cref{alg:bayesapprox2}), and compare each method's  performance to the original spectral reprojection in \cref{alg:gegapprox}.  We also test our new  methods for varying signal-to-noise ratio (SNR), which we define as
\begin{equation}\label{eq: SNR}
\text{SNR} = 10 \log_{10} \left( \frac{\Vert \bF \f \Vert_2^2} {N \alpha^{-1}}\right).
\end{equation}
In each example the noisy Fourier data $\wh{\bb}$  is given by \eqref{eq:noisyfourier1d} with noise variance $\alpha^{-1}$ corresponding to \cref{eq: SNR}. The noiseless $N$ Fourier coefficients $\wh{\f}$  are constructed via trapezoidal rule \cref{eq: trapezoid} for a highly resolved mesh ($8N$).  

\subsection{An illustrative example}
\label{sec:methodperformance}
We first seek to make a general comparison of the Fourier partial sum approximation given by \eqref{eq:partial_four_sum1 mv2}, the (Gegenbauer)  spectral reprojection  in \cref{alg:gegapprox}, the BSR approach in \cref{alg:bayesapprox}, and finally the GBSR  method in \cref{alg:bayesapprox2}.  We choose an example that is readily recovered for large enough $N$ (number of given Fourier coefficients) in the noise-free case using spectral reprojection for suitable parameter choices $m$ and $\lambda$. In particular, we consider the function

\begin{equation} \label{ex:example1}
f(x)=e^x \sin(5x), \quad x \in [-1,1]. 
\end{equation}
We assume we are given the first $N = 128$ noisy Fourier samples $\wh{\bb}$ in \cref{eq:noisyfourier1d} with noise variance $\alpha^{-1} = 2 \times 10^{-3}$ (equivalently SNR $\approx 5.95$). 
The spectral reprojection parameters   $m=9$ and $\lambda=4$ are chosen to satisfy \Cref{thm:GSSV92} and also to ensure accurate and robust recovery in the noise-free case. 

\begin{figure}[h!]
\centering
\includegraphics[width=0.24\textwidth]{Numerical1/3_fN_N128_mu4_NGeg9.pdf}
\includegraphics[width=0.24\textwidth]{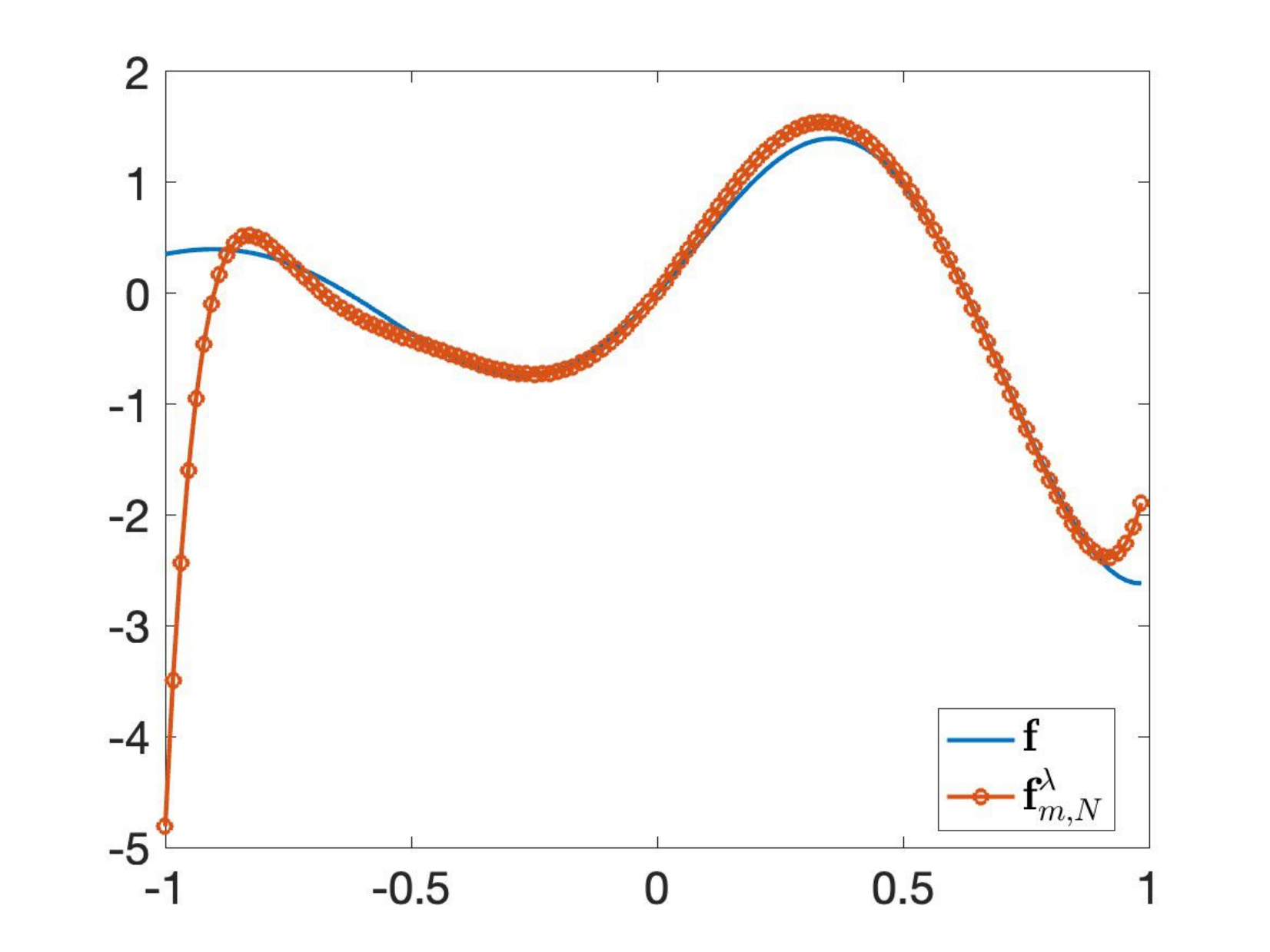}
\includegraphics[width=0.24\textwidth]{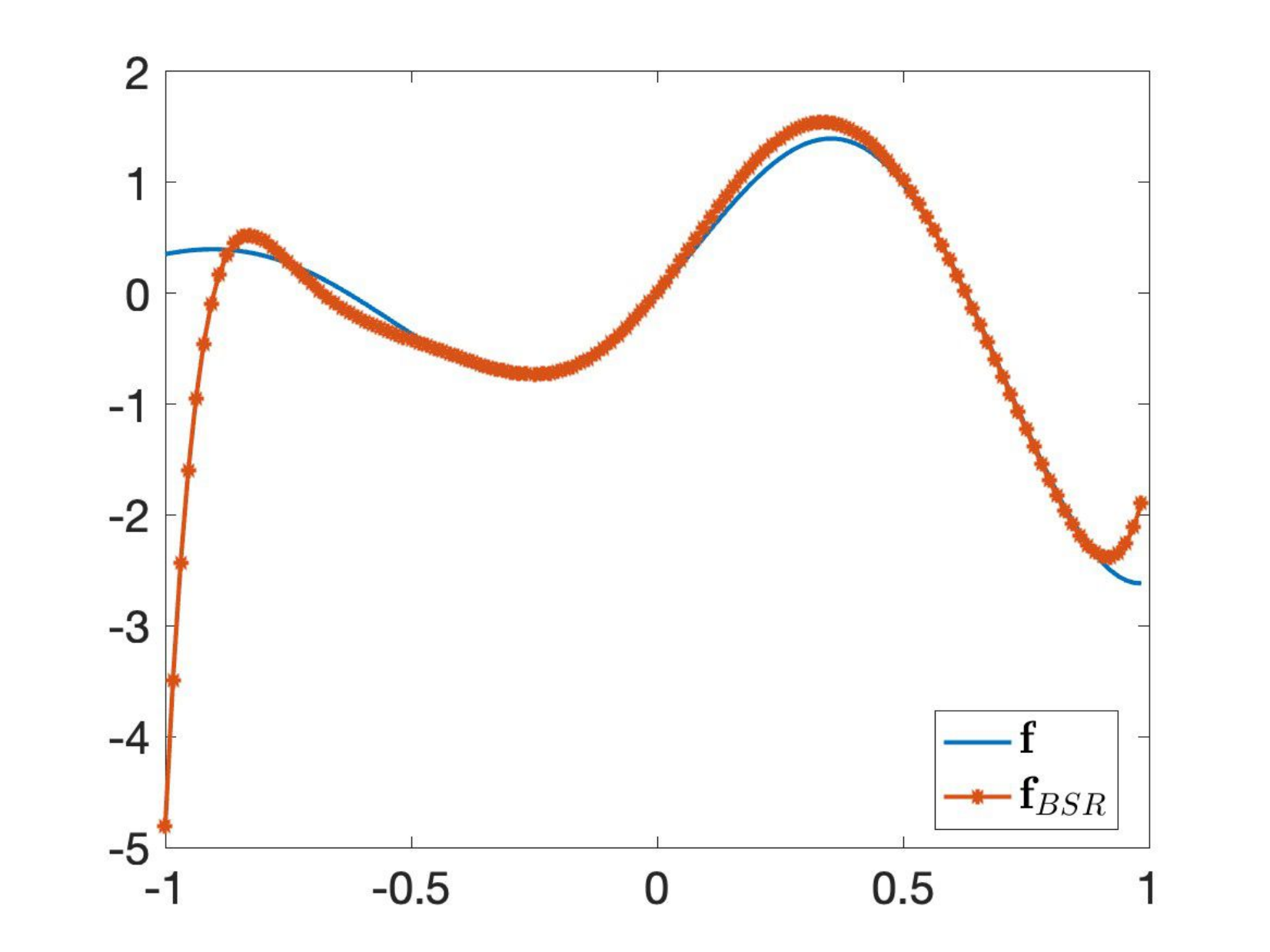} 
\includegraphics[width=0.24\textwidth]{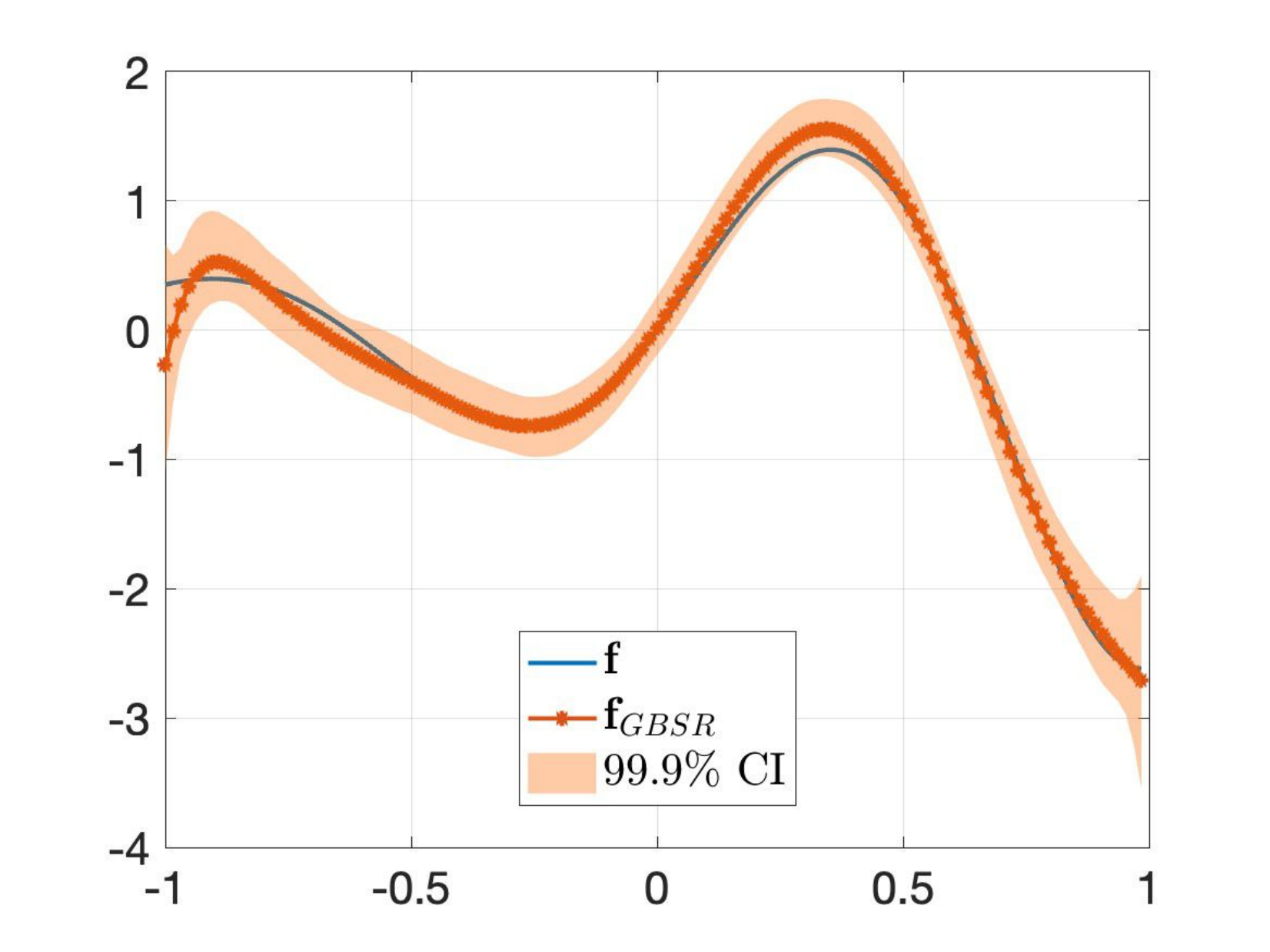}
\caption{Recovery of $f(x) = e^x\sin{(5x)}$. (left) Fourier reconstruction $\f_N$ in \eqref{eq:partial_four_sum1 mv2}; (center-left) Spectral reprojection $\f_{m,N}^{\lambda}$ in \cref{alg:gegapprox}; (center-right) $\f_{BSR}$ in \cref{alg:bayesapprox}; and (right) $\f_{GBSR}$ in \cref{alg:bayesapprox2}.}.
\label{fig: 1}
\end{figure}

Figure \ref{fig: 1} compares the recovery of \cref{ex:example1} using \cref{eq:partial_four_sum1 mv2}, \cref{alg:gegapprox}, \cref{alg:bayesapprox} and \cref{alg:bayesapprox2}.  The $99.9\%$ credible intervals of the solution posterior $p(\f \,|\, \wh{\bb})$ using the final (fixed) estimates of the hyperparameters $\wt{\gamma}$ and $\beta$, as described in \cref{rem:fixedgammabeta}, are also displayed. In addition to the Gibbs phenomenon artifacts manifested by the over/undershoots at the boundaries,
the Fourier partial sum $\f_N$ clearly yields a noisy approximation.   Standard filtering would reduce the oscillations, but also cause excessive blurring \cite{HGG}. \cref{alg:gegapprox} improves accuracy away from the boundaries, but due to the added noise there are large deviations at the boundaries, especially near $x=-1$. In this regard we note that it was shown in \cite{ArchGelb2002}  that the variance near the boundaries using the Gegenbauer reconstruction method is {\em not} reduced when compared to the standard Fourier partial sum. As expected, the point estimate in \cref{alg:bayesapprox} yields nearly identical results to those obtained in \cref{alg:gegapprox}.  The corresponding credible intervals obtained  for fixed hyperparameters $\gamma$ and $\beta$ (see \cref{rem:fixedgammabeta}) are very tight and therefore not shown.  This is also to be expected following the discussion at the beginning of \Cref{sec:Bayes}, and indeed served to motivate the GBSR approach.  Finally, we observe that \cref{alg:bayesapprox2}  yields the overall  best results, with better performance throughout the domain. The credible interval is also displayed in \Cref{fig: 1}(right).

\begin{rem}\label{rem:freudpolys}
Due to issues of ill conditioning leading to the Runge phenomenon, the Gegenbauer spectral reprojection may perform poorly for large $\lambda$ even though \Cref{thm:GSSV92} is satisfied \cite{Boyd2005}. The robust Gibbs complementary (Freud) basis proposed in \cite{GelbTanner} as the reprojection basis $\{\psi_l^\lambda\}_{l = 0}^m$ in \cref{eq:orthogsum} helps alleviate these issues, but requires more points per wave (large $N$) to fully resolve the approximation.  The approach moreover was considered only for noise-free data. Using the Freud reprojection basis did not yield any noticeable improvement in our experiments.  Therefore, given the additional cost incurred in using large $N$, we do not consider it further.
\end{rem}

\subsection{Noise analysis}
\label{sec:SNRcase}

Now we investigate a more typical scenario in computational inverse problems, that is, with fewer observables and possibly lower SNR.   Since it was  used in \cite{Gottlieb97},  we also consider the non-periodic analytic function $f(x)=\cos(1.4\pi(x+1))$ on domain $[-1,1]$. We assume we are given $N = 48$ noisy Fourier coefficients $\wh{\bb}$ in \eqref{eq:noisyfourier1d} with different SNR levels, where $\wh{\f}$  are constructed from \cref{eq: trapezoid} using a refined mesh $8N$. %\AG{Following \cite{Gottlieb97},} we choose parameters $m = 9$ and $\lambda = 4$.

\begin{remark}[Noisy versus noise-free data]
For the {\em noiseless} case, it was shown in \cite{Gottlieb97} that choosing $m = \lambda = 9$  yields more accurate results given $N = 48$ Fourier coefficients $\wh{\f}$. % than our chosen values $m = 9$, $\lambda = 4$.  
However, when the Fourier coefficients contain even small amounts of noise, the performance quality of the spectral reprojection method rapidly deteriorates. In particular, the ill-conditioning attributed to larger values of $\lambda$ exacerbates the impact of noise in the approximation. To mitigate these effects,  we choose $\lambda = 4$ in our numerical experiments.  Note that $m$, which  corresponds to the number of Gegenbauer polynomials needed to resolve the underlying function,  is less impacted by noise. Hence $m = 9$ is still appropriate.
\end{remark}

\begin{figure}[h!]
\centering
\includegraphics[width=0.24\textwidth]{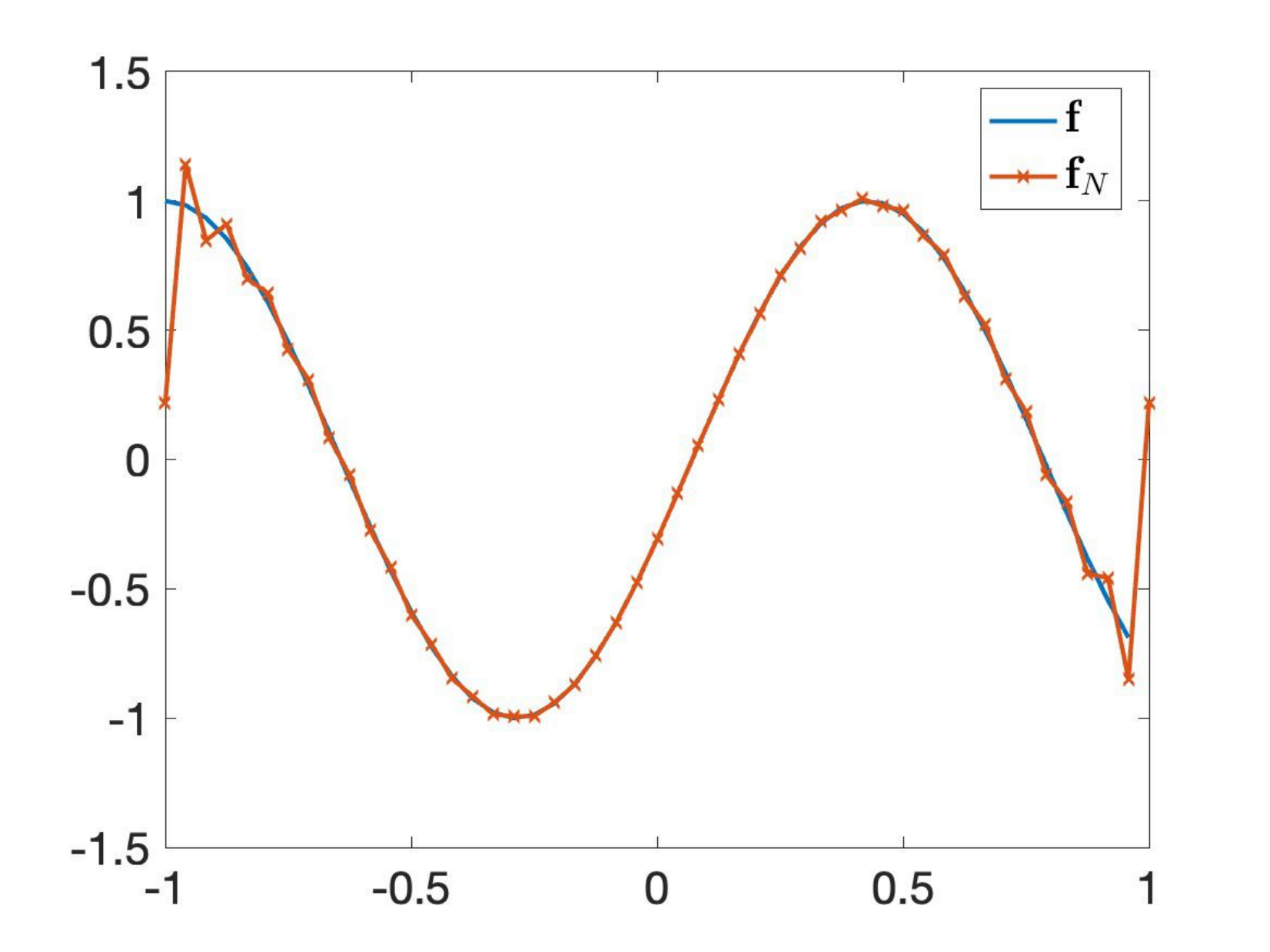}
\includegraphics[width=0.24\textwidth]{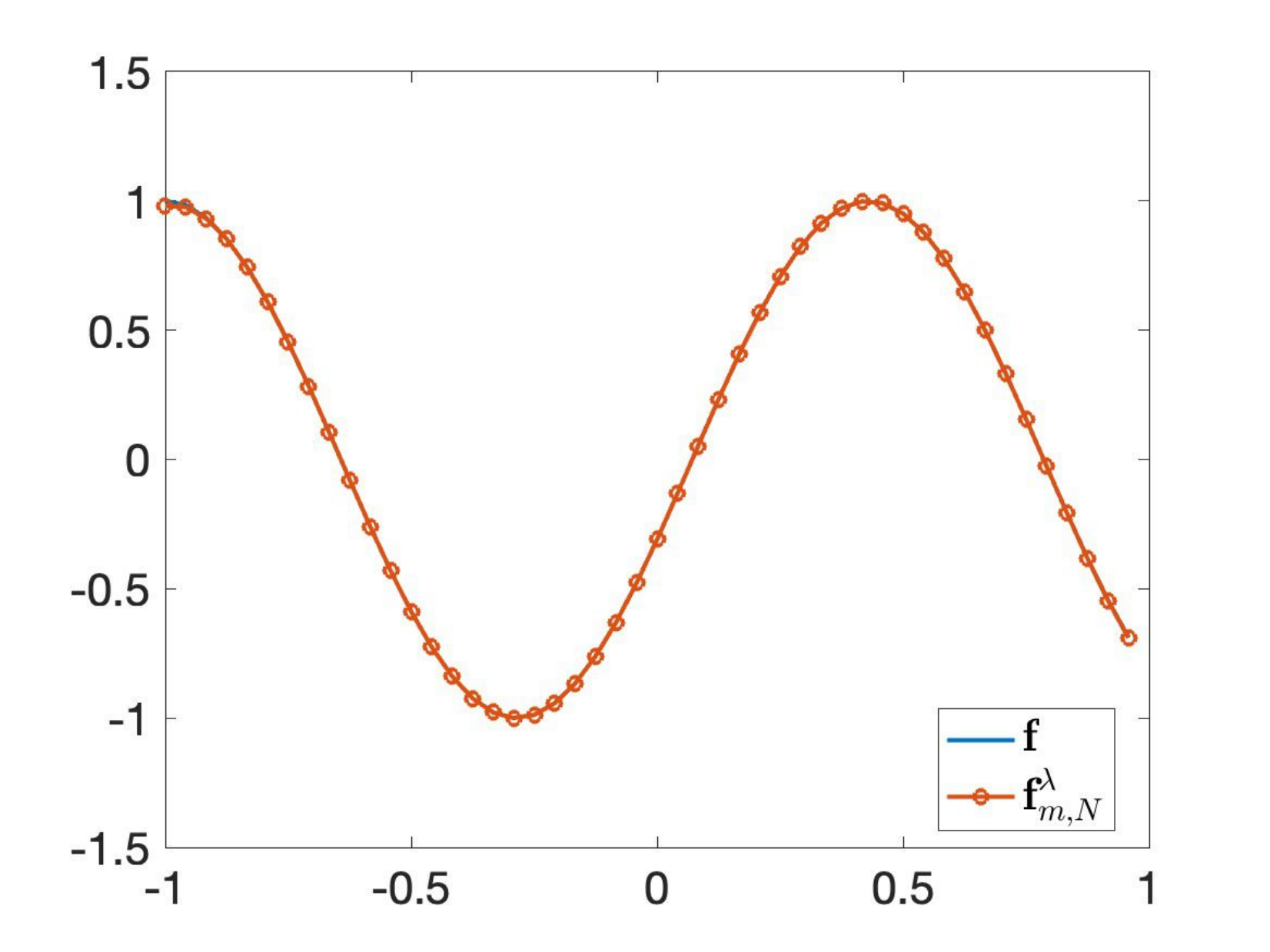}
\includegraphics[width=0.24\textwidth]{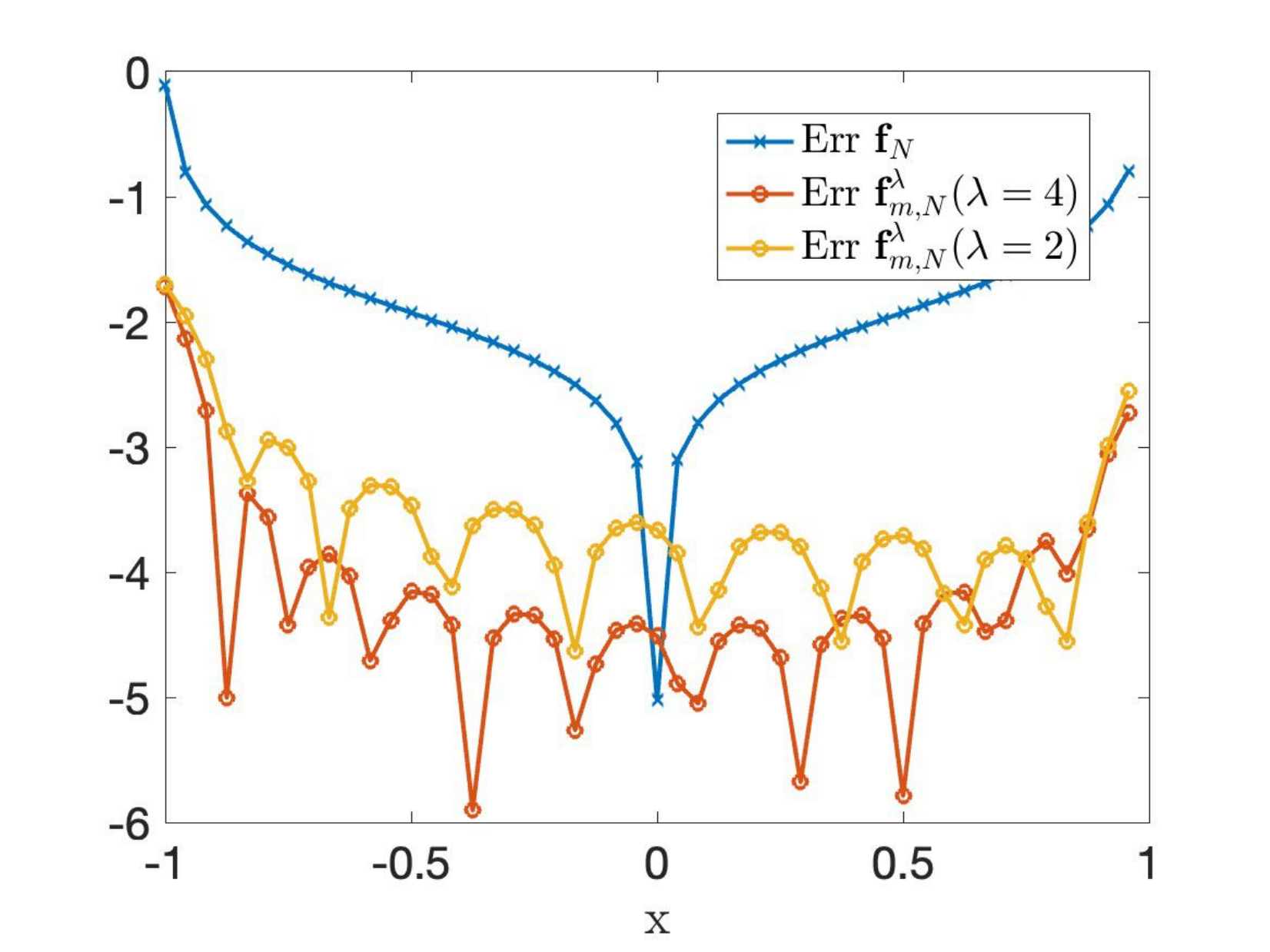}
\\
\includegraphics[width=0.24\textwidth]{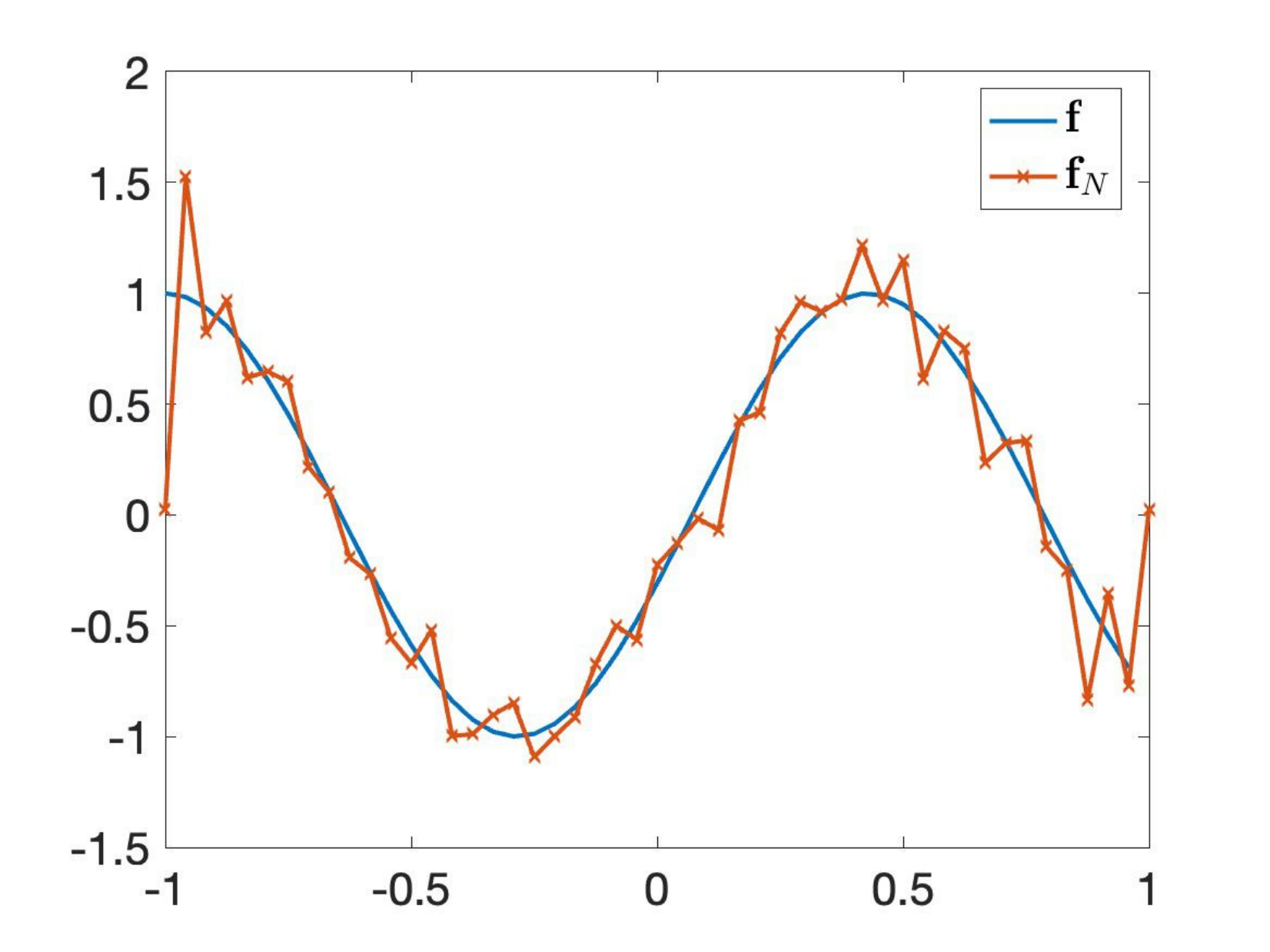}
\includegraphics[width=0.24\textwidth]{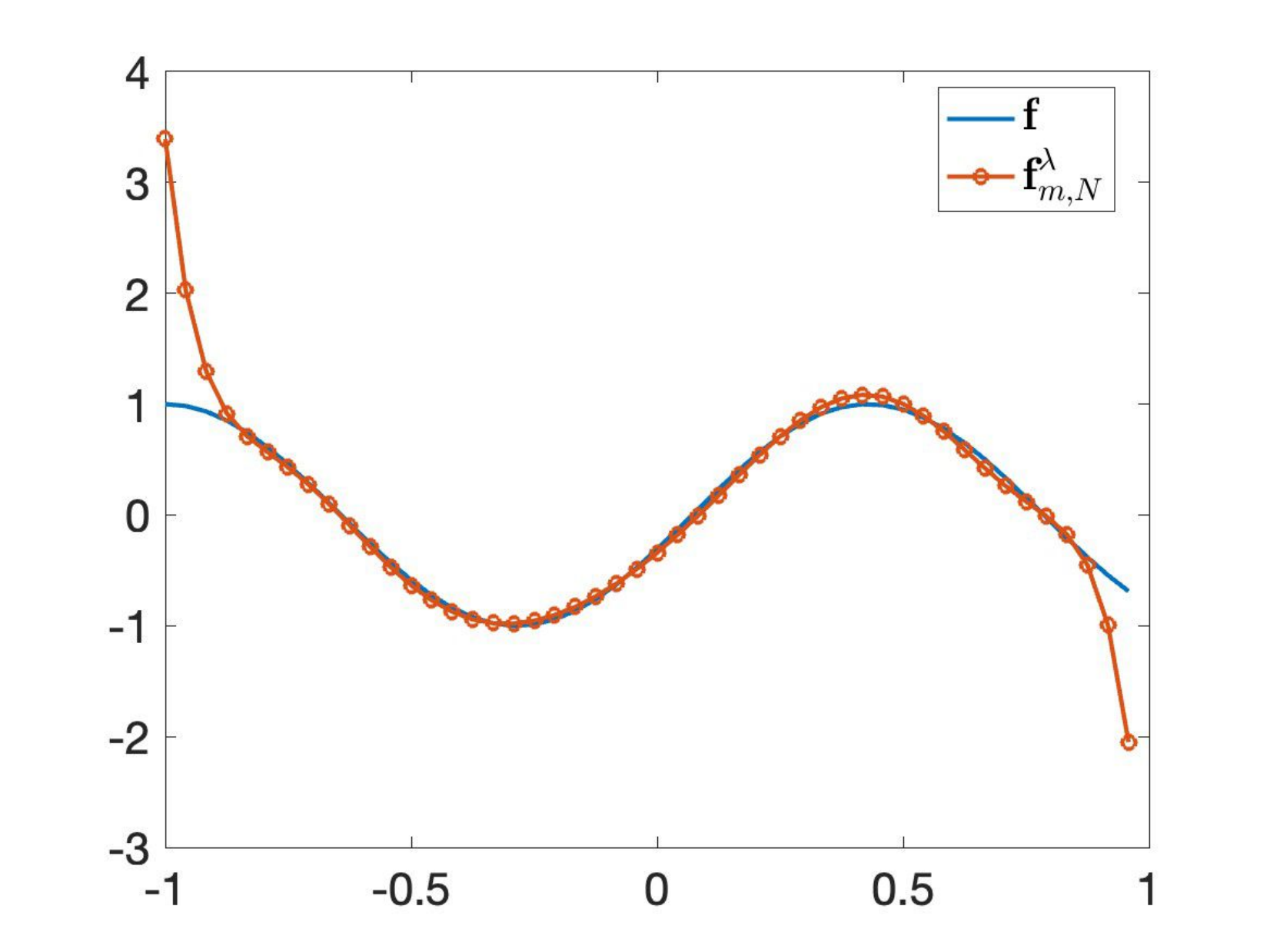}
\includegraphics[width=0.24\textwidth]
{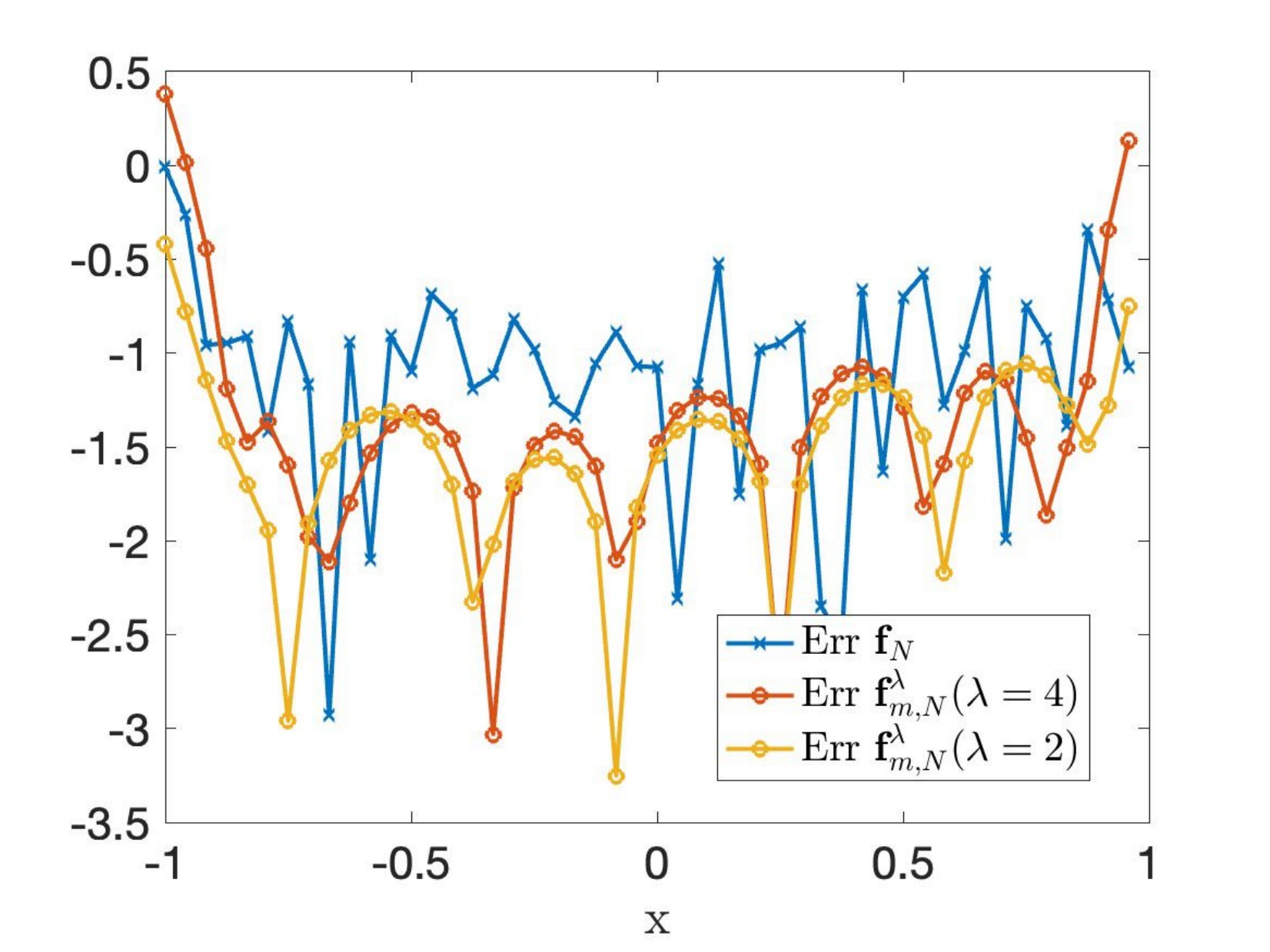}
\caption{(Top-left) Fourier reconstruction $\f_N$ in \cref{eq:partial_four_sum1 mv2}; (top-middle) Spectral reprojection $\f_{m,N}^{\lambda}$ computed via \cref{alg:gegapprox}; (top right) pointwise errors for {\em noiseless} data for $\lambda=4$ and $\lambda=2$. (Bottom) Same experiments with SNR $= 10$.}
\label{fig: 2}
\end{figure}

To serve as baseline comparisons for the GBSR method, Figure \ref{fig: 2}(top) shows the Fourier partial sum in \cref{eq:partial_four_sum1 mv2} and the Gegenbauer reconstruction obtained from \cref{alg:gegapprox} given  noiseless Fourier data. The pointwise errors for $\lambda=4$ and $\lambda=2$ are also shown. Clearly, Gegenbauer reconstruction successfully resolves the Gibbs phenomenon for appropriate choices of $\lambda$ and $m$ according to \cref{thm:GSSV92}, with $\lambda = 4$ providing a more accurate approximation than $\lambda = 2$. Figure \ref{fig: 2}(bottom) shows the same results for SNR $= 10$. The performance of Gegenbauer reconstruction notably deteriorates, especially near the boundaries, which is attributable to the Runge phenomenon  being exacerbated by noise. In this case we observe that  $\lambda=2$ yields better overall performance, and noticeably at the boundaries.  \Cref{sec:numparams} provides a more in depth  discussion regarding the relationship between SNR and $\lambda$.

\begin{figure}[h!]
\centering
\includegraphics[width=0.24\textwidth]{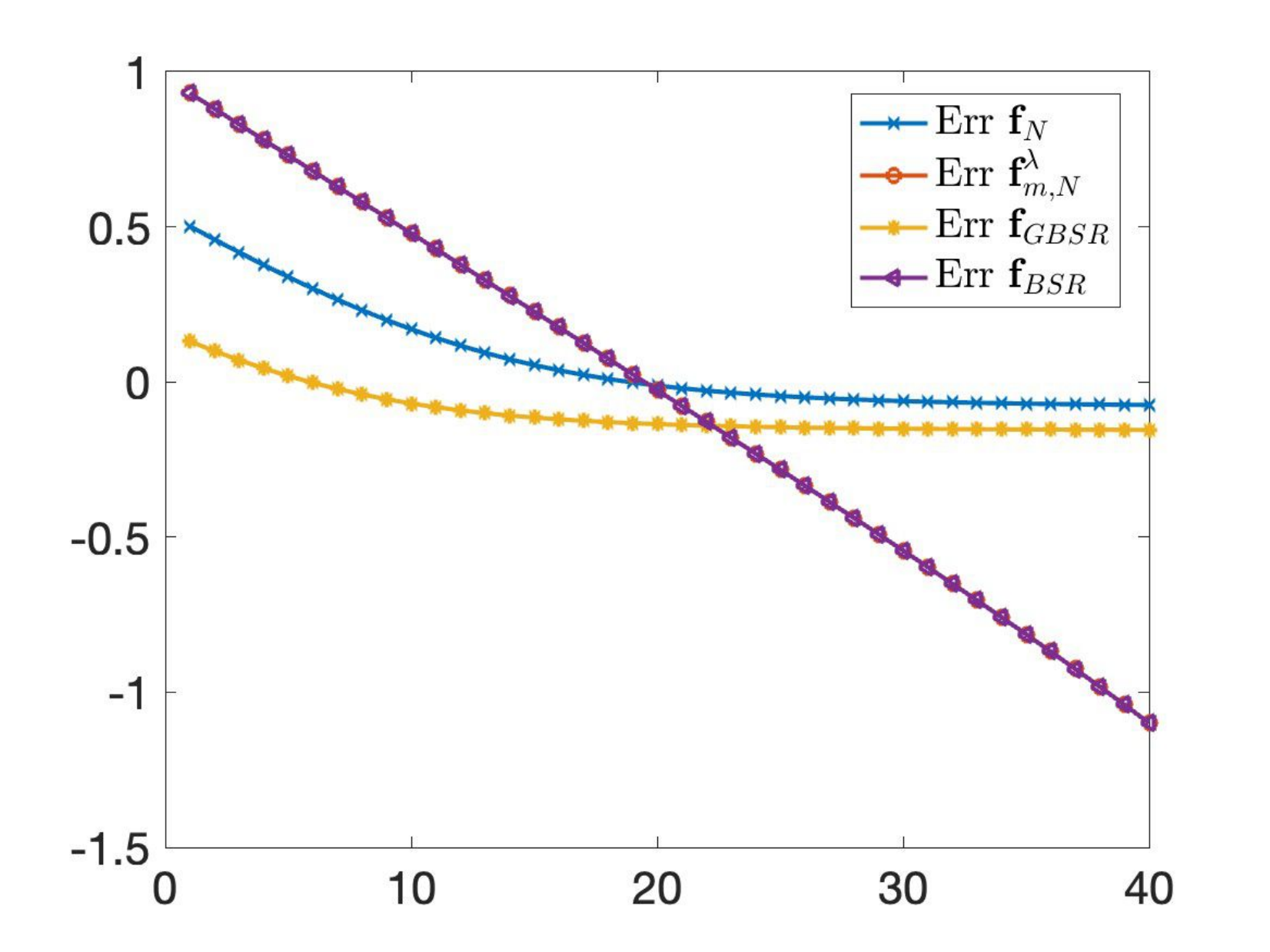}
\includegraphics[width=0.24\textwidth]{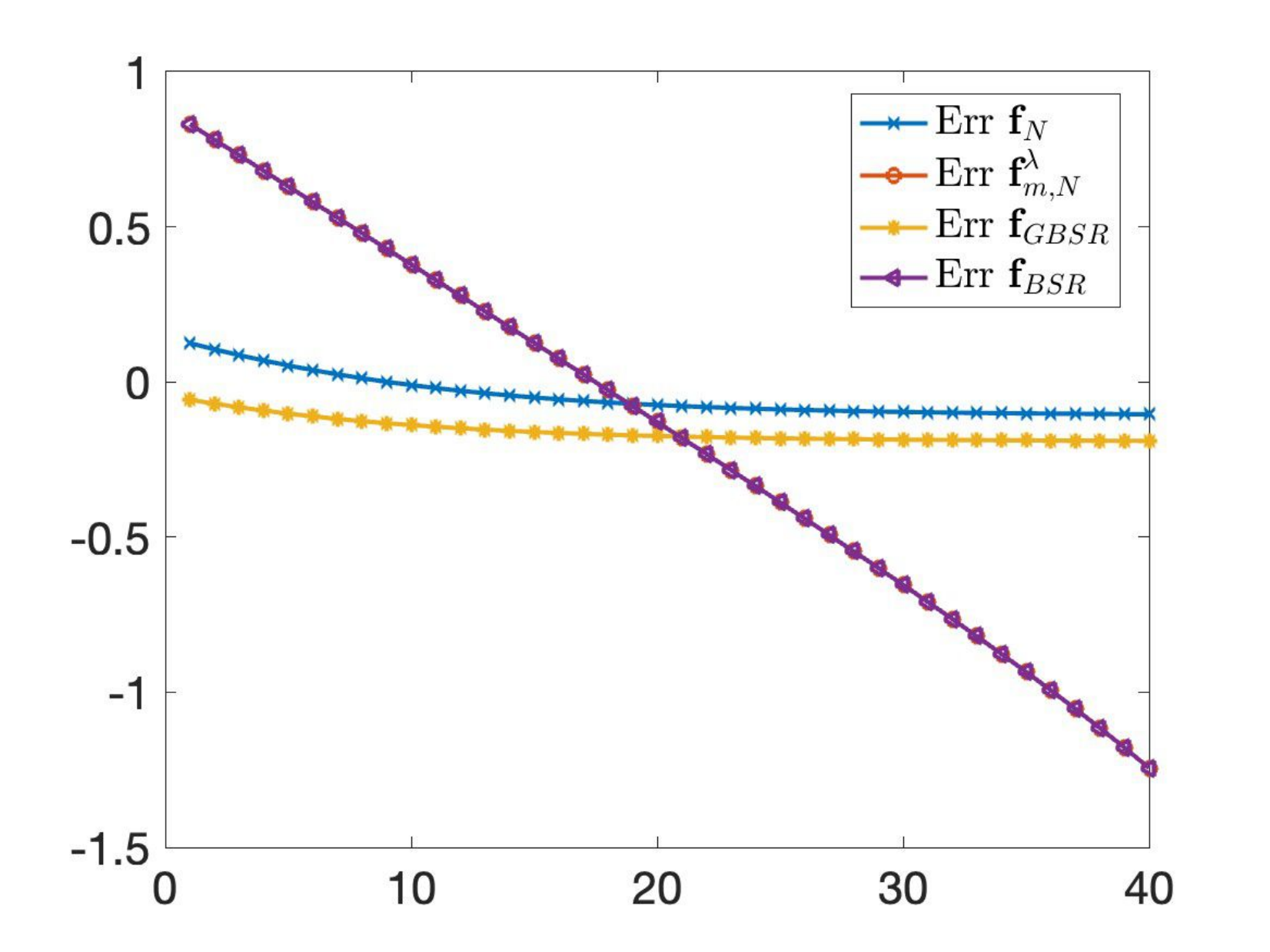}
\includegraphics[width=0.24\textwidth]{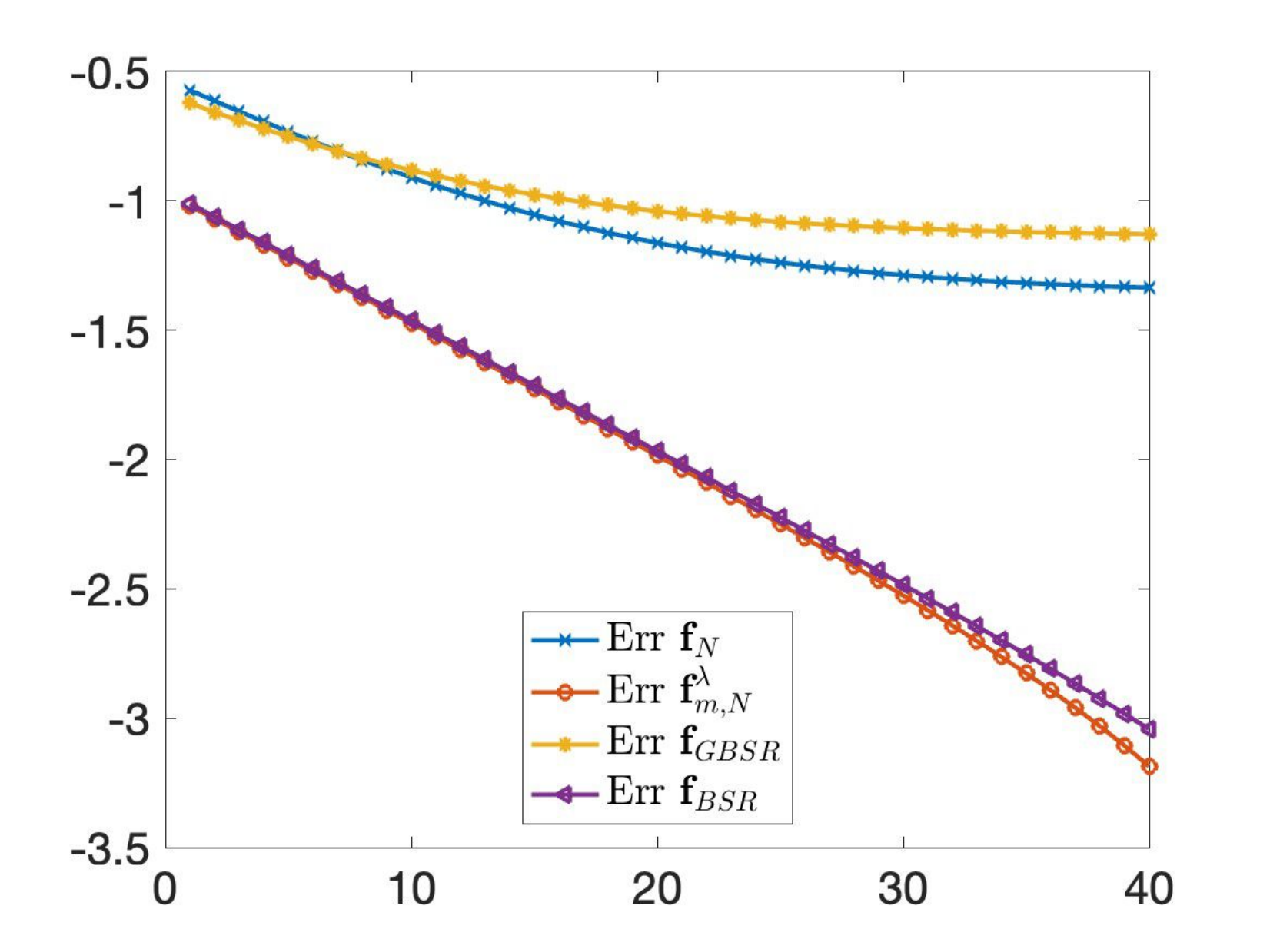}
\includegraphics[width=0.24\textwidth]{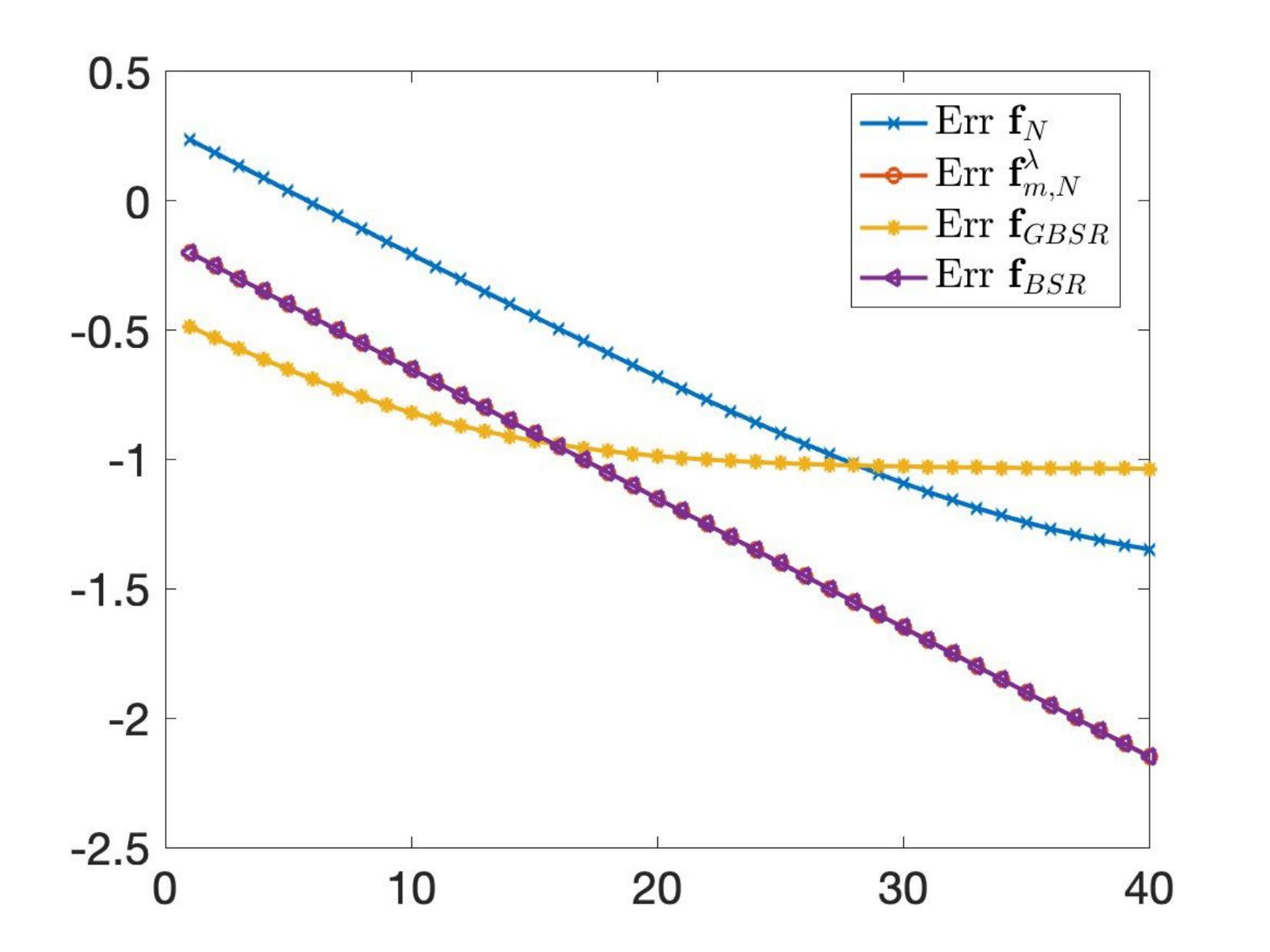} 
\caption{Log error plots for varying SNR: (left) $l_2$ error in $[-1,1]$, (middle-left) at $x=-1$; (middle-right) at $x=-0.8$; and (right) $l_2$ error in $[-.5,.5]$ for $\f_N$, $\f_{m,N}^\lambda$ $\f_{BSR}$, and $\f_{GBSR}$.}
\label{fig: 3}
\end{figure}

Figure \ref{fig: 3} compares errors in logarithmic scale as a function of SNR throughout the interval $[-1,1]$, at the boundary $x=-1$, at the point $x=-0.8$, and in the interior region from $x=-0.5$ to $x=-0.5$ (chosen to show overall error away from the boundaries). %Since the point estimate for the BSR method in \cref{alg:bayesapprox} is essentially equivalent to the original Gegenbauer spectral reprojection $\f_{m,N}^\lambda$ in \cref{alg:gegapprox}, we do not include these results here.  
We observe that with increasing SNR, the Gegenbauer approximation $\f_{m,N}^\lambda$ in \cref{alg:gegapprox} starts to converge, as is consistent with \cref{fig: 2}(right). The point estimate for the BSR method in \cref{alg:bayesapprox} is essentially equivalent to the original Gegenbauer reconstruction method. The Fourier partial sum and the GBSR method also produce more accurate approximations with increasing SNR, although as expected each method's convergence rate is not clearly established. The  $l_2$ error in $[-1,1]$ for all approximations is obviously dominated by the maximum error occurring  at $x = -1$. The advantage in using spectral reprojection becomes more apparent away from the boundary,  especially with increasing SNR.  At  $x = -.8$ the GBSR does not yield improvement from standard Fourier reconstruction,  even in low SNR environments.  By contrast, once we are far enough from the boundary, using \cref{alg:bayesapprox2} yield  more accurate solutions for low SNR. To summarize, in the low SNR regime, the Gegenbauer reconstruction method does not satisfy the conditions for obtaining numerical convergence, especially near the boundary. In this case, the GBSR method provides a relatively reasonable approximation, as well as uncertainty quantification. % In this case we speculate that error contributions from the noise are likely more than those coming from the trapezoidal rule, so that the choice of $\lambda$ does not negatively impact the overall performance.  These results also suggest that while the Gegenbauer reprojection method is {\em not robust} with respect to $\lambda$ or noise, the GBSR approach generally is.   

%%%%%%%%%%%%%%%%%%%%%%%%%

\subsection{Relationship between SNR and Gegenbauer parameters}
\label{sec:numparams}

\Cref{thm:GSSV92} establishes the relationship between the number of Fourier coefficients $N$, the expansion degree in the Gegenbauer polynomial basis $m$ and the corresponding weighting parameter $\lambda$. It does not discuss the relationship between these parameters and SNR, however. By design, the BSR method yields MAP estimates that behave similarly to the Gegenbauer reconstruction solution.  Specifically,  parameters that satisfy \Cref{thm:GSSV92} yield the best approximations, with $\lambda = \kappa N$ required for small error $E_D$ being analogous to the likelihood term \cref{eq:likelihooddistribution1}. On the other hand, $\lambda = \kappa N$ is not required for the GBSR method according to its likelihood distribution \cref{eq:likelihooddistribution2}.  For {\em both} methods, when the observables are practically noise-free, there is some non-negligible contribution in using the trapezoidal rule \cref{eq: trapezoid} for small $\lambda$, especially for small $N$.  In all cases, $m$ must be chosen large enough to resolve the underlying function.  In this regard, it is important to point out that the parameter $m$ is what describes our {\em prior} belief. That is, for $m = 0$ we assume the function is constant, for $m = 1$, the function is linear, etc.  Indeed, $m = 0$ and $m = 1$ correspond more closely to what is often assumed for image recovery.  By contrast, the BSR and GBSR enable the recovery of functions with much greater variability (as dictated by $m$).

While the choice of $m$ is somewhat independent of SNR, how $\lambda$ might be chosen changes significantly as SNR decreases.  First, any error made using trapezoidal rule is insignificant in low SNR environments.  Second, as already discussed, the variance of the Gegenbauer reconstruction method increases at the boundaries.  Coupled with the ill-conditioning effects for larger choices of $\lambda$, it is not clear how the value of $\lambda$ affects the performance of each method. 

To emphasize the effects of noise on the choice of $\lambda$, we consider the same example in \cref{sec:SNRcase} for varying (integer) values of  $\lambda$ and fixed SNR of $2, 10,$ and $30$. %Particularly, we still consider $f(x)=\cos(1.4\pi(x+1))$ on  $[-1,1]$ and 
As before we are given $N = 48$ noisy Fourier coefficients $\wh{\bb}$ in \eqref{eq:noisyfourier1d} and we again  choose $m = 9$ for $\lambda \in [1,8]$.

\begin{figure}[h!]
\centering
\includegraphics[width=0.24\textwidth]{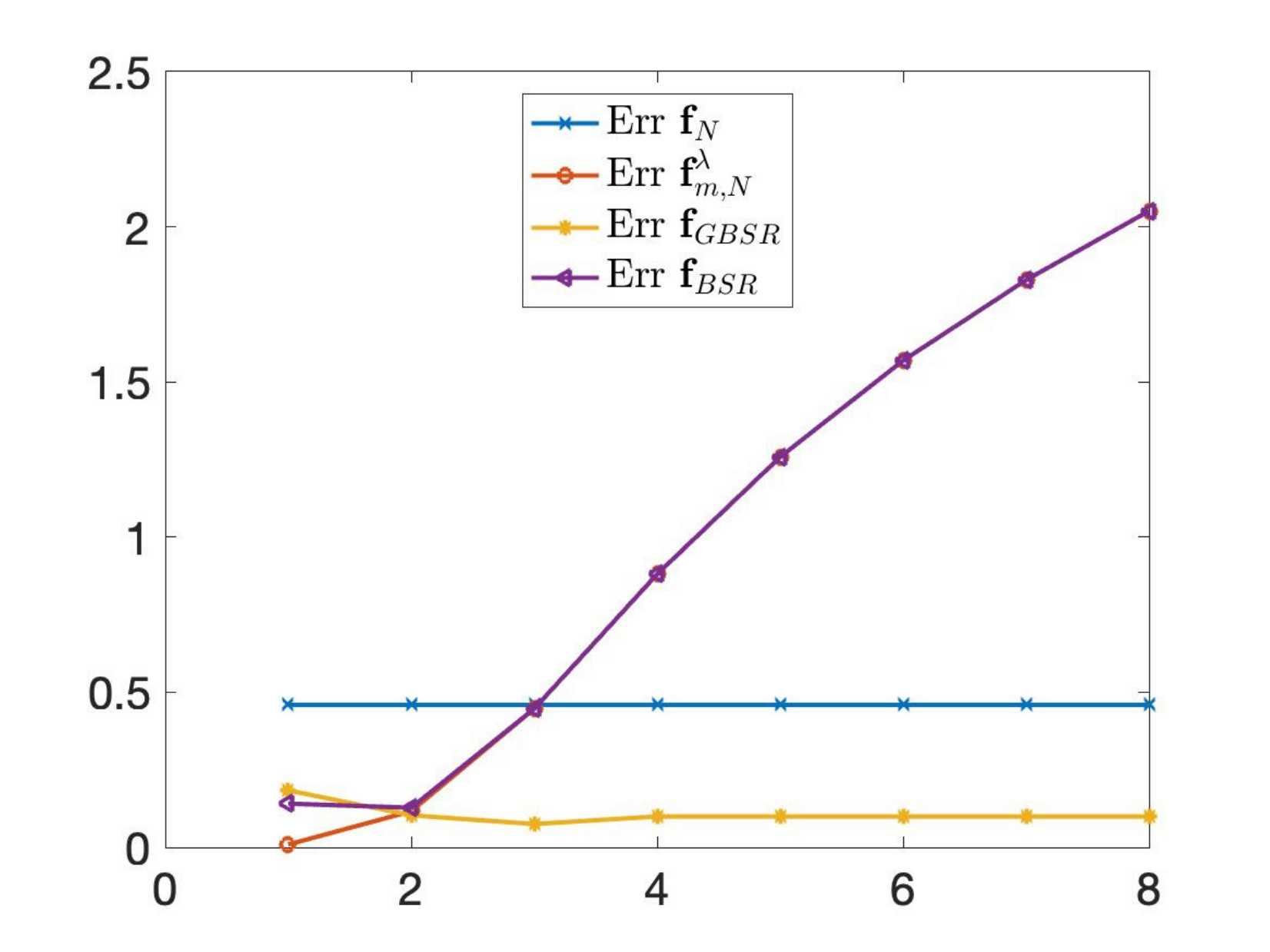}
\includegraphics[width=0.24\textwidth]{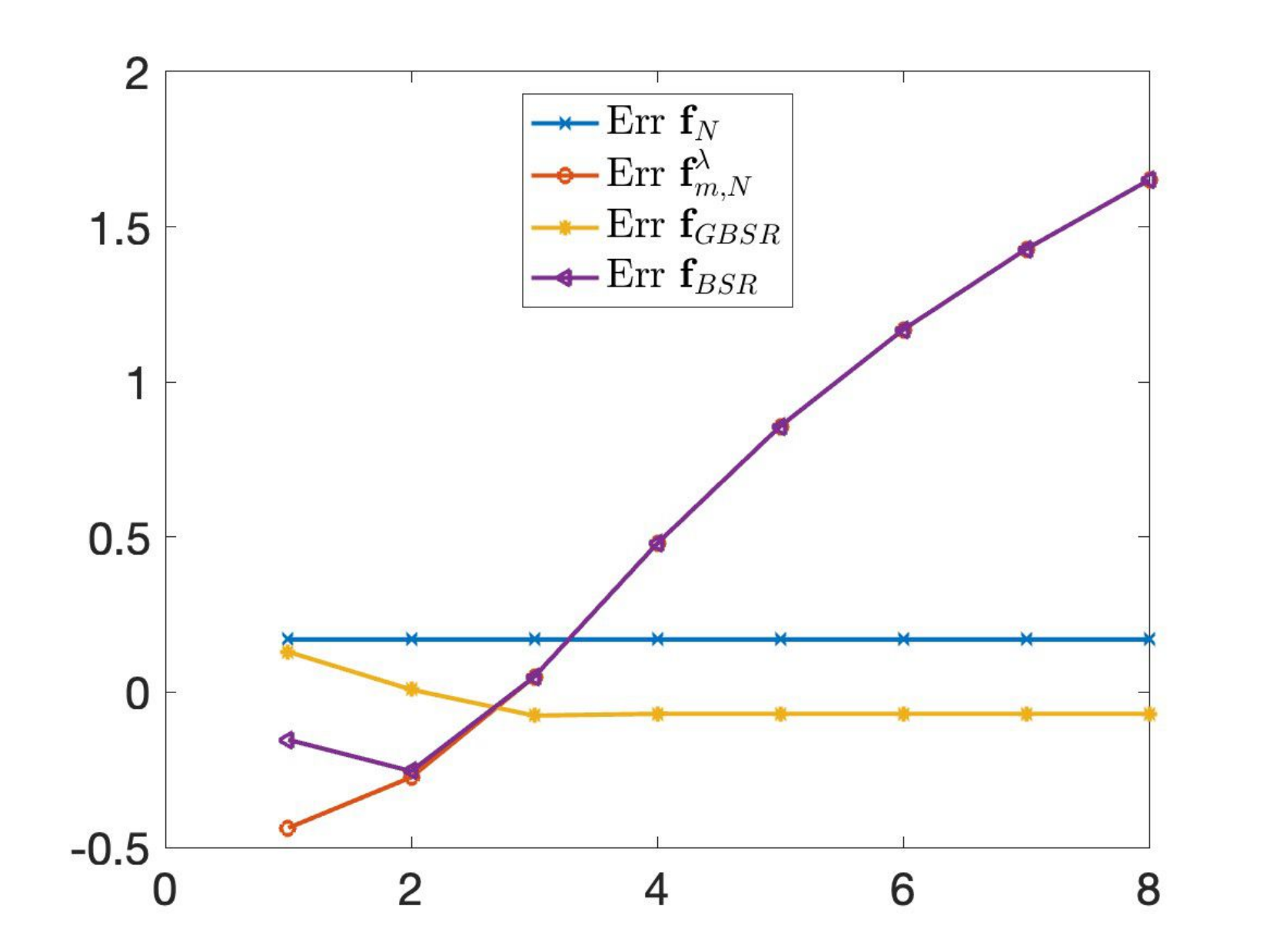}
\includegraphics[width=0.24\textwidth]{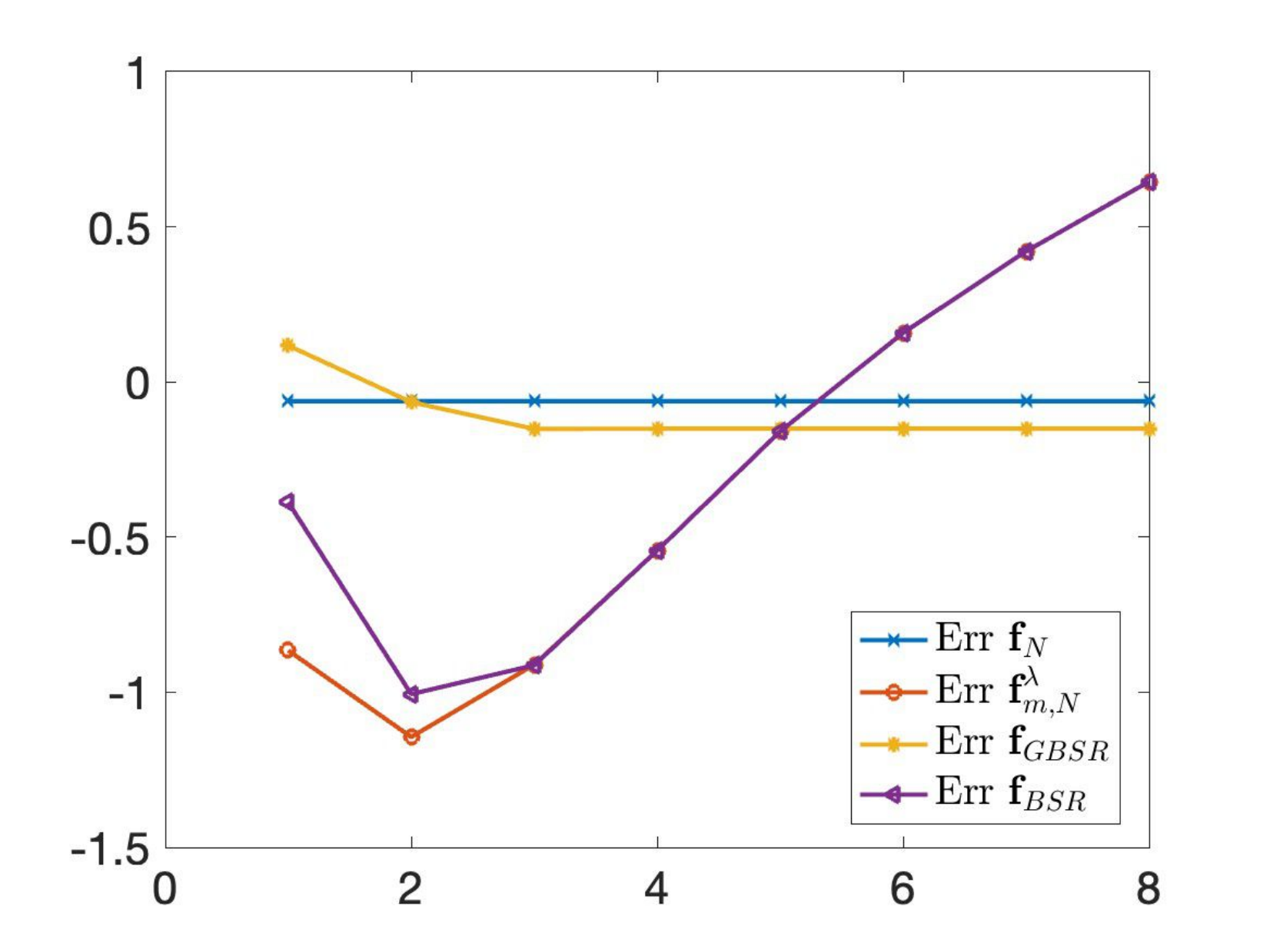}
\caption{Log error plots for varying $\lambda$: $l_2$ error in $[-1,1]$ with SNR value 2 (left), 10 (middle) and 30 (left) for $\f_N$, $\f_{m,N}^\lambda$ $\f_{BSR}$, and $\f_{GBSR}$.}
\label{fig: 4}
\end{figure}

Figure \ref{fig: 4} compares the $l_2$ error in $[-1,1]$ in logarithmic scale with respect to different $\lambda$ for SNR $= 2, 10$ and $30$. No theoretical justification is available  for applying any of these methods in low SNR environments. As  SNR increases, it becomes more evident that an optimal  $\lambda$ for both the Gegenbauer reconstruction  and BSR methods exists, and is a cumulative result of all the possible considerations already mentioned. On the other hand, the GSBR method is clearly robust for various $\lambda$ values in low SNR environments, mainly because the likelihood term does not involve its use. Of course, none of these methods are a panacea in low SNR environments.  The GBSR and BSR methods provide additional uncertainty quantification for fixed hyperparameters, although by design the BSR credible interval is quite narrow. %\help{[Anne: Tongtong,  I  think it would be useful to provide CI graphs corresponding to some of the choices in Figure 6.4, perhaps all three choices of SNR with $\lambda = 4$ for BSR and GSBR (to drive home the point that BSR CI is not quite as useful)?]}
We illustrate this in Figure \ref{fig: 5}, where the BSR and GBSR approximations are provided along with their respective credible intervals for $\lambda = 2, 4$. Consistent with \cref{fig: 4}, the BSR provides the most accurate results for $\lambda = 2$. The credible interval is very narrow, even in low SNR environments and almost impossible to discern for SNR = $10$ and $30$. By contrast, we observe nearly identical estimations for $\lambda = 2, 4$ when using GBSR.  There is furthermore additional information to be gleaned from their credible intervals.  Moreover the error near the boundaries appears to be dominated by the likelihood term also having a larger error there, attributable to the over/undershoots caused by the Gibbs phenomenon, as opposed to simply noise (see \cref{fig: 2}(left)). % fact that GBSR weights likelihood more heavily particularly near the boundary deteriorates its overall $l_2$ error, as the error is dominated by that near the boundaries. }

\begin{figure}[h!]
\centering
\includegraphics[width=0.24\textwidth]{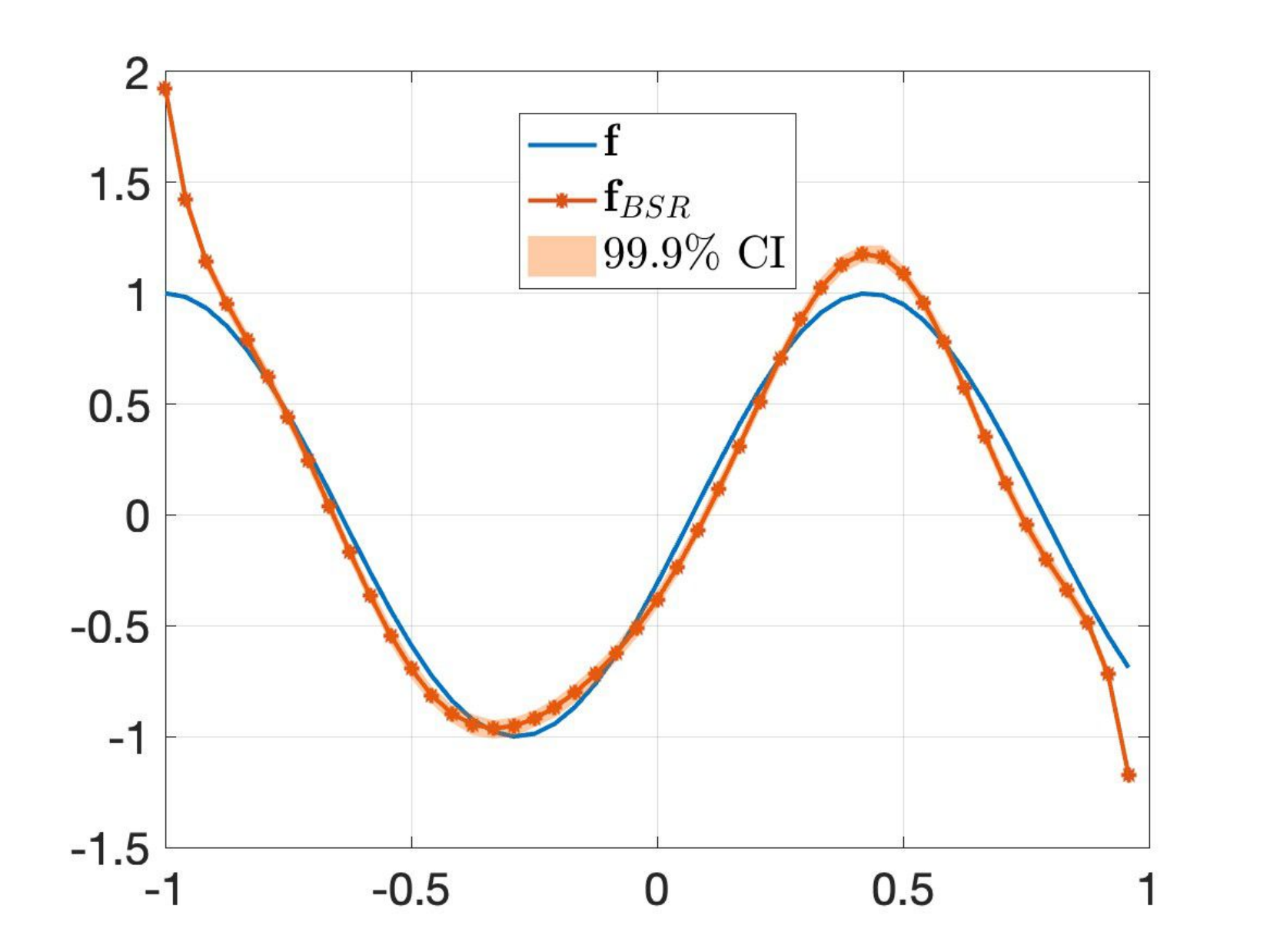}
\includegraphics[width=0.24\textwidth]{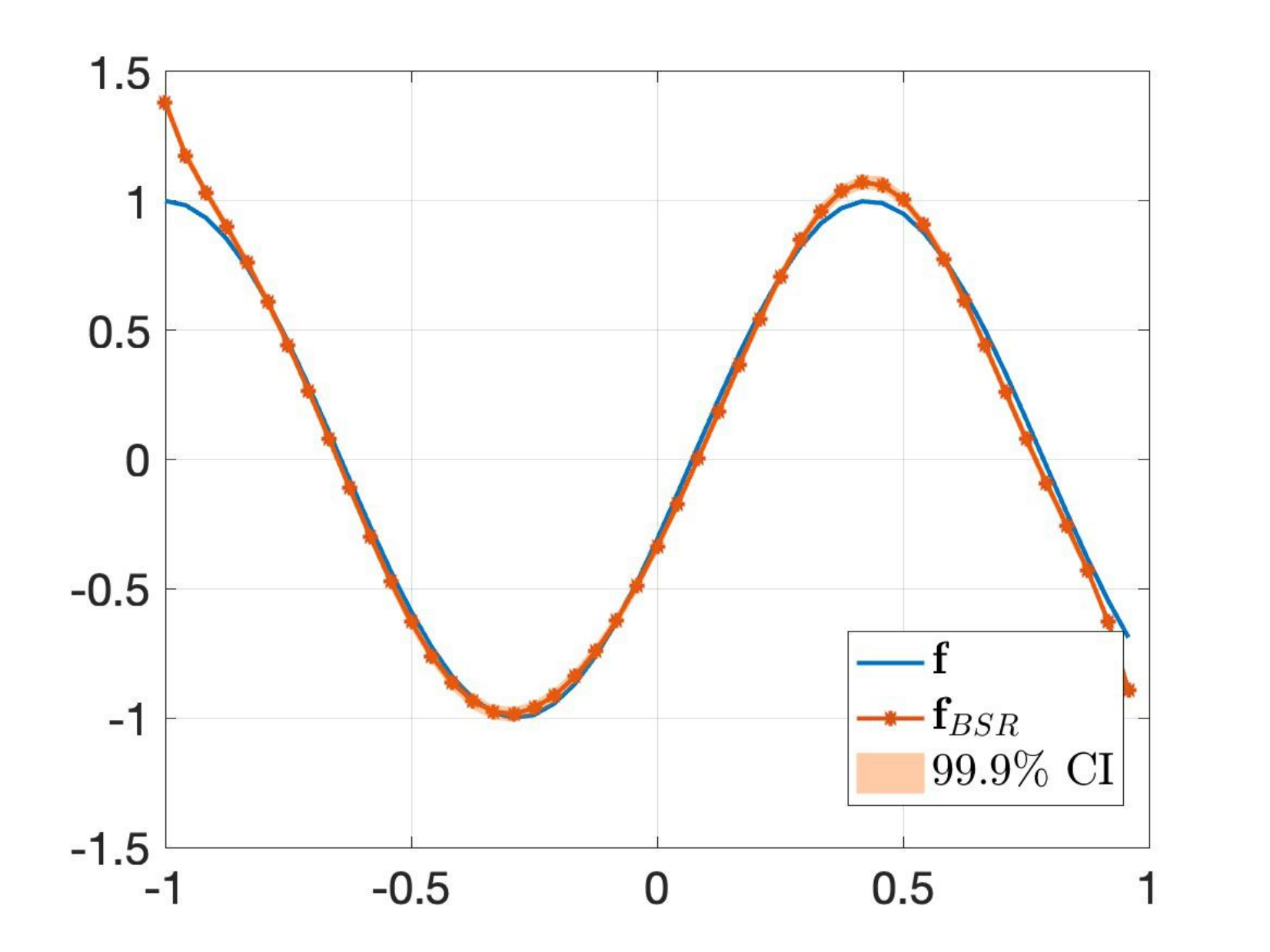}
\includegraphics[width=0.24\textwidth]{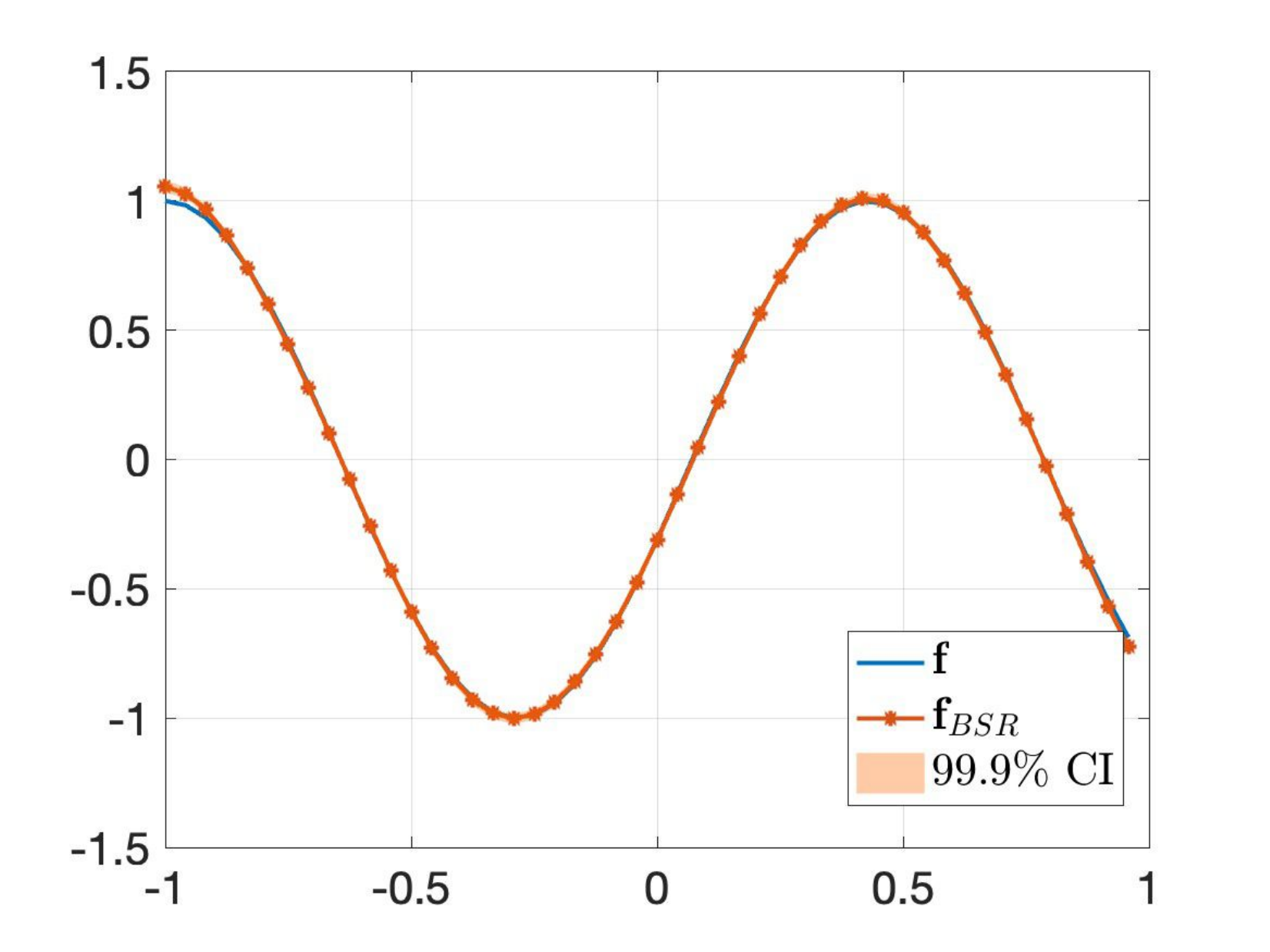} \\
\includegraphics[width=0.24\textwidth]{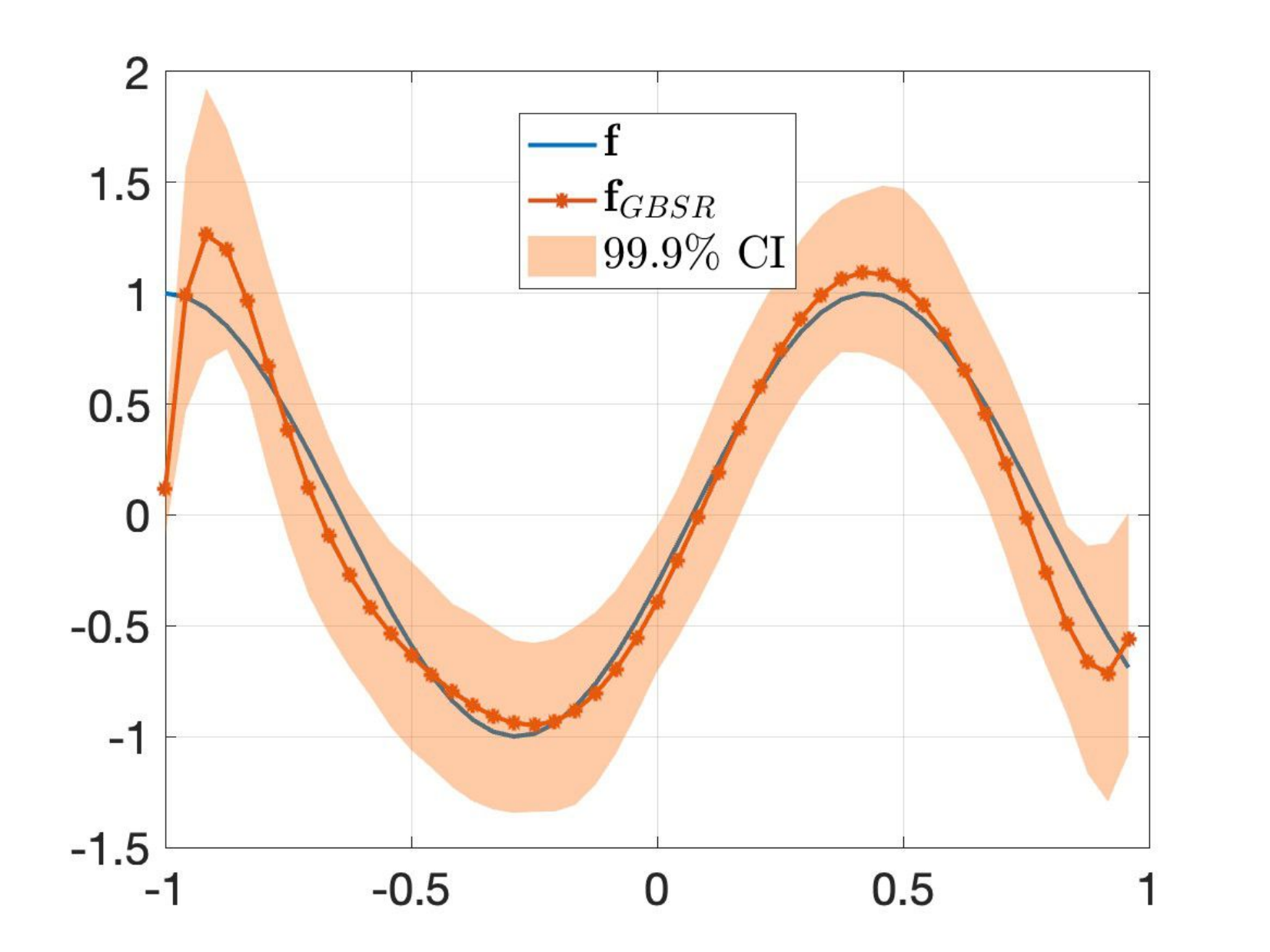}
\includegraphics[width=0.24\textwidth]{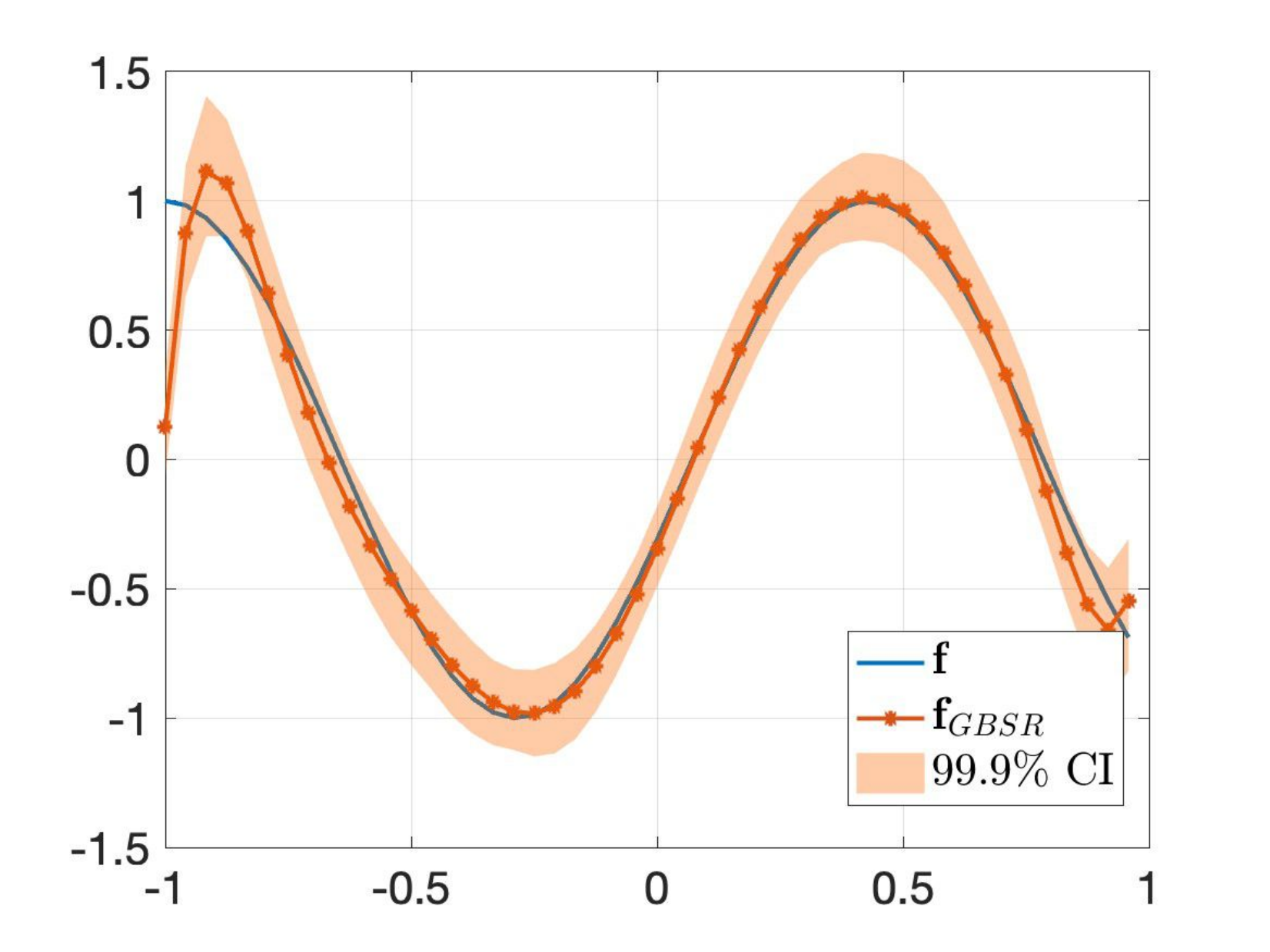}
\includegraphics[width=0.24\textwidth]{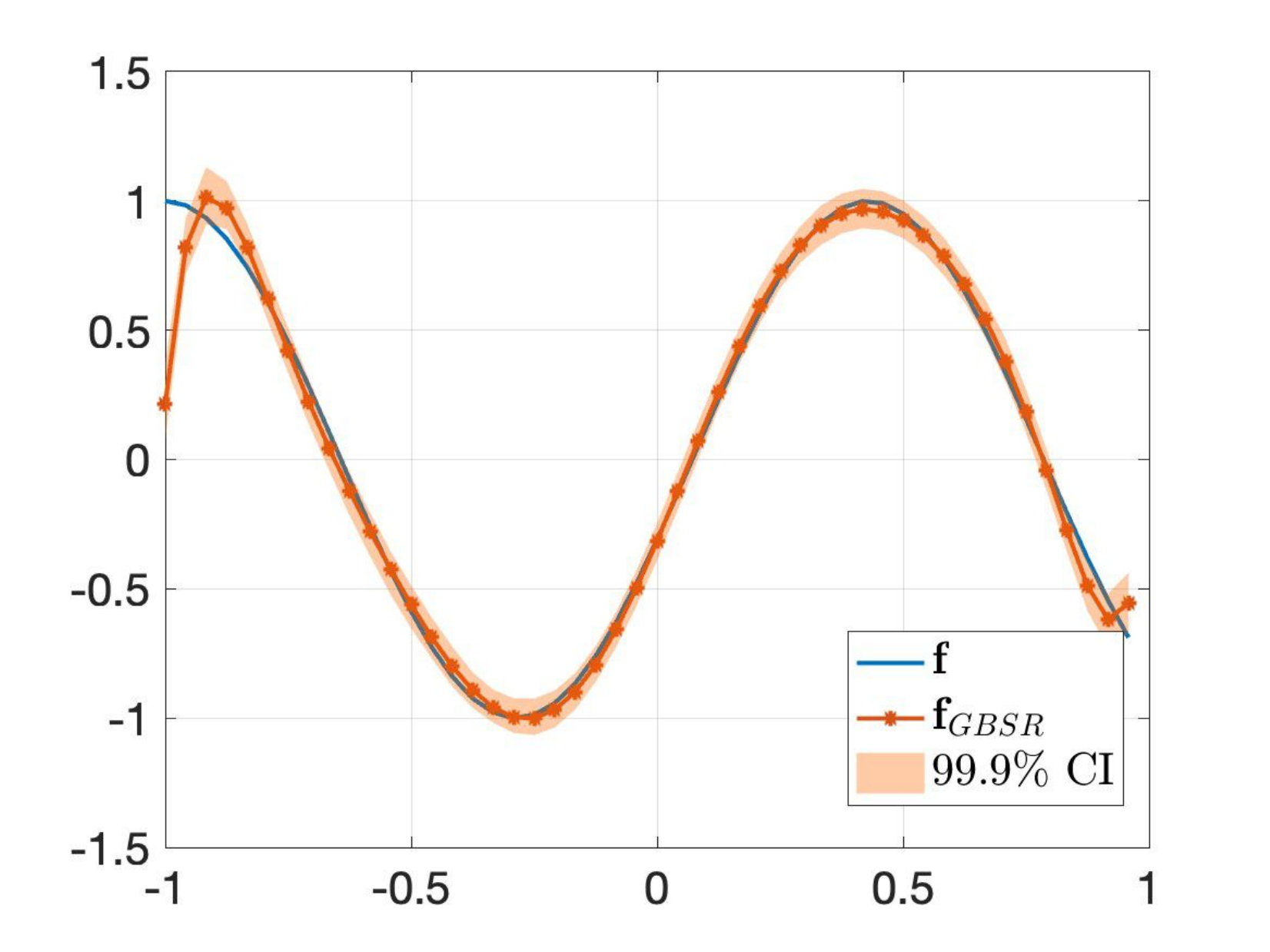}\\
\includegraphics[width=0.24\textwidth]{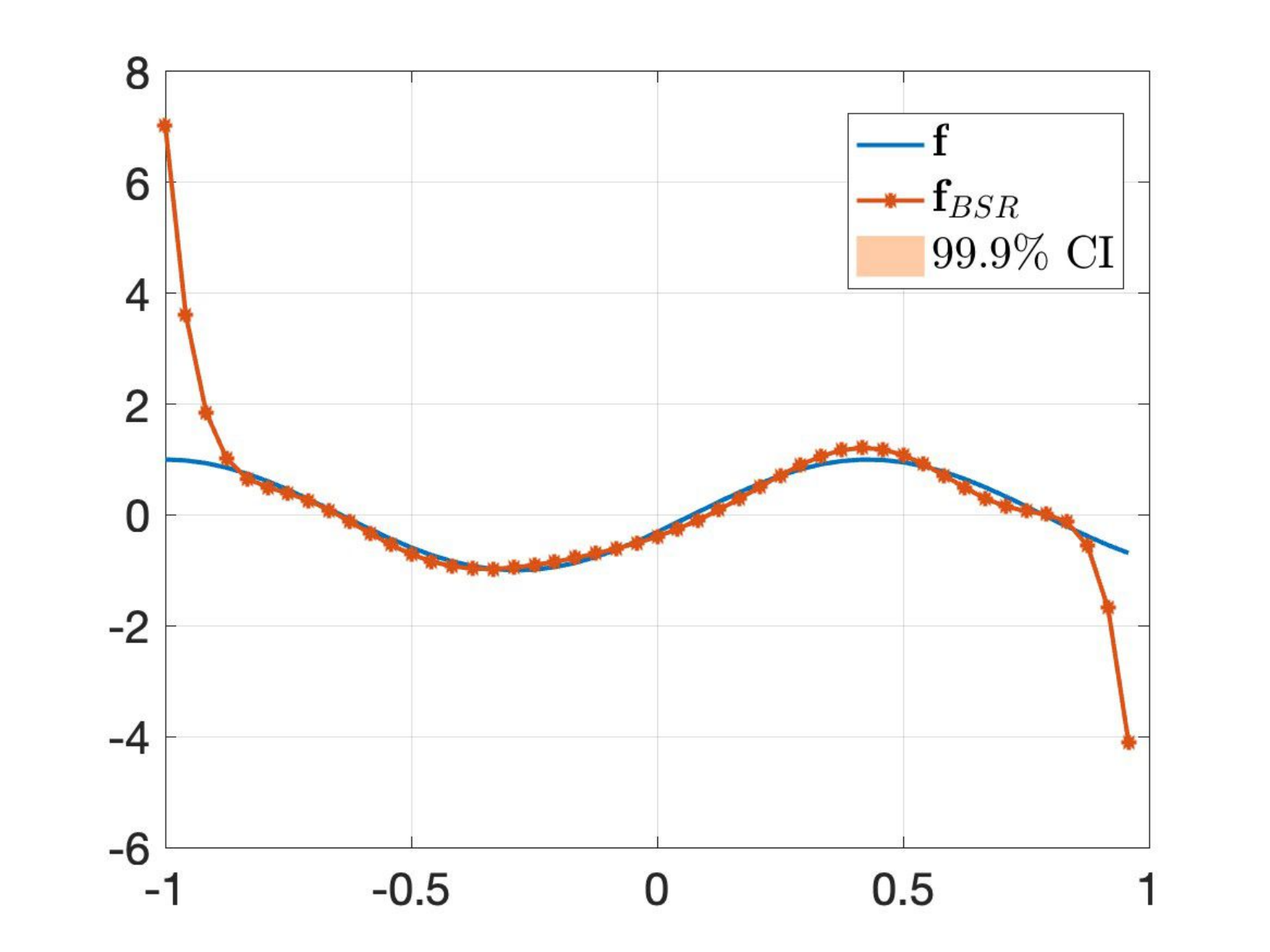}
\includegraphics[width=0.24\textwidth]{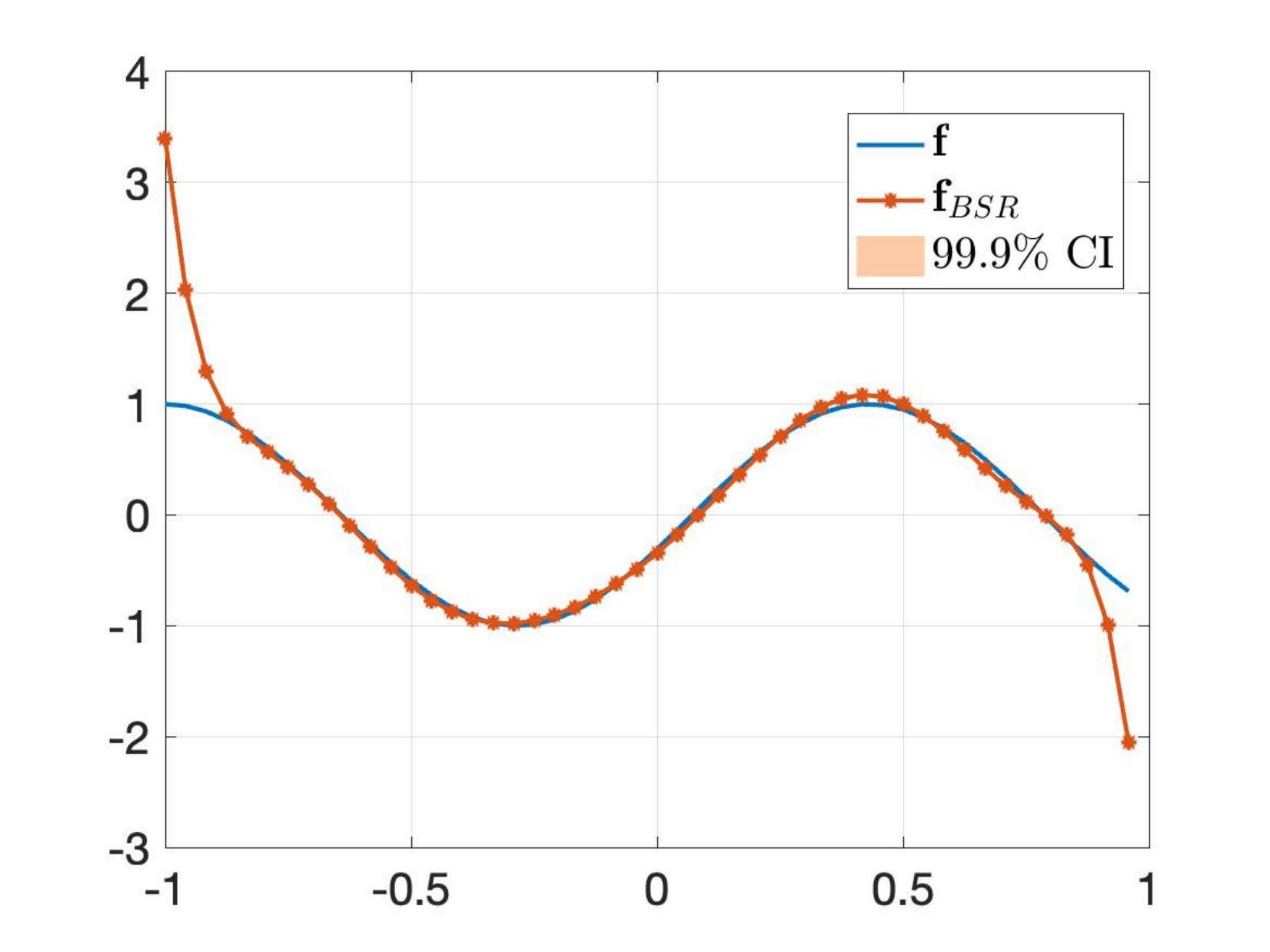}
\includegraphics[width=0.24\textwidth]{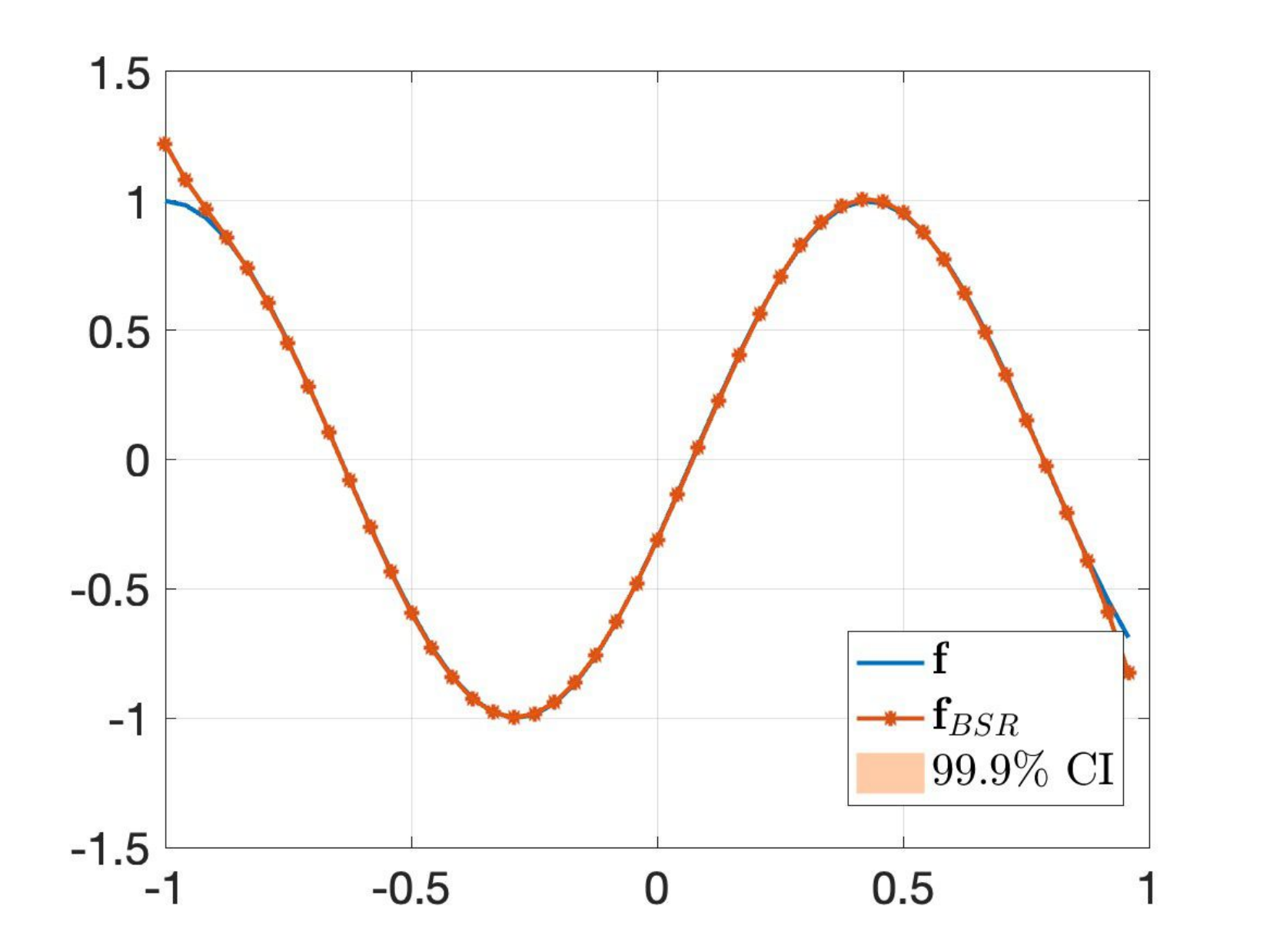} \\
\includegraphics[width=0.24\textwidth]{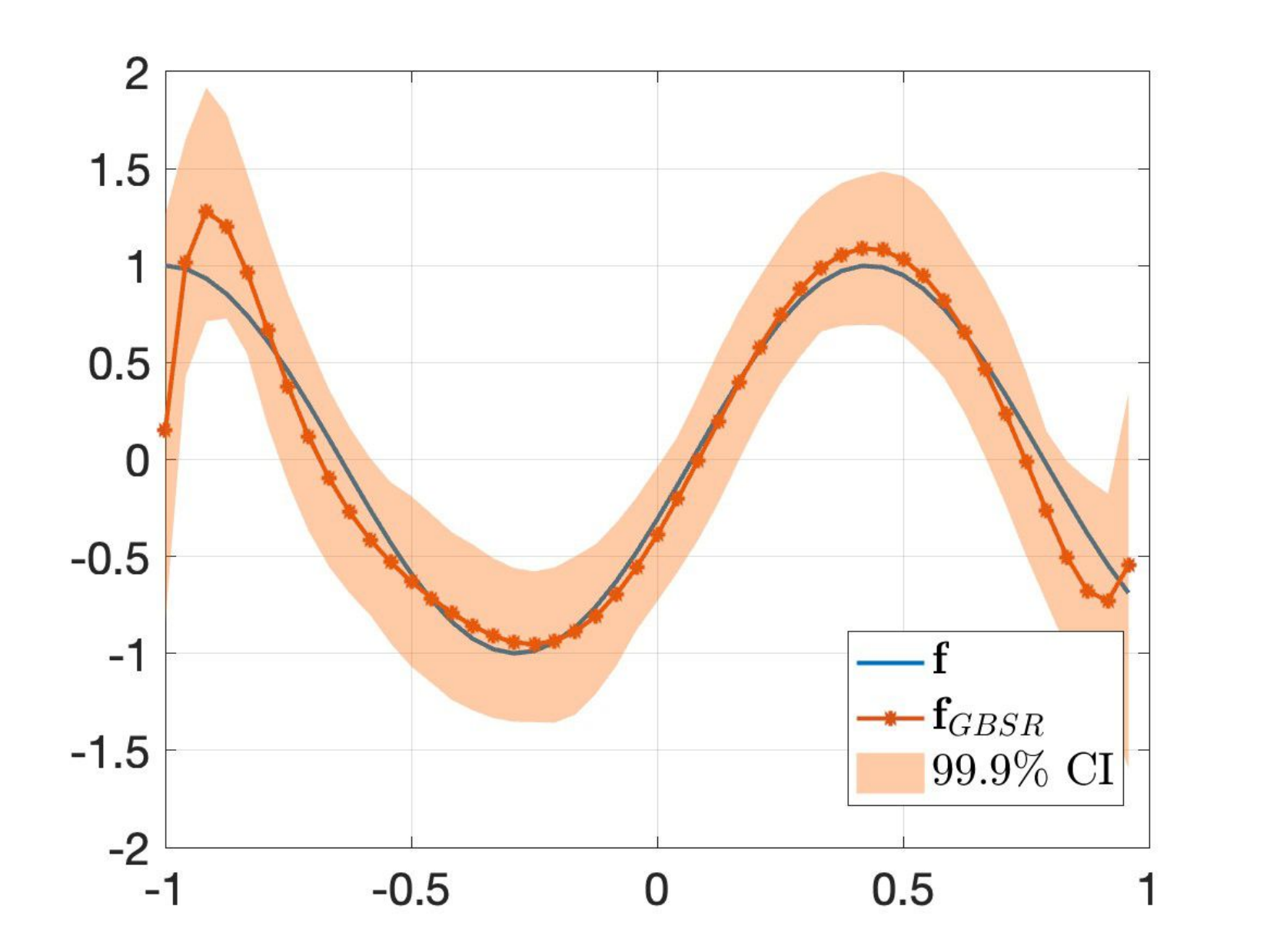}
\includegraphics[width=0.24\textwidth]{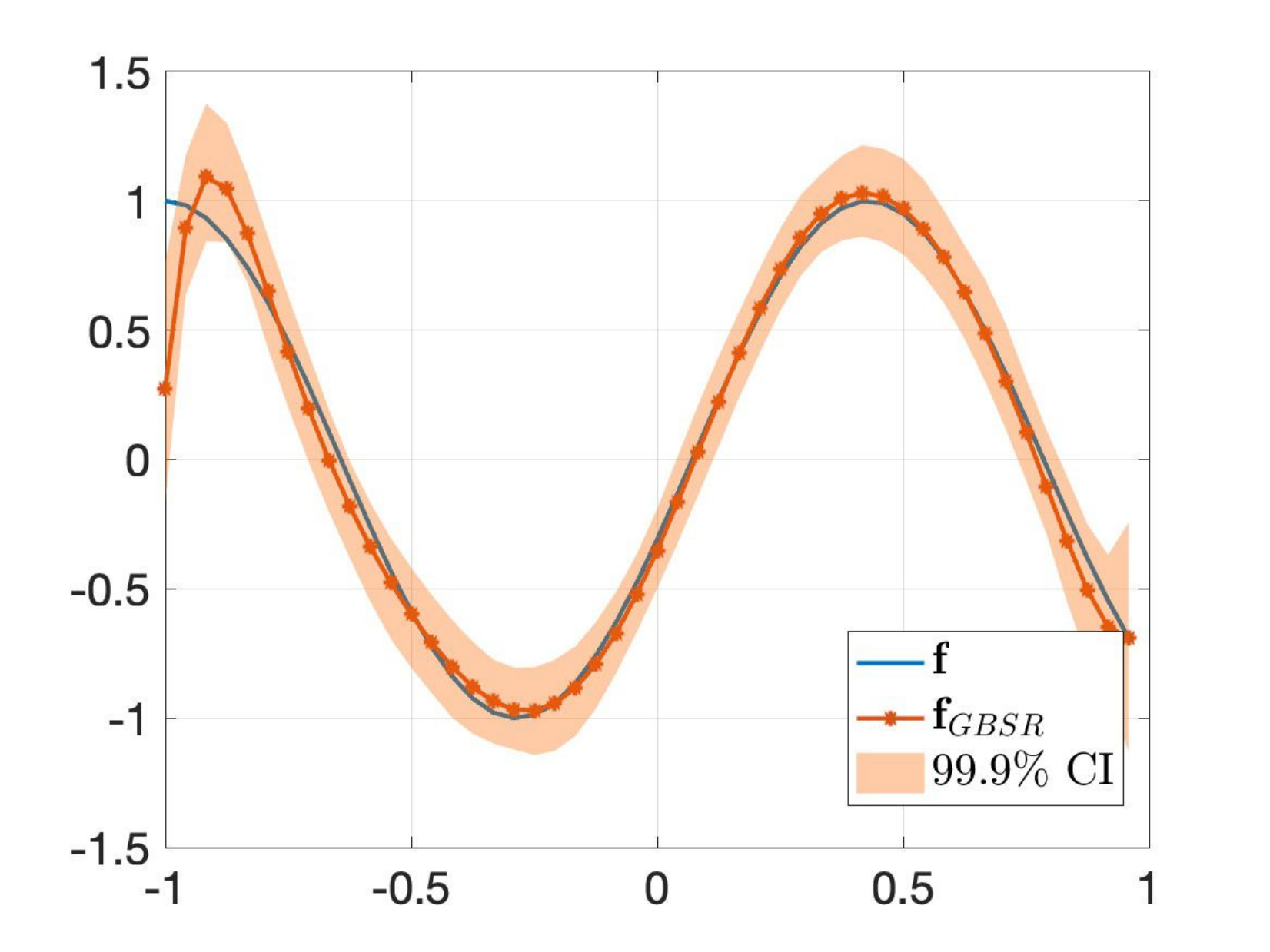}
\includegraphics[width=0.24\textwidth]{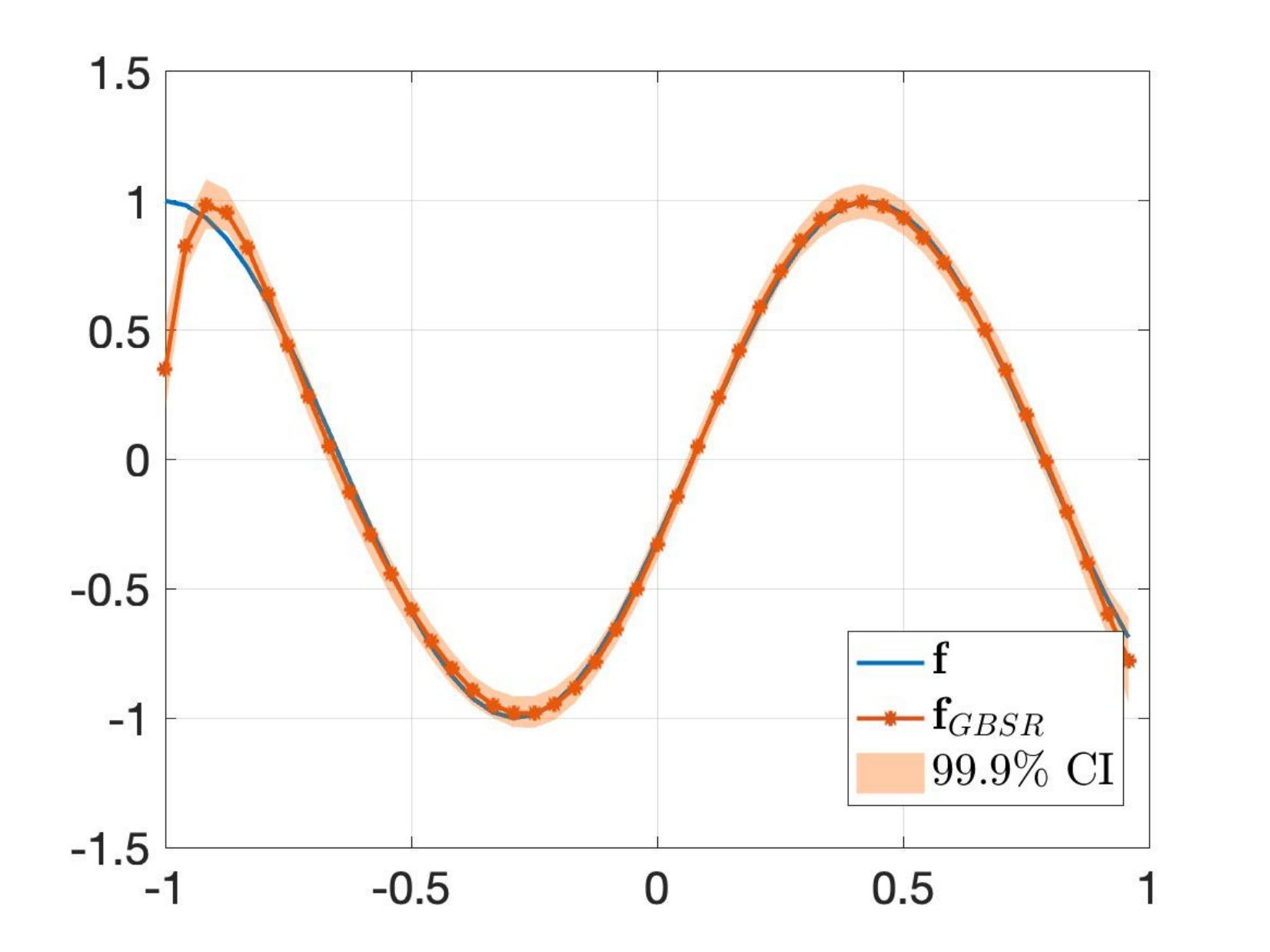}
\caption{Credible intervals with $\f_{BSR}$ and $\f_{GBSR}$ for SNR value 2 (left), 10 (middle) and 30 (left) and $\lambda=2$ (top and middle-top), $\lambda=4$ (middle-bottom and bottom).}
\label{fig: 5}
\end{figure}

We emphasize that our numerical experiments are not meant to be prescriptive -- that is, it is not the case that one method always outperforms another in all environments.  This is to be expected, since \cref{thm:GSSV92} has already identified the best solution in the noiseless case.  Rather we demonstrate how spectral reprojection  inspires two Bayesian approaches, BSR and GSBR, both of which can recover analytic non-periodic functions from noisy Fourier data and provide some uncertainty quantification. Indeed, our results suggest that there is justification for choosing a hybrid BSR/GSBR approach, one that weighs the BSR more heavily near the boundaries and GSBR more heavily in the interior.  It is also possible to modify the BSR method so that the likelihood and prior use different parameters for $m$ and $\lambda$, which also can be analyzed in the context of the proof of \cref{thm:GSSV92}.  These and other ideas  will be explored in future investigations. 
\section{Concluding Remarks}
\label{sec:summary}

In this investigation we proposed the Bayesian spectral reprojection (BSR) and generalized  Bayesian spectral reprojection (GBSR) inference methods to recover smooth  non-periodic functions from given noisy Fourier data.  Both method are motivated by the theoretical convergence analysis provided for the  spectral reprojection method given noiseless data. These insights are brought into a statistical framework to effectively combine numerical analysis with uncertainty quantification.

The BSR approach provides a statistical framework for the spectral reprojection method, enabling a better understanding of the limitations for using spectral reprojection in noisy environments. The likelihood and prior terms follow directly from the error analysis established for the original (noiseless) spectral reprojection method in \cite{GSSV92}. By design, the BSR yields estimates  similar to the Gegenbauer reconstruction method, while also  providing uncertainty quantification for fixed hyperparameters. As both the likelihood and prior terms are related to Gegenbauer polynomials, the approach limits the space of possible solutions, resulting in very narrow credible intervals. Furthermore, the approximation quality for both spectral reprojection and BSR  rapidly deteriorate with decreasing SNR. 

The GBSR approach is then proposed to expand the space of possible solutions, which is accomplished by incorporating the Fourier data directly into the likelihood term (as opposed to first transforming the data to the space of $m$th-degree Gegenbauer polynomials). In this regard the GBSR provides a more natural construction for computational linear inverse problems. The resulting approximation is more robust  with respect to the Gegenbauer polynomial weighting parameter $\lambda$, especially in low SNR environments.  By design it cannot achieve the same accuracy for high SNR, as the error from the likelihood term is still large due to the effects of the Gibbs phenomenon. The best accuracy for the limiting case of no noise  is  already obtained using spectral reprojection. % It weights the likelihood more heavily, especially near the boundary. Thus it provides better results in the low SNR regime, when the spectral reprojection method and BSR method deteriorates, yet cannot achieve numerical convergence with increasing SNR as the other two methods.

There are clearly benefits and drawbacks in using both the BSR and GBSR, depending on SNR and Gegenbauer weight $\lambda$. Importantly, these benefits and drawbacks are highly predictable  by the construction of each method and theoretical motivation  in \cite{GSSV92} (along with the well-known properties of Fourier reconstruction for non-periodic functions).  Because of this, developing a hybrid approach that does not require additional user inputs is completely attainable.  

\section*{Acknowledgements} 
%This work is partially supported by the NSF grant DMS \#1912685, AFOSR grant \#F9550-22-1-0411, DOE ASCR \#DE-ACO5-000R22725, and DOD ONR MURI grant \#N00014-20-1-2595.

The authors would like to thank Guohui Song and Jack Friedman for their thoughtful insights during the early stages of this research. 

% \section*{Data availability}

% All data generated or analyzed during this study are included in this published article. MATLAB codes used for obtaining results are available upon request to the authors.

% \section*{Declarations}

% This work is partially supported by the NSF grant DMS \#1912685 (AG), AFOSR grant \#F9550-22-1-0411 (TL, AG), DOE ASCR \#DE-ACO5-000R22725 (AG), and DOD ONR MURI grant \#N00014-20-1-2595 (TL, AG). Anne Gelb is a member of the editorial board of the Journal of Scientific Computing. The authors have no other relevant financial or non-financial interests to disclose.

\appendix 

\small
\bibliographystyle{siamplain}
\bibliography{mybibfile}

\end{document}